\magnification=\magstep1
\input amstex
\documentstyle{amsppt}
\pagewidth{6.5truein}
\pageheight{8.9truein}
\ifx\refstyle\undefinedZQA\else\refstyle{C}\fi
%
%
%
%
%
%
%
\count255=\catcode`\!
\catcode`\!=11
\ifx\plot!loaded\relax
   \catcode`\!=\count255\endinput\else\let\plot!loaded=\relax\fi
\chardef\plot!savecc=\count255
%
%
\def\plot!zero{0}
\def\plot!one{1}
\newdimen\plotunitx
\newdimen\plot!unitxu
\def\plot!figscalex{1}
\plotunitx=1truebp
\plot!unitxu=1truebp
\newdimen\plotunity
\newdimen\plot!unityu
\def\plot!figscaley{1}
\plotunity=1truebp
\plot!unityu=1truebp
\def\plot!sepscfac{1}
\newtoks\plot!symbol
\newbox\plot!figurebox
\newdimen\plotcurrx
\newdimen\plotcurry
\newdimen\plotlinewidth
\plotlinewidth=.6bp
\def\plot!savemem{0}
\let\plot!global=\relax
\newif\ifplot
\plottrue
\newif\ifplotPS
\newif\ifplotseparate
\newif\ifplot!infig
\newif\ifplot!trans
\newif\ifplot!loctrans
\newif\ifplot!sepscaled
\newif\ifplot!nowraw
\newif\ifplot!globdef
\newif\ifplot!rawbounds
\newif\ifplot!spseq
\newif\ifplot!stacklock
\newif\ifplot!lhcs
\plot!spseqtrue
%
%
\mathchardef\plot!ci="220E
\mathchardef\plot!bu="220F
\mathchardef\plot!cd="2201
\def\plot!sp{ }
\def\plot!em{}
\let\plot!bg={
\let\plot!eg=}
\def\plot!ni{\prevdepth=-1000pt }
\def\plot!sm#1{{\setbox0=\hbox{#1}\ht0=0pt \dp0=0pt \box0 }}
\def\plot!lo#1\plot!re{\def\plot!bo{#1}\plot!it}
\def\plot!it{\plot!bo \let\plot!ne=\plot!it \else\let\plot!ne=\relax\fi
   \plot!ne}
\let\plot!re=\fi
%
%
\long\def\plot!firstlet#1#2\plot!endofarg{#1}
\let\plot!PStrueold=\plotPStrue
\def\plotPStype#1{
   \plot!PStrueold
   \ifcase#1
      \plotPSfalse
   \or
      \def\plot!local{" }
      \def\plot!global{! }
      \def\plot!rawbegin##1{ps:SDict begin ##1\plot!sp end}
      \def\plot!raw##1{ps:SDict begin ##1\plot!sp end}
      \def\plot!rawend##1{ps:SDict begin ##1\plot!sp end}
      \def\plot!trans{currentpoint /p!s1 2 index def /p!s2 1 index def
         translate exec 0 0 moveto p!s1 neg p!s2 neg translate }
      \def\plot!setorig{}
      \def\plot!rawsetcurr{}
      \def\plot!rawsetorig{currentpoint translate }
      \def\plot!PSfile##1{\includegraphics{##1}}
      \plot!globdeftrue
      \plot!stacklocktrue
      \plot!lhcstrue
   \or
      \def\plot!local{ps:: }
      \def\plot!global{ps::[global] }
      \def\plot!rawbegin##1{ps::[inline,begin] ##1}
      \def\plot!raw##1{ps::[inline] ##1}
      \def\plot!rawend##1{ps::[inline,end] ##1}
      \def\plot!trans{Xpos Ypos translate
         exec Xpos neg Ypos neg translate }
      \def\plot!setorig{Xpos Ypos translate }
      \def\plot!rawsetcurr{Xpos Ypos moveto }
      \def\plot!rawsetorig{Xpos Ypos translate }
      \def\plot!PSfile##1{\plot!rawstart{gsave Xpos Ypos translate}%
         \special{ps: plotfile ##1 inline}\plot!rawfinish{grestore}}
      \plot!rawboundstrue
   \or
      \def\plot!local{empty.ps }
      \def\plot!rawbegin##1{empty.ps ##1}
      \def\plot!raw##1{empty.ps ##1}
      \def\plot!rawend##1{empty.ps ##1}
      \def\plot!trans{pop }
      \def\plot!setorig{}
      \def\plot!rawsetcurr{0 0 moveto }
      \def\plot!rawsetorig{}
      \def\plot!PSfile##1{\special{##1}}
      \plot!spseqfalse
   \else
      \immediate\write16{plotPStype: Illegal type specified -- type set to 0}
      \plotPSfalse
   \fi
   \ifplotPS\ifplot!globdef\else\ifx\plot!global\relax\else
      \ifcase\plot!savemem\or
         \plot!globdeftrue\plot!gwarning
      \or
         \immediate\write16{ }
         \immediate\write16{Do you want me to try to conserve TeX memory?}
         \immediate\write16{(If so, the resulting pages must be printed in %
            order from the beginning.)}
         \message{? }
         \read-16 to\plot!inputstr
         \edef\plot!inputstr{\plot!firstlet\plot!inputstr xx\plot!endofarg}
         \if y\plot!inputstr\plot!globdeftrue\fi
         \if Y\plot!inputstr\plot!globdeftrue\fi
      \fi
   \fi\fi\else\global\let\plot!gwarning=\relax\fi
   \ifplot!globdef\global\let\plot!gwarning=\relax\fi
   \plot!PSinit   
}
\def\plotPSask{
   \immediate\write16{ }
   \immediate\write16{Which DVI-to-PS converter are you using?}
   \immediate\write16{  1. Rokicki dvips (Radical Eye)}
   \immediate\write16{  2. ArborText dvips (DVILASER/PS)}
   \immediate\write16{  3. OzTeX [partial functionality]}
   \immediate\write16{  0. None of the above %
      [PostScript insertions replaced or suppressed]}
   \message{Enter a number: }
   \read-16 to\plot!inputstr
   \edef\plot!inputstr{\plot!firstlet\plot!inputstr xx\plot!endofarg}
   \ifcat0\plot!inputstr\else\def\plot!inputstr{Z}\fi
   \ifnum\expandafter`\plot!inputstr<`0\def\plot!inputstr{99}\fi
   \ifnum\expandafter`\plot!inputstr>`9\def\plot!inputstr{99}\fi
   \plotPStype{\plot!inputstr}
}
\def\plotPStrue{\plotPSask}

\let\plot!septrueold=\plotseparatetrue
\def\plotseparatetrue{\plot!septrueold\plottrue}
\let\plot!sepfalseold=\plotseparatefalse
\def\plotseparatefalse{\plot!sepfalseold\plotnosepscale}
\let\plot!falseold=\plotfalse
\def\plotfalse{\plot!falseold\plotseparatefalse}
\def\plotask{
   \immediate\write16{ }
   \immediate\write16{How do you want figures to appear?}
   \immediate\write16{  1. Normally within text}
   \immediate\write16{  2. On separate pages (blank spaces left in text)}
   \immediate\write16{  0. Not at all (blank spaces left in text)}
   \message{Enter a number: }
   \read-16 to\plot!inputstr
   \if 0\plot!firstlet\plot!inputstr xx\plot!endofarg
      \plotfalse
   \else\if 1\plot!firstlet\plot!inputstr xx\plot!endofarg
      \plottrue \plotseparatefalse
   \else\if 2\plot!firstlet\plot!inputstr xx\plot!endofarg
      \plotseparatetrue
   \else
      \immediate\write16%
         {plotask: Illegal option specified -- figures suppressed}
      \plotfalse
   \fi\fi\fi
}
\gdef\plot!gwarning{
   \immediate\write16{ }
   \immediate\write16{plot WARNING: This DVI-to-PS converter does not %
      fully support global}
   \immediate\write16{ \plot!sp\plot!sp definitions -- printing pages out %
      of order may cause problems.}
   \global\let\plot!gwarning=\relax}
%
%
\def\plot!zerocp{\plotcurrx=0pt\plotcurry=0pt\relax}
{\toksdef\abc=7 \catcode`p=12 \catcode`t=12 \global\abc={pt}}
\expandafter\def\expandafter\plot!unpt\expandafter#\expandafter1\the\toks7{#1 }
\def\plotthedim#1{\expandafter\plot!unpt\the#1 }

\def\plotpdim#1{\plotthedim{#1}\ifnum1000=\mag 1.00375 \else\ifx\plot!unmag
   \plot!em 1003.75 \the\mag\plot!sp div \else 1.00375 \fi\fi div }
\def\plotfdims#1#2{\plotthedim{#1}\plotthedim{\plotunitx}div
   \plotthedim{#2}\plotthedim{\plotunity}div }
\def\plot!recip#1{{\dimen6=1in \dimen8=#1in \dimen0=0sp \global\toks7={}
   \plot!lo
      \count255=\dimen6
      \divide\count255 by\dimen8
      \global\toks7=\expandafter\expandafter\expandafter{\expandafter\the
         \expandafter\toks\expandafter7\expandafter\plot!sp\the\count255}
      \ifdim\dimen0=0sp\global\toks7=\expandafter{\the\toks7.}\fi
      \dimen4=\dimen8
      \multiply\dimen4 by\count255
      \advance\dimen6 by -\dimen4
   \ifdim\dimen0<6sp
      \advance\dimen0 by 1sp
      \multiply\dimen6 by 10
   \plot!re}}
\def\plot!applyssf{\plotunitx=\plot!sepscfac\plot!unitxu
   \plotunity=\plot!sepscfac\plot!unityu}
\def\plot!PSssf{\ifplot!sepscaled\plot!sepscfac\plot!sp dup scale \fi}
\def\plot!errmsg{\errmessage
   {plot error: command allowed only within \string\plotfigbegin
   ...\string\plotfigend}}
\def\plot!rawstart#1{\plot!bg\ifplot!nowraw\special{\plot!raw{#1}}\else
   \special{\plot!rawbegin{#1}}\plot!nowrawtrue\plot!rawboundsfalse\fi}
\def\plot!rawfinish#1{\plot!eg\ifplot!nowraw\special{\plot!raw{#1}}\else
   \special{\plot!rawend{#1}}\fi}
\def\plot!rawstartnull{\ifplot!rawbounds\plot!rawstart{}\else
   \plot!bg\plot!nowrawtrue\fi}
\def\plot!rawfinishnull{\plot!eg\ifplot!rawbounds\ifplot!nowraw
   \else\special{\plot!rawend{}}\fi\fi}
%
%
\def\plotfigscales#1 #2 {\def\plot!figscalex{#1} \plot!unitxu=#1truebp
   \def\plot!figscaley{#2} \plot!unityu=#2truebp \plot!applyssf}
\def\plotfigscale#1 {\plotfigscales{#1} {#1} }
\def\plotrotate#1 {\def\plot!symangle{#1}\plot!loctranstrue\plot!transtrue}
\def\plotscale#1 {\def\plot!symscale{#1}\plot!loctranstrue\plot!transtrue}
\def\plotcancelmag{\def\plot!unmag{\ifnum1000=\mag \else1000 \the\mag\plot!sp
   div dup scale \fi}\plot!loctranstrue\plot!transtrue}
\def\plotresumemag{\def\plot!unmag{}}
\long\def\plotgentrans#1\plotPSend{\plot!loctranstrue\plot!transtrue
   \long\def\plot!transcode{\ifplot!lhcs 1 -1 scale \fi
   #1\ifplot!lhcs\plot!sp 1 -1 scale\fi}}
\def\plotnotrans{\def\plot!transcode{\ifx\plot!symangle\plot!zero\else
   \plot!symangle\plot!sp\ifplot!lhcs neg \fi rotate \fi
   \ifx\plot!symscale\plot!one\else\plot!symscale\plot!sp dup scale \fi}
   \plotresumemag\plot!loctransfalse\ifplot!sepscaled
   \else\plot!transfalse\fi\def\plot!symscale{1}\def\plot!symangle{0}}
\def\plotsepscale#1 {\ifplotseparate\def\plot!sepscfac{#1}
   \ifx\plot!sepscfac\plot!one\plotnosepscale\else
   \plot!sepscaledtrue\plot!transtrue\plot!applyssf \fi\fi}
\def\plotnosepscale{\plot!sepscaledfalse\ifplot!loctrans\else\plot!transfalse
   \fi\def\plot!sepscfac{1}\plot!applyssf}
%
%
\def\plotfigbegin{\setbox\plot!figurebox=\vbox\plot!bg
   \plot!infigtrue\plotnotrans\plot!zerocp}
\def\plotfigend{\xdef\plot!gtemp{\plot!sepscfac}\plot!eg
   {\dimen2=\ht\plot!figurebox
   \dimen4=\wd\plot!figurebox
   \dimen6=\dp\plot!figurebox
   \ifplotseparate
      \shipout\vbox to9truein{\vss\hbox to6.5truein{\hss\box\plot!figurebox
         \hss}\vss}
      \plot!recip{\plot!gtemp}
      \setbox\plot!figurebox=\hbox to \the\toks7\dimen4{\hfil
         \vrule width0pt height\the\toks7\dimen2 depth\the\toks7\dimen6}
   \else\ifplot\else
      \setbox\plot!figurebox=\hbox to \dimen4{\hfil
         \vrule width0pt height\dimen2 depth\dimen6}
   \fi\fi
   \box\plot!figurebox}}
\def\plot!plot{\ifplot
   \vbox to0pt{\kern-\dimen8
   \hbox{%
      \kern\dimen6
      \let\plot!temp=\plot!em
      \ifplotPS\ifplot!trans\ifplot!spseq
         \edef\plot!temp{\plot!unmag\plot!PSssf\plot!transcode}\fi\fi\fi
      \ifx\plot!temp\plot!em\else\plot!rawstart{{\plot!temp} \plot!pla}\fi
      \plot!sm{\hbox to0pt{\hss\the\plot!symbol\hss}}%
      \ifx\plot!temp\plot!em\else\plot!rawfinish{grestore}\fi
   }
   \kern\dimen8}
   \plot!ni\fi}
\def\plot#1 #2 {\ifplot!infig
   {\dimen6=#1\plotunitx\dimen8=#2\plotunity\plot!plot}
   \else\plot!errmsg\fi}
\def\plothere{\ifplot!infig
   {\dimen6=\plotcurrx\dimen8=\plotcurry\plot!plot}
   \else\plot!errmsg\fi}
\def\plotmove#1 #2 {\ifplot!infig
   \plotcurrx=#1\plotunitx\plotcurry=#2\plotunity
   \else\plot!errmsg\fi}
\def\plotrmove#1 #2 {\ifplot!infig
   \advance\plotcurrx by #1\plotunitx \advance\plotcurry by #2\plotunity
   \else\plot!errmsg\fi}
%
%
\def\plottext#1{\plot!symbol={\lower\fontdimen22\textfont2\hbox{#1}}}
\def\plotfmla#1{\plottext{$#1$}}
\def\plotsfmla#1{\plot!symbol={\lower\fontdimen22\scriptfont2\hbox
   {$\scriptstyle #1$}}}
\def\plotssfmla#1{\plot!symbol={\lower\fontdimen22\scriptscriptfont2\hbox
   {$\scriptscriptstyle #1$}}}
\def\plotcentered#1{\plot!symbol={\vbox to0pt{\vss\hbox{#1}\vss}}}

\def\plotclosed{
   \ifplotPS
      \plot!symbol={\special{\plot!local \plot!bullet}}
   \else\plotfmla{\plot!bu}\fi}

%
%
\long\def\plotPSbegin#1\plotPSend{
   \ifplot!infig
   \ifplot\ifplotPS
      \let\plot!temp=\plot!em
      \ifx\plot!one\plot!figscalex\else\let\plot!temp=\plot!one\fi
      \ifx\plot!one\plot!figscaley\else\let\plot!temp=\plot!one\fi
      \ifx\plot!temp\plot!em\else\def\plot!temp{\plot!figscalex\plot!sp
         \ifdim\plot!unitxu=\plot!unityu dup dup scale \else
         \plot!figscaley\plot!sp 2 copy scale dup mul exch
         dup mul 2 div sqrt \fi div }\fi
      \special{\plot!local
         gsave \plot!setorig
         \ifnum1000=\mag\else1000 \the\mag\plot!sp div dup scale \fi
         \ifplot!sepscaled\plot!sepscfac\plot!sp dup scale \fi
         \plotpdim{\plotlinewidth}\plot!temp setlinewidth
         \plotfdims{\plotcurrx}{\plotcurry}moveto #1\plot!sp grestore}
   \fi\fi
   \else\plot!errmsg\fi}
\long\def\plotPSglobal#1\plotPSend{%
   \ifplotPS
      \ifx\plot!global\relax
         \errmessage{plotPSglobal: This PS type does not support global %
            insertions}%
      \else
         \special{\plot!global #1}%
         \plot!gwarning
      \fi
   \fi }
\long\def\plotPSmath#1\plotPSend{{%
   \ifplotPS
      \mathchoice
         {\special{\plot!raw{1.00001}}}%
         {\special{\plot!raw{1.0}}}%
         {\special{\plot!raw{0.7}}}%
         {\special{\plot!raw{0.5}}}
      \special{\plot!local gsave \plot!setorig #1\plot!sp
         \ifplot!stacklock 0 \fi grestore}
      \ifplot!stacklock\special{\plot!raw{count 0 ne {pop} if}}\fi
   \fi}}
\long\def\plotPSgdef#1#2\plotPSend{{%
   \ifplotPS
      \escapechar=-1 %
      \ifplot!globdef
         \special{\plot!global /p!0\string#1 {#2} def}%
         \xdef#1{p!0\string#1 }%
      \else
         \xdef#1{#2\plot!sp}%
      \fi
   \else\global\let#1=\relax\fi}}
\def\plot!befsetup{\setbox0=\hbox{%
   \hskip 1in\raise1in\hbox{\special{\plot!raw{\plot!bfa}}}}%
   \ht0=0pt \wd0=0pt \box0}

\long\def\plotPSbefore#1#2\plotPSend{{%
   \setbox8=\hbox{#1}%
   \ifplotPS
      \hbox{\plot!befsetup\special{\plot!raw{\plot!bfb #2\plot!sp
      grestore}}\box8}%
   \else \box8\fi }}
\long\def\plotPSduring#1#2\plotPSend{{%
   \toks0={{#1}}%
   \ifplotPS\ifplot!spseq
      \hbox{\plot!rawstart{gsave #2}\the\toks0 \plot!rawfinish{grestore}}%
   \else \the\toks0\fi \else \the\toks0\fi }}
\long\def\plotPSbefaft#1#2\plotPSend#3\plotPSend{{%
   \toks0={{#1}}\setbox8=\hbox{\the\toks0}\def\plot!bef{#2}\def\plot!aft{#3}%
   \ifplotPS\ifplot!spseq
      \setbox8=\hbox{\vrule width\wd8 height\ht8 depth\dp8}%
      \hbox to\wd8{\plot!befsetup
         \plot!rawstart{\plot!bab
            \ifx\plot!bef\plot!em\else
            {\plot!bac \plot!bef\plot!sp \plot!bad} \plot!tra \fi}%
         \hbox to0pt{\the\toks0 \hss}%
         \plot!rawfinish{\ifx\plot!aft\plot!em\else {\plot!bac
            \plot!aft} \plot!tra \fi grestore}\hfil}%
   \else \box8\fi \else \box8\fi }}
\long\def\plotPSraw#1\plotPSend{%
   \ifplotPS\ifplot!spseq\special{\plot!raw{#1}}\fi\fi }
%
%
\newdimen\plotdotspc
\plotdotspc=2bp
\def\plotDRline{\ifplot!infig
   \ifplot\ifplotPS
      {\plotfigscale 1
      \edef\plot!litemp{\plotfdims{\dimen0}{\dimen2}}%
      \plotPSbegin\plot!litemp\plot!lia \plotPSend}
   \else\plot!TeXline\fi\fi
   \else\plot!errmsg\fi}
\def\plot!TeXline{{
   \plotlinewidth=\plot!sepscfac\plotlinewidth
   \plotdotspc=\plot!sepscfac\plotdotspc
   \ifdim\dimen0=\plotcurrx
      \dimen6=\dimen0
      \dimen4=\dimen2
      \advance\dimen4 by -\plotcurry
      \ifdim\dimen4<0sp\dimen4=-\dimen4\fi
      \advance\dimen4 by \plotlinewidth
      \plotcentered{\vrule width\plotlinewidth height\dimen4}
      \dimen8=\dimen2
      \advance\dimen8 by \plotcurry
      \divide\dimen8 by 2
      \plot!plot
   \else\ifdim\dimen2=\plotcurry
      \dimen8=\dimen2
      \dimen4=\dimen0
      \advance\dimen4 by -\plotcurrx
      \ifdim\dimen4<0sp\dimen4=-\dimen4\fi
      \advance\dimen4 by \plotlinewidth
      \plotcentered{\vrule width\dimen4 height\plotlinewidth}
      \dimen6=\dimen0
      \advance\dimen6 by \plotcurrx
      \divide\dimen6 by 2
      \plot!plot
   \else\plot!TeXdotline\fi\fi}}
\def\plot!TeXdotline{      
   \plotcentered{\vrule width\plotlinewidth height\plotlinewidth}
   \ifdim\dimen0<\plotcurrx
      \dimen6=\dimen0
      \dimen8=\dimen2
      \dimen0=\plotcurrx
      \dimen2=\plotcurry
   \else
      \dimen6=\plotcurrx
      \dimen8=\plotcurry
   \fi
   \advance\dimen2 by -\dimen8
   \dimen4=\dimen0
   \advance\dimen4 by -\dimen6
   \global\dimen1=\dimen2
   \ifdim\dimen1<0sp\global\dimen1=-\dimen1\fi
   \global\advance\dimen1 by \dimen4
   \global\divide\dimen1 by \plotdotspc
   \global\advance\dimen1 by 1sp
   \global\dimen7=\dimen4
   \global\divide\dimen7 by \dimen1
   \global\dimen3=\dimen7
   \global\multiply\dimen3 by -\dimen1
   \global\advance\dimen3 by \dimen4
   \global\dimen9=\dimen2
   \global\divide\dimen9 by \dimen1
   \global\dimen5=\dimen9
   \global\multiply\dimen5 by -\dimen1
   \global\advance\dimen5 by \dimen2
   \ifdim\dimen5<0sp
      \global\advance\dimen5 by \dimen1
      \global\advance\dimen9 by -1sp
   \fi
   \dimen2=0sp
   \dimen4=0sp
   \plot!lo
      \plot!plot
   \ifdim\dimen6<\dimen0
      \advance\dimen2 by \dimen3
      \advance\dimen4 by \dimen5
      \advance\dimen6 by \dimen7
      \advance\dimen8 by \dimen9
      \ifdim\dimen2<\dimen1\else
         \advance\dimen6 by 1sp
         \advance\dimen2 by -\dimen1
      \fi
      \ifdim\dimen4<\dimen1\else
         \advance\dimen8 by 1sp
         \advance\dimen4 by -\dimen1
      \fi
   \plot!re
}
\def\plotline#1 #2 {
   {\dimen0=#1\plotunitx \dimen2=#2\plotunity
   \plotDRline}
   \plotmove {#1} {#2}
}
\def\plotrline#1 #2 {
   {\dimen0=\plotcurrx \dimen2=\plotcurry
   \advance\dimen0 by #1\plotunitx \advance\dimen2 by #2\plotunity
   \plotDRline}
   \plotrmove {#1} {#2}
}
%
%
\def\plotvskip#1 {\vskip #1\plotunity\relax}
\def\plothphant#1 {\vbox to0pt{\hbox to#1\plotunitx{\hfil}}\plot!ni}
%
%
\def\plot!PSinit{
   \plotPSgdef\plot!circ
      gsave
      \plot!setorig
      newpath
      0 0 2.16 0 360 arc fill
      1 setgray
      newpath
      0 0 1.44 0 360 arc fill
      grestore
   \plotPSend
   \plotPSgdef\plot!bullet
      gsave
      \plot!setorig
      newpath
      0 0 2.16 0 360 arc fill
      grestore
   \plotPSend
   \plotPSgdef\plot!cdot
      gsave
      \plot!setorig
      newpath
      0 0 0.44 0 360 arc fill
      grestore
   \plotPSend
   \plotPSgdef\plot!tra \plot!trans \plotPSend
   \plotPSgdef\plot!pla gsave \plot!tra \plotPSend
   \plotPSgdef\plot!bfa \plot!rawsetcurr currentpoint transform\plotPSend
   \plotPSgdef\plot!bfb gsave \plot!rawsetorig itransform
      72 div exch 72 div exch scale\plotPSend
   \plotPSgdef\plot!bab gsave \plot!rawsetorig itransform
      72 div /p!s4 exch def 72 div /p!s3 exch def grestore gsave\plotPSend
   \plotPSgdef\plot!bac p!s3 p!s4 scale\plotPSend
   \plotPSgdef\plot!bad 1 p!s3 div 1 p!s4 div scale\plotPSend
   \plotPSgdef\plot!lia currentpoint newpath moveto lineto stroke\plotPSend
   \plot!reclaim
}
\def\plot!reclaim{
   \def\plotPSask{\immediate\write16{plot: PStype has already been set}}
   \def\plotPSfalse{\plotPSask}
   \def\plotPStype##1{\plotPSask}
   \let\plot!PSinit=\relax
   \ifplotPS
      \let\plot!TeXline=\relax
      \let\plot!TeXdotline=\relax
   \fi
   \let\plot!reclaim=\relax
}
\catcode`\!=\plot!savecc
%
%
\plottrue
\plotPStype 1

\define\R{{\bold R}}
\define\Z{{\bold Z}}
\define\ve{{\bold e}}
\define\vg{{\bold g}}
\define\vt{{\bold t}}
\define\vu{{\bold u}}
\define\vv{{\bold v}}
\define\vw{{\bold w}}
\define\vx{{\bold x}}
\define\vy{{\bold y}}
\define\vz{{\bold z}}
\define\vzero{{\bold 0}}
\define\vdot{\cdot}
\define\grad{\operatorname{grad}}
\define\QED{ \null\nobreak\hfill$\blacksquare$\enddemo}
\define\QNED{ \null\nobreak\hfill$\blacksquare$\csname endproclaim\endcsname}
\define\procl#1{\smallskip{\it {#1}. }}
\define\eps{\varepsilon}
\define\dist{\delta}
\define\octa#1#2{{S_{#2}}}
\define\tetra#1#2{{S'_{#2}}}
\define\Rocta#1#2{{\bar S_{#2}}}
\define\Rtetra#1#2{{\bar S'_{#2}}}
\define\indx#1:#2|{\vphantom{#1#2}\left|\smash{#1}:\smash{#2}\right|}
\define\volume{\operatorname{vol}}
\define\badadjustI{@!@!@!@!@!@!}
\define\badadjustII{\!}

\define\AguFioGar{1}
\define\AnnBau{2}
\define\BerComHsu{3}
\define\BoeWan{4}
\define\CheJia{5}
\define\Chung{6}
\define\ChuFabMan{7}
\define\CohKarMatSch{8}
\define\CohLitLobMat{9}
\define\DinHaf{10}
\define\DouJan{11}
\define\DraFab{12}
\define\Fary{13}
\define\FidForZit{14}
\define\ForLam{15}
\define\Gritzmann{16}
\define\Hoylman{17}
\define\Rogers{18}
\define\Sabidussi{19}
\define\StaCow{20}
\define\Urakawa{21}
\define\WonCop{22}
\define\YebFioMorAle{23}

\topmatter
\title The degree-diameter problem for several varieties of
Cayley graphs, I:  the Abelian case \endtitle
\leftheadtext{Randall Dougherty and Vance Faber}
\rightheadtext{The degree-diameter problem for Abelian Cayley graphs}
\author Randall Dougherty and Vance Faber\endauthor
\affil Los Alamos National Laboratory \\ Ohio State University
\\ LizardTech, Inc.\endaffil
\date September 8, 2000 \enddate
\ifx\refstyle\undefinedZQA
   \thanks
   The first author was supported by NSF grant number DMS-9158092 and
   by a fellowship from the Sloan Foun\-da\-tion\endthanks
\else
   \thanks
   The first author was supported by NSF grant number DMS-9158092 and
   by a fellowship from the Sloan Foun\-da\-tion.\endthanks
\fi
\address LizardTech, Inc., 1008 Western Ave., Seattle, WA 98104\endaddress
\email rdougherty\@lizardtech.com \endemail
\address LizardTech, Inc., 1008 Western Ave., Seattle, WA 98104\endaddress
\email vxf\@lizardtech.com \endemail
\abstract We address the degree-diameter problem for Cayley graphs of
Abelian groups (Abelian graphs), both directed and undirected.  The
problem turns out to be closely related to the problem of finding
efficient lattice coverings of Euclidean space by shapes such as
octahedra and tetrahedra; we exploit this relationship in both
directions.
For\/ $2$~generators (dimensions) these methods yield optimal
Abelian graphs with a given diameter~$k$.  (The results in
two dimensions are not new; they
are given in the literature of distributed loop networks.)
We find an
undirected Abelian graph with\/ $3$~generators and a given diameter~$k$
which we conjecture to be as large as possible; for the directed case, we
obtain partial results. These results are connected to efficient lattice coverings
of~$\R^3$ by octahedra or by tetrahedra; computations on Cayley graphs
lead us to such lattice coverings which we conjecture to be optimal.
(The problem of finding such optimal coverings can be reduced to a finite
number of nonlinear optimization problems.)
We discuss the asymptotic behavior of the Abelian degree-diameter problem
for large numbers of generators.  The graphs obtained here are substantially
better than traditional toroidal meshes, but, in the simpler undirected
cases, retain certain desirable features such as good routing
algorithms, easy constructibility, and the ability to host
mesh-connected numerical algorithms without any increase in
communication times.  \endabstract

\endtopmatter
\document

\subhead Introduction \endsubhead
The degree-diameter problem for graphs is the following question: What
is the largest number of vertices a graph (undirected or directed)
can have if one is given upper bounds on the degree of each vertex
and on the diameter of the graph (the maximum path-distance from
any vertex to any other)?  One application for such graphs is in
the design of interconnection networks for parallel processors,
where one wants to have a large number of processors without
requiring a large number of wires at a single processor or
incurring long delays in messages from one processor to another.  For more
information on the (undirected) degree-diameter problem,
see Dinneen and Hafner~\cite{\DinHaf}; up-to-date results can be
found at the following Web site:
$$\text{\tt http://www-mat.upc.es/grup\_de\_grafs/grafs/taula\_delta\_d.html}$$

A desirable extra property of such networks is that they look
identical from any processor; this means that the graphs one uses
should be vertex-transitive, i.e., for any two vertices $x$ and $y$
there is an automorphism of the graph which maps $x$ to $y$.  Here we
will restrict ourselves to a special class of vertex-transitive graphs,
called Cayley graphs.  A Cayley graph is specified by a group and a set
of generators for
this group; the vertices of the graph are the elements of the group,
and there is an edge from $x$ to $y$ if and only if there is a
generator $g$ such that $y = xg$.  (It can be shown that every
vertex-transitive graph is isomorphic to a generalized form
of Cayley graph called a Cayley coset graph~\cite{\Sabidussi}.  In
this paper, though, we will look only at Cayley graphs.  For Abelian
groups, this is no loss of generality since every Cayley coset graph of
an Abelian group is isomorphic to a Cayley graph of an Abelian group.)
In a
directed Cayley graph from a group on $d$ generators, every vertex has
in-degree and out-degree~$d$; if $d$~generators are used to form
an undirected Cayley graph, then the degree of each vertex is
the number of generators of order~$2$ plus twice the number of
generators of order greater than~$2$
(unless there are redundant generators).  So we will usually talk about
Cayley graphs on a given number of generators rather than of a
given degree; the cases where some generators have order~$2$ and
hence only contribute $1$ rather than~$2$ to the degree of an
undirected Cayley graph will be handled separately.

The most straightforward approach for trying to find large Cayley
graphs of small diameter on a given number $d$ of generators is to
examine various groups, look at some or all possible sets of $d$
generators for such a group, and check whether each such set
in fact generates the group and, if so, what the diameter of
the graph is.  But this can be a very large task even for relatively
small groups and generating sets.  In this paper, we will use a
different approach which facilitates studying many groups and
generating sets at once; it yields provably optimal results
for some families of groups, and good lower and upper bounds
for others.

In its most general form, the idea is as follows.  Let $F$ be the free
(universal) group on $d$~generators.  Then, for any group $G$ and any
set of $d$ generators for $G$, there is a homomorphism $\pi\colon F \to
G$ which maps the canonical generators for $F$ to the given generators
for $G$; clearly $\pi$ is surjective.  Let $N$ be the kernel of $\pi$.
Then $N$ is a normal subgroup of $F$, and $\indx F:N| = |G|$; in fact,
$G$ is isomorphic to~$F/N$, and the Cayley graph of $G$ with the given
generators is isomorphic to the Cayley graph of~$F/N$ with the
canonical generators for $F$.  Let $S$ be the set of elements of~$F$
which can be expressed as a word of length at most $k$ in the
generators.  (If one is interested in undirected Cayley graphs, then
such words may use inverse generators $g^{-1}$ as well as generators;
for the directed case, only words using generators and not inverse
generators are allowed.)

\proclaim{Proposition 1} The Cayley graph for $G$ on the given
generators has diameter at most~$k$ if and only if $SN = F$. \endproclaim

\demo{Proof} First, suppose $SN = F$.  Let $a$ be an arbitrary
element of $G$; then $a = \pi(x)$ for some $x$, and $x$ can be
written in the form $wy$ with $w \in S$ and $y \in N$.
Hence, $a = \pi(wy) = \pi(w)\pi(y) = \pi(w)$.  Now $w$ can be written
as a word of length at most $k$ in the generators of $F$, so $a = \pi(w)$
can be written as the same word in the corresponding generators of $G$.
Since $a$ was arbitrary, the Cayley graph has diameter at most $k$.

Now suppose the Cayley graph has diameter at most $k$.  Let $x$
be any element of $F$; then $\pi(x)$ can be written as a word $w'$ of
length at most $k$ in the generators of $G$.  Let $w$ be the corresponding
word in the generators of $F$; then $\pi(w) = w' = \pi(x)$, so
$\pi(w^{-1}x)$ is the identity of $G$, so $w^{-1}x \in N$.
Hence, $x = w(w^{-1}x) \in SN$, as desired. \QED

So finding a Cayley graph on $d$ generators with diameter $k$ whose
size is as large as possible is equivalent to finding a normal
subgroup $N$ of $F$ such that $SN = F$ and $\indx F:N|$ is as large
as possible.  Of course, we immediately get the upper bound
$\indx F:N| \le |S|$, but this is probably not attainable.

Unfortunately, the collection of normal subgroups of $F$ is so
large and varied as to be unmanageable.  So what we will do instead
is restrict ourselves to certain families (usually varieties) of
groups; this allows us to replace $F$ with a free group for the
family in question, which may be much easier to work with.
For instance, if we only consider the Cayley graphs of Abelian
groups, then we can replace $F$ with the free {\sl Abelian} group
on $d$ generators, which is simply $\Z^d$; the normal subgroups
of $\Z^d$ are well understood and relatively easy to work with.
It turns out that this reduces the degree-diameter problem for Abelian
Cayley graphs to interesting problems about lattice coverings of
Euclidean space by various shapes.  Some of these problems can
be solved completely, giving optimal Abelian Cayley graphs; others
are still open.

A simple path-counting argument gives upper bounds for the size $n$ of
a Cayley graph with $d$~generators and diameter limit $k$:  in the
directed case, $$n \le 1 + d + d^2 + \dots + d^k = {d^{k+1} - 1 \over d
- 1},$$ and in the undirected case, $$n \le 1 + 2d + 2d(2d-1) + \dots
+ 2d(2d-1)^{k-1} = {d(2d-1)^k-1 \over d-1}.$$  (The formulas for $d=1$
are $k+1$ and $2k+1$.)  These limits are well-known and actually
apply to the degree-diameter problem for arbitrary graphs;
they are known as the Moore bounds. For $d=1$,
these limits are attained by simple cycle graphs, which are Cayley
graphs of cyclic groups.

In most cases, we will find that the
attainable values for $n$ using Cayley graphs in the families we
consider do not approach these upper bounds; the equations defining the
families force many paths to be redundant.  But the extra structure
provided by the groups may provide compensating advantages in parallel
computers, such as good routing algorithms, easy constructibility, and
the ability to map common problems onto the architecture. In
particular, many of the Cayley graphs of Abelian groups that we discuss
in this paper are multi-dimensional rectangular meshes with
additional connections at the boundary.  Thus, mesh calculations with
natural boundary conditions are trivially mapped into these graphs,
while the extra connections are utilized only when global
communications are being carried out.  In addition, the mesh nature of
these graphs allows the physical construction of the network to be
carried out with relatively short wires.  This will be discussed
further below.

In a separate paper, we will examine some other varieties
of groups for which similar analyses of Cayley graphs
are feasible.

We would like to (and hereby do) thank Michael Dinneen, whose computer
work determining the diameters of many Cayley graphs yielded some of
the specific groups listed here and suggested promising families of
groups to examine.  Thanks also to Francesc Comellas, for pointing us to
the existing literature on directed loop networks (in particular,
the very useful survey paper of Bermond, Comellas, and
Hsu~\cite{\BerComHsu}), and to Alexander Hulpke for discussions
on automorphisms of Abelian groups.

\subhead Abelian groups \endsubhead
In the rest of this paper, we will examine the Cayley graphs arising
from Abelian groups.  Toroidal meshes and hypercubes are examples of
such graphs.

The degree-diameter problem for Abelian Cayley graphs has been
considered by others.  In particular, Annexstein and
Baumslag~\cite{\AnnBau} show that the number of
generators $d$, diameter $k$, and size $n$ of a directed Abelian Cayley graph
satisfy $$k\geq \Omega (n^{1/d});$$ in fact, if $d \le n^{1/d}$, then
$$k\geq \Omega (dn^{1/d}).$$  They also discuss similar results for
Cayley graphs of nilpotent groups.

In addition, Chung~\cite{\Chung} has constructed directed Abelian
Cayley graphs $G$ with $n=p^t - 1$, $d=p$, and $$k\le \left\lceil 2t +
\frac{4t\log t}{\log p - 2 \log (t-1)}\right\rceil$$ for any positive
integer $t<\sqrt{p}+1$, where $p$ is a prime.  (Note that Chung's
examples have diameters which are small compared to the number of
generators; in contrast, we will concentrate here on graphs with a
small fixed number of generators and relatively large diameters.)  Her
methods involve estimates of the second eigenvalue of the Laplacian of
the adjacency matrix of $G$.  For more about estimating the diameter of
general graphs from knowledge of this eigenvalue, see Chung, Faber, and
Manteuffel~\cite{\ChuFabMan}.  This eigenvalue is also
connected to the sphere {\sl packing} problem for real lattices; see
Urakawa~\cite{\Urakawa}.

The more specific case of Cayley graphs of {\sl cyclic}\/ groups
has been studied more extensively; such graphs are usually referred
to by some variant of the phrase `loop networks.'  The survey
paper of Bernard, Comellas, and Hsu~\cite{\BerComHsu} is an excellent
guide to the literature in this area.

We start by taking care of some generalities and notational matters.
In this paper we will use the symbol $+$ for the group operation(s).
Let $\Z_m$ be the cyclic group of order $m$ (for definiteness, the
set $\{0,1,\dots,m-1\}$ with the operation of addition modulo $m$).
The groups $\Z_m$ and the infinite cyclic group $\Z$ each have
a canonical generator, $1$; but they also have other sets of generators,
and some of these will be important later.

When one has groups $G_1,G_2,\dots,G_l$ and a set of generators for
each, then one can get a set of generators for the product
$G_1 \times G_2 \times \dots \times G_l$ by putting together
the given generator sets.  More precisely, for each $i \le l$ and each
generator $g$ of $G_i$, let $\ve_i(g)$ be the element of the product group
which has the identity element of $G_j$ as its $j$\snug'th coordinate
for all $j\le l$ except~$i$; the $i$\snug'th coordinate is~$g$.  Then
the set of all elements~$\ve_i(g)$ is a natural generating set for the
product group.  (The resulting diameter for the product group is the
sum of the diameters of the groups~$G_i$.)
In the case where the groups~$G_i$ are cyclic groups
with the canonical single generators, we write simply $\ve_i$ for
the $l$\snug-tuple with $1$~at the $i$\snug'th coordinate and
$0$~elsewhere.

A two-dimensional toroidal mesh is simply the Cayley graph of the
group $\Z_m \times \Z_n$ with the canonical generators $\ve_1=(1,0)$
and $\ve_2=(0,1)$; higher-dimensional meshes are obtained from longer
products.  For the two-dimensional case, the number of vertices
is~$mn$, the degree is $2$~in the directed case and $4$~in the
undirected case (assuming $m,n \ge 3$), and the diameter is
$m+n-2$ in the directed case, $\lfloor m/2 \rfloor + \lfloor n/2 \rfloor$
in the undirected case.  The calculations in three or more dimensions
are analogous.

The $d$\snug-dimensional hypercube is the Cayley graph of the
group $\Z_2^d$ with the canonical generators; since one gets
bidirectional links in any case, we may as well just talk about the
undirected version.  In this case the size of the graph is exponential
in the diameter, but only because the degree also grows with $d$:
the size is $2^d$, the degree is $d$ (not $2d$, because the generators
have order $2$), and the diameter is also $d$.

We will see that with a fixed small number $d$ of generators, one can
obtain nearly optimal results for undirected Abelian Cayley graphs
by using a twisted toroidal mesh; the twist allows one to multiply
the number of nodes in the ordinary toroidal mesh by $2^{d-1}$
without increasing the diameter.  For two dimensions we can get
exactly optimal results by slightly adjusting this graph;
for higher dimensions, the optimal size is still open.  For the
directed case, we can again get optimal results in two dimensions,
but the higher-dimensional case is again unsolved (and rather strange
even in three dimensions).

To get these results, we argue as in Proposition~1, but specialize to
the case of Abelian groups.  Hence, instead of the free group on
$d$~generators, we can use the free Abelian group on $d$~generators, which
is simply $\Z^d$ with the canonical generators $\ve_i$, $1 \le i \le d$.
For any Abelian group~$G$ generated by $g_1,\dots,g_d$, there is a
unique homomorphism from $\Z^d$ onto~$G$ which sends $\ve_i$ to~$g_i$ for
all~$i$.  Let~$N$ be the kernel of this homomorphism; then $G$ is
isomorphic to~$\Z^d/N$, and the Cayley graph of~$G$ with the given
generators is isomorphic to the Cayley graph of~$\Z^d/N$ with the
canonical generators for~$\Z^d$.

Given a diameter limit~$k$, let $\octa dk$ be the set of elements of~$\Z^d$
which can be expressed as a word of length at most $k$ in the
generators~$\ve_i$, which are allowed to occur positively or negatively.
(The dimension~$d$ will be clear from context.)  Then $\octa dk$~can also
be described as the set of points in~$\Z^d$ at distance at most~$k$
from the origin under the $\l^1$ (Manhattan) metric:
$$\octa dk = \{(x_1,\dots,x_d) \in \Z^d \colon |x_1|+\dots+|x_d|\le k\}.$$
Let $\tetra dk$ be the subset of~$\octa dk$ consisting of those elements
whose coordinates are all nonnegative; these are the elements which can
be expressed as words of length at most~$k$ in the generators~$\ve_i$
where only positive occurrences of the generators are allowed.
Then $\octa dk$~looks like a regular dual $d$\snug-cube (a square for
$d=2$, an octahedron for $d=3$), while $\tetra dk$~looks like
a right $d$\snug-simplex (a triangle for $d=2$, a tetrahedron for
$d=3$).

Now, by the same proof as for Proposition 1, we get:

\proclaim{Proposition 2} Let $G$, $N$, and $g_1,\dots,g_d$ be as above.
Then the undirected Cayley graph for $G$ and $g_1,\dots,g_d$ has
diameter at most $k$ if and only if $\octa dk + N = \Z^d$, and
the directed Cayley graph for $G$ and $g_1,\dots,g_d$ has
diameter at most $k$ if and only if $\tetra dk + N = \Z^d$. \QNED

So $|\octa dk|$ and $|\tetra dk|$ give upper bounds for the undirected
and directed versions of this case of the degree-diameter problem.
It is not hard to show that 
$$|\tetra dk| = {k+d \choose d},$$
so $|\tetra dk| = k^d / d! + O(k^{d-1})$ for fixed $d$.  For
$|\octa dk|$, we easily get the asymptotic form
$|\octa dk| = k^d 2^d / d! + O(k^{d-1})$ for fixed $d$, but exact formulas
are harder; Stanton and Cowan~\cite{\StaCow} give several, such as
$$|\octa dk| = \sum_{i=0}^d 2^i {d \choose i} {k \choose i}.$$
In particular, when $d$ is $1$, $2$, or $3$, the formula for
$|\octa dk|$ is $2k+1$, $2k^2+2k+1$, or $(4k^3+6k^2+8k+3)/3$,
respectively.

\subhead Lattice coverings and tilings \endsubhead
Proposition 2 tells us that, to find an optimal undirected (directed)
Cayley graph of diameter $k$ on $d$ generators, we should look for a
subgroup~$N$ of~$\Z^d$ such that $\octa dk + N$ ($\tetra dk + N$) is
all of~$\Z^d$ and the index $\indx \Z^d : N|$ is as large as possible;
the largest index we can hope for is $|\octa dk|$ ($|\tetra dk|$).  But
the structure of subgroups of~$\Z^d$ of finite index (which are all
normal, of course) is well known; they are precisely the
$d$\snug-dimensional lattices in $\Z^d$.  Because of this, we will
use the letter $L$ instead of $N$ for such subgroups of $\Z^d$ for
the rest of this paper.

A $d$\snug-dimensional lattice $L$ in $\Z^d$ is specified by
$d$ linearly independent vectors $\vv_1,\dots,\vv_d$ in $\Z^d$;
$L$ is the set of all integral linear combinations of these
vectors.  We have $\indx \Z^d:L| = |\det M|$, where $M$
is the $d\times d$ matrix whose $i$\snug'th row is $\vv_i$,
for $i=1,\dots,d$.

Note that any bounded set contains only finitely many members of $L$.
It follows that, if $S$ is a bounded subset of $\Z^d$ and $\vx$ is a
point in $\Z^d$, then there are only finitely many $\vv \in L$ such
that $\vx \in S + \vv$.

Also note that $L$, or indeed the entire group $\Z^d$, has a linear
ordering $\prec$ which is compatible with addition: $\vx \prec \vy$
implies $\vx+\vz \prec \vy+\vz$.  To define $\prec$, first choose a
direction (a nonzero vector $\vv$ in~$\R^d$), and put $\vx \prec \vy$
if $\vy$ is farther in this direction than~$\vx$~is ($\vx \cdot \vv <
\vy \cdot \vv$).  If two vectors are at the same distance in this
direction, then compare them in a second direction; repeat until all
ties are broken.  One example of this is lexicographic order:  compare
according to the first coordinate, then according to the second
coordinate if the first coordinates are equal, and so on.  Or, in this
discrete case, one can choose the initial direction so that it
distinguishes all points and no tie-breaking is necessary;
for instance, if $\vv = (1,\pi,\pi^2,\dots,\pi^{d-1})$,
then we never have $\vx \cdot \vv = \vy \cdot \vv$ for distinct~$\vx,\vy$
in~$\Z^d$, so we can just define $\vx \prec \vy$ to mean
$\vx \cdot \vv < \vy \cdot \vv$.

A {\it lattice covering} of $\Z^d$ by a set $S \subseteq \Z^d$ is
a collection of translates of $S$ by members of a lattice~$L$
(i.e., $\{S+\vv\colon \vv \in L\}$) which covers $\Z^d$.  If the
translates are disjoint, so that each point of~$\Z^d$ is covered
exactly once, then we have a {\it lattice tiling} of $\Z^d$
by $S$.

If we have a lattice covering as above, then $|S| \ge \indx \Z^d:L|$;
if it is a tiling, then $|S| = \indx \Z^d:L|$.  So we can measure
the extent to which a covering is `almost' a tiling by one of two
numbers: the {\it density} of the covering, which is
$\indx \Z^d:L|/|S| \ge 1$ (this is the average number of sets in the
covering to which a random point of $\Z^d$ belongs), or the
{\it efficiency} of the covering, which is $|S|/\indx \Z^d:L| \le 1$.
So Proposition~2 tells us that, in order to get the best
possible Abelian Cayley graph on $d$ generators with diameter $k$,
we must find a lattice covering of $\Z^d$ by $\octa dk$ or
$\tetra dk$ whose density (efficiency) is as small (large) as
possible.  We now give one more reformulation of the question.

\proclaim{Lemma 3} Suppose we have a lattice covering of\/ $\Z^d$ using
a bounded set~$S$ and the lattice~$L$.  Then there is a set $T
\subseteq S$ such that the translates of\/~$T$ by~$L$ form a lattice
tiling of\/~$\Z^d$. \endproclaim

\demo{Proof} Let $\prec$ be a linear order of~$L$ compatible with
addition.  Now let $T$ be the set of all points in $S$ which are
not in any of the sets $S + \vv$ for $\vv \in L$, $\vv \succ \vzero$.
We will see that every point~$\vx$ is in exactly one of the sets
$T + \vv$ for $\vv \in L$.

Fix $\vx$.  As noted before, $\vx$ is in only finitely many
of the translates $S + \vv$, so let $\vw$ be the $\prec$\snug-greatest
member of $L$ such that $\vx \in S+\vw$.  Let $\vy = \vx-\vw \in S$.
Then $\vy$ cannot be in $S+\vv$ for $\vv \succ \vzero$ in $L$,
because, if it were, we would have $\vx = \vy + \vw
\in S + \vv + \vw$ and $\vv + \vw \succ \vw$, contradicting
the maximality of $\vw$.  So $\vy \in T$ and $\vx \in T+\vw$.

Now suppose $\vx = \vy + \vw = \vy' + \vw'$ where $\vw$ and~$\vw'$ are
distinct members of~$L$ and $\vy$~and~$\vy'$ are in~$T$ (and hence in
$S$).  Then $\vv = \vy - \vy' = \vw' - \vw$ is a nonzero member of~$L$,
so either $\vv \succ \vzero$ or $\vv \prec \vzero$.  In the former
case, $\vy = \vy'+\vv \in S+\vv$ contradicts $\vy \in T$; in the latter
case, $\vy' = \vy - \vv \in S + (-\vv)$ contradicts $\vy' \in T$. \QED

Note that such a set $T$ must be of cardinality $\indx \Z^d:L|$, which
is the size of the Cayley graph of $\Z^d/L$.  So the size of the
largest undirected (directed) Abelian Cayley graph on $d$ generators
with diameter $k$ is equal to the size of the largest subset $T$ of
$\octa dk$ ($\tetra dk$) such that there is a lattice tiling of
$\Z^d$ using $T$.

\subhead Approximation by lattice coverings of real space \endsubhead
The study of lattice coverings and lattice tilings is more familiar
for~$\R^d$ than for~$\Z^d$.  We will show that real coverings can
be approximated to some extent by integer coverings, and vice versa,
so that known results from the real context can be transferred to
the integer lattices we are interested in.

The definitions of lattice, lattice covering, and lattice tiling
are the same in~$\R^d$ as in~$\Z^d$, except that we allow boundaries
to be shared in the definition of a tiling.  This lets us work throughout
with closed sets (usually polyhedra with their interiors) instead
of having to keep some of the boundary points and discard others.
Most of the results above for integer lattices go through verbatim
for real lattices, including Lemma~3 and its proof (although we
will probably use the closure of~$T$ rather than $T$~itself in practice).
The main difference is that the absolute
determinant $|\det M|$ of the matrix formed
from the generators of a lattice~$L$ is not $\indx \R^d:L|$
(which is infinite).  Instead, this determinant is the measure
of the parallelepiped determined by the generating vectors; this
parallelepiped gives a lattice tiling of~$\R^d$ using the lattice~$L$.
It follows easily that any other set~$S$ which gives a
lattice tiling using the lattice~$L$ must have measure $|\det M|$
(barring pathological cases of non-measurable sets or
positive-measure boundaries).  Such a set~$S$ is called a
{\it fundamental region} for the lattice~$L$.
One can now define the density or
efficiency of a covering by dividing the measure of the covering set
by this determinant or vice versa.

One can transform a lattice covering of $\Z^d$ using $S$ into a
lattice covering of $\R^d$ by replacing each point of $S$ with
a unit cube (i.e., replace $S$ with $S+U$ where $U$ is a fixed unit
$d$\snug-cube with edges parallel to the coordinate axes); the two
coverings will have the same density.  However, transforming results in the
other direction is harder, because the real results usually involve
actual triangles, octahedra, etc. rather than polycube approximations,
and the lattices used will often not be integer lattices.  We will
now present results that allow us to get around these difficulties.

In $\R^d$, let $\Rocta dk$ be the closed $\l^1$\snug-ball of radius $k$ at
the origin:
$$\Rocta dk = \{(x_1,\dots,x_d)\colon |x_1|+\dots+|x_d|\le k\}.$$
Let $\Rtetra dk$ be the set of nonnegative points in $\Rocta dk$:
$$\Rtetra dk = \{(x_1,\dots,x_d)\colon x_1,\dots,x_d\ge 0, x_1+\dots+
x_d\le k\}.$$
Let $L$ be any lattice in $\R^d$.

\proclaim{Proposition 4}
\roster
\item"(a)" If $\octa dk + L$ covers\/ $\Z^d$, then
$\Rocta d{k+d/2} + L$ covers\/ $\R^d$.
\item"(b)" If $\tetra dk + L$ covers\/ $\Z^d$, then
$\Rtetra d{k+d} + L$ covers\/ $\R^d$.
\endroster\endproclaim

\demo{Proof}
(a) By the triangle inequality for $l^1$~distance, we have
$\Rocta dk + \Rocta d{d/2} \subseteq \Rocta d{k+d/2}$, so
$\octa dk + L + \Rocta d{d/2} \subseteq \Rocta d{k+d/2} + L$;
therefore, it suffices to show that $\Z^d + \Rocta d{d/2} = \R^d$.
For any $\vx \in \R^d$, let $\vy$ be the element of~$\Z^d$ nearest
to~$\vx$ (i.e., round each coordinate of $\vx$ to the nearest integer);
then $\vx - \vy$ is in $\Rocta d{d/2}$, so $\vx$ is in $\Z^d + \Rocta d{d/2}$.

(b) Similarly, this follows from the fact that $\Z^d + \Rtetra dd = \R^d$,
which is proved in the same way as above (round each coordinate of $\vx$
downward instead of to the nearest integer).
\QED

One can argue in the same way within $\Z^d$ to get: 

\proclaim{Proposition 5} If $L$ is a lattice in\/ $\Z^d$ and $m$ is a
positive integer, then:
\roster
\item"(a)" If $\octa dk + L$ covers\/ $\Z^d$, then
$\octa d{mk+\lfloor m/2\rfloor d} + mL$ covers\/ $\Z^d$.
\item"(b)" If $\tetra nk + L$ covers\/ $\Z^d$, then
$\tetra d{mk+(m-1)d} + mL$ covers\/ $\Z^d$.
\endroster\endproclaim

\demo{Proof}
For (a), clearly $\octa d{mk} + mL$ covers $m\Z^d$, so, as in the preceding
proposition, it suffices to note that $m\Z^d + \octa d{\lfloor m/2\rfloor d}$
covers $\Z^d$, because we can just round any member of $\Z^d$ to the
nearest member of $m\Z^d$.  Similarly, (b) holds because
$m\Z^d + \tetra d{(m-1)d}$ covers $\Z^d$.
\QED

Of course, one can get similar results for sets other than $\octa dk$
and $\tetra dk$.

Using Proposition~4, it is easy to get from a covering of $\Z^d$ using
an integer lattice to a covering of $\R^d$ using a real lattice
(and using the real shape $\Rocta dk$ or $\Rtetra dk$); if $k$ is
large relative to $d$, then the two coverings have about the
same efficiency.  We will now show
that one can move in the other direction as well.

First, we give a useful criterion for deciding whether one has a
lattice covering of $\R^d$.

\proclaim{Proposition 6} Suppose $S$ is a nonempty subset of\/ $\R^d$ and $L$
is a lattice in\/ $\R^d$.  If there is an $\eps > 0$ such that
$S+L$ covers all points within distance $\eps$ of $S$, then
$S+L$ covers\/~$\R^d$. \endproclaim

\demo{Proof}  Clearly there is some point $\vx_0$ in $S+L$.  We will
show that, if $\vx \in S+L$ and the distance $\dist(\vx,\vy)$ is less
than $\eps$, then $\vy \in S+L$.  Applying this once shows that all
points within distance~$\eps$ of~$\vx_0$ are in $S+L$; applying it again
shows that all points within distance~$2\eps$ of~$\vx_0$ are in $S+L$;
since this can be repeated forever, we find that all points of~$\R^d$
are in $S+L$.

Let $\vx$ and $\vy$ be as above.  Find $\vv \in L$ such that
$\vx \in S+\vv$.  Then $\vy-\vv$ is within distance~$\eps$ of
$\vx-\vv \in S$, so there is $\vv' \in L$ such that
$\vy - \vv \in S+\vv'$.  Hence, $\vy \in S+\vv'+\vv$, so
$\vy \in S+L$, as desired. \QED

Using this, we can now show that, if one has a lattice covering using
a bounded subset of $\R^d$, then one can perturb the lattice slightly
and still get a lattice covering using a slightly larger subset of $\R^d$.

\proclaim{Proposition 7} Let $S$ be a bounded subset of\/ $\R^d$, and
let $L$ be a lattice in\/ $\R^d$ such that $S+L = \R^d$; let\/
$\vv_1,\dots,\vv_d$ be a list of generators for $L$.  Then there are
positive numbers $\eta$ and~$\rho$ such that, for all $r \in (0,1)$, if
the distance $\dist(\vv_i,\vv'_i)$ is less than $r\eta$ for all $i \le
d$, then $S^+ +L' = \R^d$, where $L'$ is the lattice generated by\/
$\vv'_1,\dots,\vv'_d$ and $S^+$ is the set of points within distance
$r\rho$ of $S$. \endproclaim

\demo{Proof} The number of members of $L$ within any bounded part
of $\R^d$ is finite, so it only takes finitely many translates of
$S$ by members of $L$ to cover any bounded part of $\R^d$.  In
particular, there is a number $M > 0$ such that the sets
$S+a_1\vv_1+\dots+a_d\vv_d$ for $(a_1,\dots,a_d) \in \Z^d$
with $|a_1|+\dots+|a_d| \le M$
cover all points within distance $\rho$ of $S$, for some $\rho > 0$.
Let $\eta = \rho/M$.

Let $r$ be any positive number less than $1$; we must see that, if
$S^+$ and $L'$ are defined as above, then $S^+ + L' = \R^d$.  By
the preceding proposition, it will suffice to show that $S^+ + L'$
covers all points within distance $(1-r)\rho$ of $S^+$.
Suppose $\vy$ is within distance $(1-r)\rho$ of $S^+$; then $\vy$
is within distance $(1-r)\rho + r\rho = \rho$ of $S$, so there exist
integers
$a_1,\dots,a_d$ with $|a_1|+\dots+|a_d| \le M$ and a point
$\vx \in S$ such that
$\vy = \vx+a_1\vv_1+\dots+a_d\vv_d$.  Let
$\vx' = \vx + \sum_{i=1}^d a_i(\vv_i-\vv'_i)$; then
$$\dist(\vx,\vx') \le \sum_{i=1}^d |a_i|\dist(\vv_i,\vv'_i)
< \sum_{i=1}^d |a_i|r\eta \le Mr\eta = r\rho,$$ so
$\vx' \in S^+$.  Since
$\vy = \vx'+a_1\vv'_1+\dots+a_d\vv'_d$, we have $\vy \in S^+ + L'$,
as desired. \QED

Note that, in the above propositions, `distance' need not be Euclidean
distance; it can be any metric arising from a norm on $\R^d$.  For our
present purposes, it will be most convenient (but not essential) to use
$\l^\infty$\snug-distance: $\dist(\vx,\vy) = \max_i |x_i-y_i|$.

\proclaim{Theorem 8} Suppose one has a lattice $L$ in\/ $\R^d$ such
that $\Rocta dk + L$ covers\/ $\R^d$; let\/ $\vv_1,\dots,\vv_d$ be
generators for $L$.  Then there is a constant $c$ such that, for all
sufficiently large real numbers $t$, if\/ $\vw_i$ is obtained from
$t\vv_i$ by rounding all coordinates to the nearest integer, and $\bar
L$ is the lattice generated by\/ $\vw_1,\dots,\vw_d$, then $\Rocta
d{tk+c} + \bar L$ covers\/ $\R^d$.  The same statement holds for
$\Rtetra dk$ and $\Rtetra d{tk+c}$ instead of $\Rocta dk$ and $\Rocta
d{tk+c}$.  \endproclaim

\demo{Proof} For the $\Rocta dk$ case, let $S = \Rocta dk$ and
find $\eta$ and $\rho$ as in the preceding proposition,
letting the distance $\dist$ be the $\l^\infty$ metric.  Let $c$ be
any fixed number greater than $d\rho/2\eta$.  Then, for any
$t > c/d\rho$, if we let $r = c/d\rho t$, then $r < 1$ and
$1/2t < r\eta$.  If we define $\vw_i$ and~$\bar L$ as above,
and let $\vv'_i = \vw_i/t$ and $L' = \bar L/t$, then
$\dist(\vv_i,\vv'_i) \le 1/2t$ for each $i$, so we can conclude that
$S^+ + L'$ covers~$\R^d$, where $S^+$ is the set of points within
distance $r\rho$ of $\Rocta dk$.  It is easy to see that
$S^+ \subseteq \Rocta d{k+dr\rho}$, so $\Rocta d{k+dr\rho} + L'$
covers~$\R^d$.  Hence, $\Rocta d{tk+tdr\rho} + tL'$
covers~$\R^d$; but $tdr\rho = c$ and $tL'=\bar L$, so we are done.

The proof for $\Rtetra dk$ is almost the same.  Let $S = \Rtetra dk$
and apply the preceding proposition (using the $l^\infty$ metric
again) to get $\eta$ and~$\rho$.  Fix $c > d\rho/\eta$.
For any $t > c/2d\rho$, if we let $r = c/2d\rho t$, then
$r < 1$ and $1/2t < r\eta$.  Now define $\vw_i$,~$\bar L$, $\vv'_i$,
$L'$, and~$S^+$ as above, and conclude again that $S^+ + L'$ covers~$\R^d$.
One can easily check that $S^+ \subseteq \Rtetra d{k+2dr\rho} - r\rho\vu$,
where $\vu = (1,\dots,1)$.  Hence,
$\Rtetra d{k+2dr\rho} - r\rho\vu + L'$ covers~$\R^d$,
so $\Rtetra d{k+2dr\rho} + L'$ covers $\R^d + r\rho\vu = \R^d$.
Now multiply by~$t$ to see that $\Rtetra d{tk+c} + \bar L$ covers~$\R^d$.
\QED

Again, it is easy to modify this proof to work for other sets
in place of $\Rocta dk$ or $\Rtetra dk$.  Also, the proof is quite
effective, allowing one to compute specific values of $c$ and $t$
which work for a given lattice~$L$ (assuming it is feasible to
compute $M$ and $\rho$).

The covering $\Rocta d{tk} + tL$ has the same efficiency as the
covering $\Rocta dk + L$; since $\Rocta d{tk+c} + \bar L$ is a
relatively slight perturbation of $\Rocta d{tk} + tL$ when $t$~is
large, it has almost the same efficiency.  Since $\bar L$~is an integer
lattice, the fact that $\Rocta d{tk+c} + \bar L$ covers~$\R^d$ implies
that $\octa d{\lfloor tk+c\rfloor} + \bar L$ covers~$\Z^d$; again the
efficiency is almost the same if $t$~is large.  Therefore, we can
construct integer lattice coverings as nearly efficient as desired to a
given real lattice covering, thus giving asymptotic results for
the present case of the degree-diameter problem.  The precise result
is as follows.

\proclaim{Theorem 9} Let $\eps_\R$ be the best possible efficiency
for a lattice covering of $\R^d$ by $\Rocta d1$, and let $\eps_\Z(k)$
be the best possible efficiency for a lattice covering of $\Z^d$
by $\octa dk$.  Then $\eps_\Z(k) = \eps_\R + O(k^{-1})$. The same
applies to $\Rtetra d1$ and~$\tetra dk$.\endproclaim

\demo{Proof}  Let $L$ be a lattice giving a lattice covering of $\R^d$
by $\Rocta d1$ with efficiency $\eps_\R$.  Let $\vv_1,\dots,\vv_d$ and
$c$ be as in Theorem~8.  Given a large integer $k$, let $t = k - c$,
and let $\bar L$ be the integer lattice approximating $tL$ as in
Theorem~8, generated by $\vw_1,\dots,\vw_d$.  Since $\bar L$ is an
integer lattice and $\Rocta d{t+c} + \bar L = \R^d$, we have $\octa dk
+ \bar L = \Z^d$.  Let $M$ and $\bar M$ be the $d \times d$ matrices
whose rows are $\vv_i$ and $\vw_i$, respectively; then $\det M =
(2^d/d!)\eps_\R$ and $\det {\bar M} \le |\octa dk|\eps_\Z(k)$, which
implies $\det(k^{-1}\bar M) \le (2^d/d! + O(k^{-1}))\eps_\Z(k)$.  But
$\bar M = tM + O(1) = kM + O(1)$, so $k^{-1}\bar M = M + O(k^{-1})$, so
$\det(k^{-1}\bar M) = \det M + O(k^{-1})$; this implies $\eps_\Z(k) \ge
\eps_\R + O(k^{-1})$.

On the other hand, Proposition 4(a) states that any lattice that
gives a covering of~$\Z^d$ by~$\octa dk$ also gives a covering of~$\R^d$
by~$\Rocta d{k+d/2}$.  If the efficiency of the former is~$\eps_\Z(k)$,
then the efficiency of the latter is $(|\octa
dk|/\volume(\Rocta d{k+d/2}))\eps_\Z(k) = (1+O(k^{-1}))\eps_\Z(k)$;
hence, $(1+O(k^{-1}))\eps_\Z(k) \le \eps_\R$, so $\eps_\Z(k) \le
\eps_\R + O(k^{-1})$.

The same argument works for $\Rtetra d1$ and~$\tetra dk$. \QED

Again, the same applies to other shapes as well.  Combining this
with the known sizes of the sets $\octa dk$ and $\tetra dk$ gives:

\proclaim{Corollary 10} Let $d$ be a fixed positive integer.

(a) If $\eps_\R$ is the best possible
efficiency for a lattice covering of~$\R^d$ by~$\Rocta d1$, then the
size of the largest possible undirected Cayley graph of an Abelian
group on $d$~generators with diameter at most~$k$ is\/
$(2^d\eps_\R/d!) k^d + O(k^{d-1})$.

(b) If $\eps'_\R$ is the best possible
efficiency for a lattice covering of $\R^d$ by $\Rtetra d1$, then the
size of the largest possible directed Cayley graph of an Abelian
group on $d$ generators with diameter at most $k$ is\/
$(\eps'_\R/d!) k^d + O(k^{d-1})$. \QNED

The above assumes that the upper limit $\eps_\R$ on the efficiency of a
lattice covering of~$\R^d$ by~$\Rocta d1$ is actually attained.  To see
that this is the case, first note that there is certainly a lattice
covering with positive efficiency, say~$\eps_0$.  Now consider the
possible sets of generating vectors $\vv_1,\dots,\vv_d$ for a lattice~$L$
giving a covering of efficiency at least~$\eps_0$.  We may assume
that, if $\vv_i$~is one of the generators and $\vw$~is an integral
linear combination of the other generators, then $|\vv_i| \le
|\vv_i+\vw|$; otherwise, just replace~$\vv_i$ with $\vv_i+\vw$ to get a
smaller set of generating vectors for~$L$, and iterate until no more
such reductions are possible.  It follows that, if, say, $\vv_1$~is the
longest of the vectors~$\vv_i$, then the angle between $\vv_1$ and the
hyperplane~$P$ spanned by $\vv_2,\dots,\vv_d$ is bounded below by a
positive number ($\pi/3$ for $d=2$, somewhat less for higher~$d$).  But
the distance from $\vv_1$ to $P$ is at most the diameter of~$\Rocta d1$
(i.e.,~$2$), since otherwise $\Rocta d1 + L$ would consist of
`hyperplanes' of copies of $\Rocta d1$ with gaps in between.  Putting
these together, we get a fixed upper bound~$B$ on the length~$|\vv_1|$,
and hence on all of the lengths~$\vv_i$.  But we also have a positive
lower bound~$b$ on the determinant $\det(\vv_1,\dots,\vv_d)$, namely
$\eps_0\volume(\Rocta d1)$.  One can now show that there is a fixed
number~$M$ such that, if $m_1,\dots,m_d \in \Z$ and $|m_1|+\dots+|m_d|
> M$, then $|m_1\vv_1+\dots+m_d\vv_d| > 3$.  Hence, the finitely many
translates $\Rocta d1 +m_1\vv_1+\dots+m_d\vv_d$ with $|m_1|+\dots+|m_d|
\le M$ will have to cover all points at distance~$\le 2$ from the
origin.  Now, the set of sequences $\vv_1,\dots,\vv_d$ with all
$|\vv_i| \le B$ such that the translates $\Rocta d1 +
m_1\vv_1+\dots+m_d\vv_d$ for $|m_1|+\dots+|m_d| \le M$ cover all
points at distance $\le 2$ from~$\vzero$ is a compact set, so there is
a sequence $\vv_1,\dots,\vv_d$ in this set for which
$\det(\vv_1,\dots,\vv_d)$ is maximal; this sequence of vectors
generates a lattice covering of~$\R^d$ by~$\Rocta d1$ (by
Proposition~6) with maximal possible efficiency.  The same argument
works for any compact shape of positive volume, such as~$\Rtetra d1$.
(Presumably this argument is well known, but the authors were not
able to find a reference for it.)

As we will see in a later section (for the case $d=3$, but the argument
is general), the determination of actual values for $\eps_\R$ and~$\eps'_\R$
for a particular dimension~$d$
can in principle be reduced to the solution of a finite number of
nonlinear optimization problems (each of which requires maximizing
a degree-\snug$d$ polynomial function over a region which is
a convex polytope in~$\R^{d^2}$).  Unfortunately, this finite number is
extremely large, so the actual values are not known for $d > 2$.
(For $d=1$ we trivially have $\eps_\R = \eps'_\R = 1$.  For
$d=2$ we will see below that $\eps_\R = 1$ and $\eps'_\R = 2/3$.)
The computations described in later sections lead us to conjecture that,
for $d=3$, $\eps_\R = 8/9$ and $\eps'_\R = 63/125$.

\subhead Two generators, undirected \endsubhead
We now begin to consider the results that can be obtained for specific
values of $d$.  As noted previously, the case $d=1$ is trivial,
so we will start with the case of undirected Abelian Cayley graphs
on two generators and a given bound~$k$ on the diameter.
The results in this subsection are not new.

As noted before, the diameter of the ordinary toroidal mesh
$\Z_m \times \Z_n$ is $\lfloor m/2\rfloor + \lfloor n/2\rfloor$.
Hence, the largest such mesh with diameter $\le k$ is the one where
$m$ is $k$ rounded up to the nearest odd integer and $n$ is
$k+1$ rounded up to the nearest odd integer.  This corresponds to
the lattice covering of $\Z^2$ by $\octa 2k$ using the lattice
$L = m\Z \times n\Z$; the efficiency of this covering is
$mn/(2k^2+2k+1)$, which tends to $1/2$ as $k$ grows large.
So one can hope for better results.

To get these, consider the real rotated square $\Rocta 2k$.  There
is obviously a lattice tiling using this square; it is just a
rotated orthogonal grid with spacing ${\sqrt 2}k$ between lines.
The lattice $L_1$ for this tiling is generated by the vectors
$(k,k)$ and $(-k,k)$; as expected, the corresponding determinant
is~$2k^2$, which is equal to the area of $\Rocta 2k$.

It now follows from the approximation results that we can get lattice
coverings of $\Z^2$ using $\octa 2k$ with efficiencies that approach
$1$ for large $k$.  However, we do not need the general results
here; since $L_1$ is already an integer lattice, we can simply
note that $\Rocta 2k+L_1 = \R^2$ implies $\octa 2k+L_1 = \Z^d$.
So this gives a lattice covering using $\octa 2k$ whose index
is $\indx \Z^2:L_1| = 2k^2$, which is better than that from the
best toroidal mesh for all $k \ge 3$; for large $k$, the
efficiency approaches $1$.

The corresponding Cayley graph $\Z^2/L_1$ turns out to be quite simple
to describe.  The $2k\times k$ rectangle $\{1,\dots,2k\} \times \{
1,\dots,k\}$ contains exactly one point from each coset of~$L_1$,
so it can serve as a set of vertices for the graph.  Adjacent points in
the rectangle (horizontally and vertically) are connected as in the
usual mesh.  Horizontally, one has the usual toroidal connections at
the ends: $(1,j)$ is connected to $(2k,j)$.  But vertically,
there is an offset of~$k$: $(i,1)$ is connected to $(i+k,k)$ if
$i \le k$, or to $(i-k,k)$ if $i > k$.  This is just like a
$2k\times k$ toroidal mesh, except that the torus is twisted
halfway around before the long edges are glued together.  This twist
allows one to double the number of vertices in a $k\times k$
toroidal mesh while increasing the diameter by at most $1$ (there
is no increase if $k$~is even).

A number of the useful properties of ordinary toroidal meshes apply
with very little change to twisted toroidal meshes.  For instance,
since the new mesh is still just a rectangular mesh with extra
connections at the boundary, it is easy to map a simple rectangular
grid into the mesh, by simply ignoring the boundary connections.

Another nice property of toroidal meshes is that it is easy to
find a shortest route from one node to another: just check for each
coordinate separately which of the two possible directions gives
a shorter path, and put the results together.  Finding optimal routes
is only slightly more complicated for the twisted toroidal mesh.
To see this, consider the given $2k\times k$ rectangle as half of a
$2k\times 2k$ rectangle; each node $(i,j)$ in the first half has
a copy $(i\pm k,j\pm k)$ in the other half.  This larger rectangle
is then copied periodically without further twists to cover $\Z^2$;
in other words, the $2k\times k$ twisted toroidal mesh is just
a $2k\times 2k$ toroidal mesh where $(i,j)$ and $(i\pm k,j\pm k)$
are identified as a single node.  Therefore, to find an optimal route
from $(i,j)$ to $(i',j')$ in the twisted mesh, apply the ordinary
$2k\times 2k$ toroidal mesh routing algorithm to find optimal
routes from $(i,j)$ to $(i',j')$ and to $(i'\pm k,j'\pm k)$,
and choose the shorter of the two.

\medskip
In the real case the lattice $L_1$ gave a perfect tiling of $\R^2$
using $\Rocta 2k$, since boundary overlap did not count; but
in the integer case the boundary overlap reduces the efficiency
slightly, from $1$ to $2k^2/(2k^2+2k+1)$.  It turns out that
if one uses a slightly modified lattice, namely the lattice $L_2$
with generating vectors $(k,k+1)$ and $(-k-1,k)$, then one gets
a covering of $\Z^2$ by copies of $\octa 2k$ with efficiency $1$
(i.e., a tiling).  See Figure 1.  We therefore get:

\proclaim{Theorem 11 \rm (multiple authors)}
The largest possible size for the undirected
Cayley graph of an Abelian group on two generators having diameter $k$
is $2k^2+2k+1$. \QNED

This result has appeared in various forms in a number of places
(usually stated so as to apply only to cyclic Cayley graphs,
but since the optimal Abelian Cayley graphs turn out to be cyclic,
the results are basically equivalent).  See, for instance,
Boesch and Wang~\cite{\BoeWan} or Yebra el~al.~\cite{\YebFioMorAle};
the Bermond-Comellas-Hsu survey~\cite{\BerComHsu} has many additional
references.

\midinsert
\centerline{
\plotfigbegin
   \plotfigscale 18
   \plotvskip 17
   \plottext{}
   \plot 23 17
   \def\1{\plotrline 0 1 }
   \def\2{\plotrline 0 -1 }
   \def\3{\plotrline 1 0 }
   \plotmove 0 2 \2\3\2\3\3
   \plotmove 0 8 \2\3\3\2\3\2\3\2\3\3\2\3\2\3\2\3\3\2\3
   \plotmove 0 14 \3\2\3\2\3\2\3\3\2\3\2\3\2\3\3\2\3\2\3\2\3\3\2\3\2\3\2\3%
         \3\2\3\2\3
   \plotmove 4 17 \3\2\3\2\3\3\2\3\2\3\2\3\3\2\3\2\3\2\3\3\2\3\2\3\2\3\3\2%
         \3\2\3\2\3
   \plotmove 12 17 \3\2\3\3\2\3\2\3\2\3\3\2\3\2\3\2\3\3\2
   \plotmove 20 17 \3\3\2\3\2
   \plotmove 0 13 \1\1\3\1\3\1\3
   \plotmove 0 5 \1\3\1\1\3\1\3\1\3\1\1\3\1\3\1\3\1\1\3\1\3
   \plotmove 2 0 \1\3\1\3\1\3\1\1\3\1\3\1\3\1\1\3\1\3\1\3\1\1\3\1\3\1\3\1\1\3
   \plotmove 8 0 \3\1\1\3\1\3\1\3\1\1\3\1\3\1\3\1\1\3\1\3\1\3\1\1\3\1\3\1\3\1
   \plotmove 14 0 \3\1\3\1\1\3\1\3\1\3\1\1\3\1\3\1\3\1\1\3\1
   \plotmove 20 0 \3\1\3\1\3\1\1
   \plotclosed
   \dimen0=1sp
   \plotcurrx=0.5\plotunitx
   \loop
      {\dimen2=1sp
      \plotcurry=0.5\plotunity
      \loop
         \plothere
         \advance\plotcurry by\plotunity
      \ifdim\dimen2<17sp
         \advance\dimen2 by 1sp
      \repeat}
      \advance\plotcurrx by\plotunitx
   \ifdim\dimen0<23sp
      \advance\dimen0 by 1sp
   \repeat
   \plotlinewidth=0.1truept
   \plotmove 11.2 8.8
   \plotline 17.8 8.8
   \plotline 17.8 6.2
   \plotline 14.8 6.2
   \plotline 14.8 5.2
   \plotline 11.2 5.2
   \plotline 11.2 8.8
\plotfigend}
\botcaption{Figure 1}
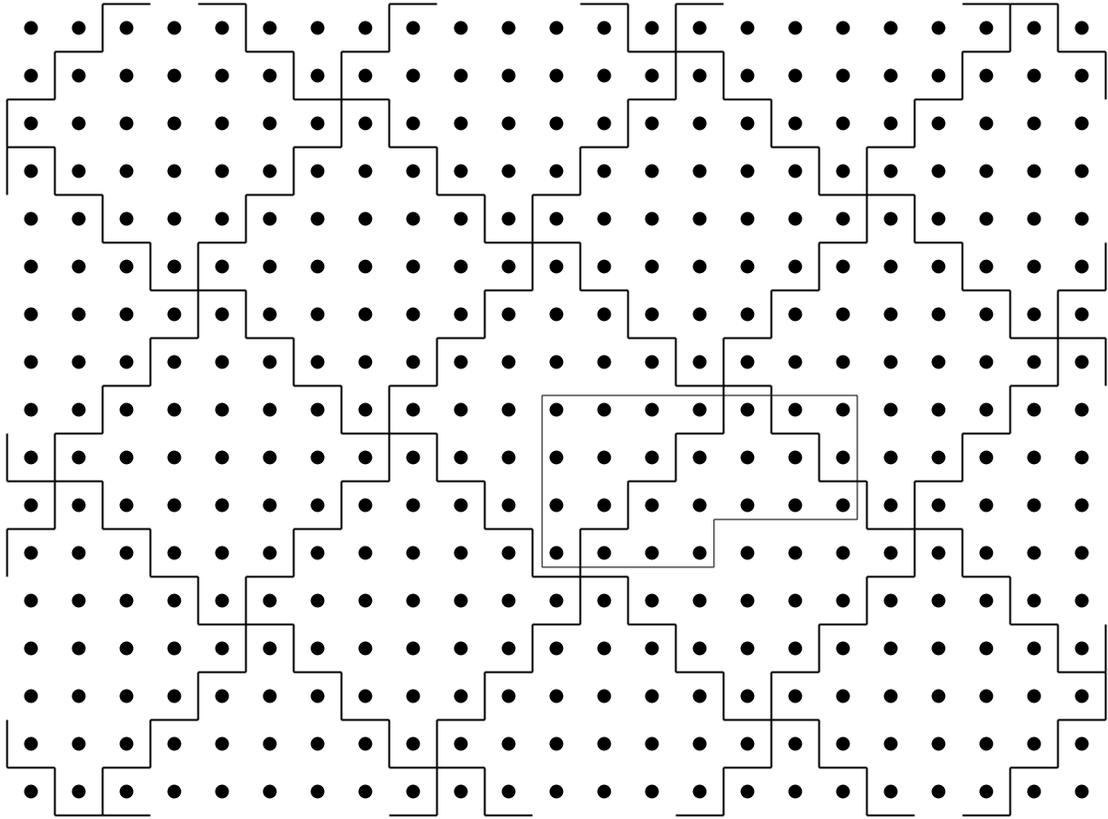
Lattice tiling of $\Z^2$ using $\octa 2k$ (shown for $k = 3$).
\endcaption
\endinsert

The tiling in Figure~1 is in Yebra el~al.~\cite{\YebFioMorAle},
among other places; it even appears in native
artwork of the southwestern United States, and may date back to
the ancient Aztecs, who did use the stepped diamond shape in
temple ornamentation.  (This shape is now commonly known as the Aztec
diamond, a term coined by J.~Propp.) 
However, it is unlikely that the Aztecs were motivated by the desire to
construct efficient parallel computation networks.

It is easy to see that the lattice tiling of $\Z^2$ by $\octa 2k$,
or of $\R^2$ by the Aztec diamond, is unique except for a possible
reflection about the line $x=y$; this just corresponds to interchanging
the two generators for the Cayley graph.  Therefore, the Cayley
graph attaining the bound in Theorem~11 is unique up to isomorphism.

Since the point $(2k+1,1)$ is in~$L_2$, we have $\ve_2 + L_2 =
(-2k-1)(\ve_1+L_2)$ in $\Z^2/L_2$.  Hence, $\Z^2/L_2$ is a cyclic group,
generated by $\ve_1+L_2$ alone.  It is isomorphic (not only as a group but
as a Cayley graph) to $\Z_{2k^2+2k+1}$ with the generating set
$\{1,2k^2\}$.  One may choose to replace the second generator by
its inverse, making the generating set $\{1,2k+1\}$; other generating
sets can be used as well.

For layout purposes, one may just arrange the nodes in the form of
the diamond $\octa 2k$ and connect the boundary nodes as specified
by $L_2$, but it is probably more convenient to use the
almost-rectangular shape outlined in Figure~1 (a $(2k+1)\times k$
rectangle with an extra partial row of length $k+1$).  The
boundary connections are similar to those for the twisted toroidal
mesh given before, but now there is also a slight twist when connecting
the short sides; there is a drop of one row when wrapping around from
right to left.  This layout shows that one can embed a rectangular
grid into this graph so as to use almost all of the nodes.

\subhead Two generators, directed \endsubhead
We now describe the largest possible directed Cayley graph of an Abelian
group on two generators with diameter bounded by $k$.
As in the preceding subsection, the two-generator results here are
already known.

The best
toroidal mesh in this case is $\Z_m \times \Z_{m'}$, where $m = \lfloor k/2
\rfloor + 1$ and $n = \lceil k/2 \rceil + 1$; this gives size
$mm' = \lfloor (k+2)^2/4 \rfloor$, which is about $1/2$ of
$|\tetra 2k|$, so one can hope to do better.

However, one is not going to get perfect efficiency in this case.
One can easily tile the plane with triangles such as $\Rtetra 2k$
if one is allowed to rotate them, but this is not possible using
only a lattice of translated copies of a triangle.  The exact
minimum density for a lattice covering of the plane by triangles
was computed by I.~F\'ary in 1950; we will give a different proof
of his result here, and then give the analogue for $\Z^2$.

\proclaim{Theorem 12 \rm(F\'ary \cite{\Fary})}  The minimum density
for a lattice covering of\/~$\R^2$ by triangles is\/ $3/2$.  Equivalently,
the maximum efficiency is\/~$2/3$. \endproclaim

\define\1{{\overline{AB}}}
\define\2{{\text{Area}}}
\demo{Proof} Since the density and efficiency of a lattice covering
are invariant under affine transformations, it does not matter
which triangle we work with, so, for slight convenience, let us
work with the isosceles right triangle $\Rtetra 21$.

One can attain the efficiency $2/3$ by using the lattice with
generating vectors $(1/3,1/3)$ and $(2/3,-1/3)$.  This corresponds
to a tiling of the plane using an L-tromino that takes up $2/3$ of
$\Rtetra 21$, as shown in Figure~2.  Or one can cut off all three
corners of the triangle to get a hexagon that tiles the plane.

\midinsert
\centerline{
\plotfigbegin
   \plotfigscale 36
   \plotvskip 3
   \plottext{}
   \plot 3 3
   \plotmove 0 0
   \plotline 2 0
   \plotline 2 1
   \plotline 1 1
   \plotline 1 2
   \plotline 0 2
   \plotline 0 0
   \plotlinewidth=0.1truept
   \plotmove 0 0
   \plotline 3 0
   \plotline 0 3
   \plotline 0 0
\plotfigend}
\botcaption{Figure 2}
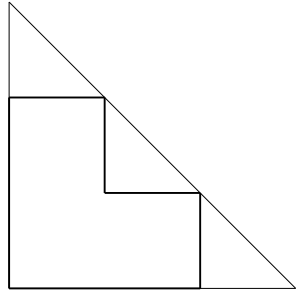
A subset of the triangle $\Rtetra 21$ which tiles the plane.
\endcaption
\endinsert

Now suppose that we have a lattice covering of $\R^2$ using $\Rtetra 21$
and the lattice $L$; we must show that the efficiency of the covering
is at most $2/3$.  Let $\prec$ be a linear ordering of $L$ compatible
with addition which is defined by primarily ordering points $(x,y)$
according to the sum $x+y$ (so points that are farther out in the
direction $(1,1)$ come later in the ordering) and breaking ties
(if any) by distance in some other direction.

Let $\1$ be the hypotenuse of $\Rtetra 21$.  Only finitely many of the
$L$\snug-translates of~$\Rtetra 21$ lie near~$\Rtetra 21$; of these,
the ones of the form $\Rtetra 21 + \vv$ for $\vv \succ \vzero$ must
cover the points which are near~$\1$ on the side away from~$\Rtetra
21$.  Since the union of finitely many translates of~$\Rtetra21$ is
closed, $\1$~itself is covered by finitely many translates $\Rtetra 21
+ \vv$ with $\vv \succ \vzero$.  Find such a covering of~$\1$ with as
few translates as possible, say $\Rtetra 21 + \vv_1,\dots,\Rtetra 21 +
\vv_m$, where $\vv_i \succ \vzero$.

Note that, since $\Rtetra 21+\vv_i$ must intersect $\1$, it contains
one of the endpoints $A$ and $B$ if and only if the coordinates of
$\vv_i$ are not both positive.  We may assume that at most one of the
vectors $\vv_i$ has both coordinates positive.  For if there are two
such, let them be $\Rtetra 21+\vv_i$ and $\Rtetra 21+\vv_j$ where
$\vv_i \prec \vv_j$.  Then $\vv_j - \vv_i \succ \vzero$; since $\vv_i$
has both coordinates positive, we have $$(\Rtetra 21 + \vv_j) \cap \1
\subseteq (\Rtetra 21 + \vv_j - \vv_i) \cap \1.$$ Furthermore, $\vv_j -
\vv_i$ cannot have both coordinates positive, because, if it did, we
would have $$(\Rtetra 21 + \vv_j) \cap \1 \subseteq (\Rtetra 21 +
\vv_i) \cap \1,$$ so $\Rtetra 21 + \vv_j$ would not have been needed in
the covering of $\1$, contradicting the minimality of $m$.  So we can
replace $\Rtetra 21 + \vv_j$ with $\Rtetra 21 + \vv_j - \vv_i$ to get
another covering of $\1$ using fewer vectors with both coordinates
positive.  Repeat this until only one such vector is left.

Since each translate $\Rtetra 21 + \vv_i$ is convex, its
intersection with $\1$ is a segment or a point.  Therefore, at most
one of these translates can contain $A$, since otherwise one of the
intersections $(\Rtetra 21 + \vv_i) \cap \1$ would include
another such intersection, making the latter translate superfluous
and contradicting the minimality of $m$.  Similarly, at most one
of the translates $\Rtetra 21 + \vv_i$ contains $B$.  Putting these
facts together, we conclude that we needed at most three of the
translates $\Rtetra 21 + \vv$ with $\vv \succ \vzero$ to cover
the segment $\1$.

In other words, there are points $P$ and $Q$ on $\1$ such that each of
the three segments $\overline{AP}$,~$\overline{PQ}$, and
$\overline{QB}$ is covered by one of the translates $\Rtetra 21 + \vv$
with $\vv \succ \vzero$.  Let $l_1,l_2,l_3$ be the lengths of these
three segments; then $l_1+l_2+l_3 = \sqrt{2}$.  Note that, if $\Rtetra
21 + \vv$ covers~$\overline{AP}$, then $\Rtetra 21 + \vv$ covers the
entire isosceles right triangle below $\overline{AP}$ whose hypotenuse
is~$\overline{AP}$; the area of this triangle is $l_1^2/4$.  Similar
statements hold for $\overline{PQ}$ and $\overline{QB}$.  So we have
three disjoint triangles included in $\Rtetra 21$ which are covered by
translates $\Rtetra 21 + \vv$ with $\vv \succ \vzero$, and the total
area of these triangles is $(l_1^2+l_2^2+l_3^2)/4$.

By the proof of (the real version of) Lemma~3, if we let $T$ be
the part of $\Rtetra 21$ which is not covered by any translate
$\Rtetra 21+\vv$ with $\vv\succ \vzero$, then $T$ gives a lattice
tiling of $\R^2$ using $L$, so the efficiency of the covering
using $\Rtetra 21$ and $L$ is $\2(T)/\2(\Rtetra21)$.
We have $\2(\Rtetra21) = 1/2$ and
$\2(T) \le 1/2 - (l_1^2+l_2^2+l_3^2)/4$.  A standard minimization
shows that, if $l_1+l_2+l_3 = \sqrt{2}$, then
$l_1^2+l_2^2+l_3^2 \ge 2/3$ (with equality only when
$l_1=l_2=l_3=\sqrt{2}/3$).  Therefore,
$\2(T) \le 1/3$, so the efficiency of the covering by $\Rtetra21$ and
$L$ is at most $2/3$. \QED
\undefine\1
\undefine\2

It now follows from Theorem~9 that the largest possible
index $\indx \Z^2:L|$ for an integer lattice $L$ giving a lattice covering
of $\Z^2$ by $\tetra 2k$ is approximately $(2/3)|\tetra 2k|$, or about
$k^2/3$, for large $k$.  However, we can actually get an exact answer
rather than an approximation.

\proclaim{Theorem 13 \rm (mainly Wong and Coppersmith~\cite{\WonCop})}
The largest possible index $\indx \Z^2:L|$ for
a lattice $L$ giving a lattice covering of $\Z^2$ by $\tetra 2k$ is
$\lfloor (k+2)^2/3 \rfloor$. \endproclaim

\define\1{{\overline{AB}}}
\demo{Proof} We will give a discrete form of the proof of Theorem~12.
Let $a$ be $(k+2)/3$ rounded to the nearest integer, and
let $b = k + 2 - 2a$ (so $b$ is also about $(k+2)/3$).  Let $T_k$ be the
set of $(i,j)$ in $\Z^2$ such that $i,j \ge 0$, $\min(i,j) < a$, and
$\max(i,j) < a+b$.  Then $T_k \subseteq \tetra 2k$, since any
$(i,j)$ in $T_k$ satisfies $i+j \le a-1+a+b-1=k$.  The set $T_k$
looks like the L-tromino from Figure~2, and it tiles $\Z^2$
using the lattice with generating vectors $(a,a)$ and
$(a+b,-b)$.  Therefore, this lattice gives a covering of
$\Z^2$ using $\tetra 2k$, and its index is $a(a+2b)$, which works
out to be $\lfloor (k+2)^2/3 \rfloor$.

Now, suppose we have a lattice covering of $\Z^2$ using $\tetra 2k$
and a lattice $L$; we must show that $\indx \Z^2:L| \le \lfloor
(k+2)^2/3 \rfloor$.  Define the linear order $\prec$ of $L$ as
before.  Let $A$ and $B$ be the points $(0,k+1)$ and $(k+1,0)$;
then the segment $\1$ contains $k+2$ integral points, which must
be covered by translates $\tetra 2k + \vv$ where $\vv \in L$ and
$\vv \succ \vzero$.

Let $\vv_1,\dots,\vv_m$ be a list of as few vectors as possible in
$L$ such that $\vv_i \succ \vzero$ and the translates
$\tetra 2{k+1}+\vv_i$ of $\tetra 2{k+1}$ cover all of the
integral points on~$\1$.  Then the same argument as for Theorem~12
shows that $m$ is at most $3$.  Hence, $\1$ can be
broken up into three segments $\overline{AP}$, $\overline{P'Q}$,
and $\overline{Q'B}$ (where $P$ and $P'$ are adjacent integral points
on~$\1$, as are $Q$ and $Q'$), each of whose integral points are
covered by one of the translates $\tetra 2{k+1}+\vv$ with $\vv \succ
\vzero$.  Let $l_1,l_2,l_3$ be the numbers of integral points
on these segments; then $l_1+l_2+l_3 = k+2$.

If $\tetra 2{k+1}+\vv$ covers the integral points on $\overline{AP}$,
then it covers all of the integral points in the isosceles right
triangle below $\overline{AP}$ and having $\overline{AP}$ as its hypotenuse.
In fact, all of these points other than those on $\overline{AP}$ itself
are covered by $\tetra 2k+\vv$; there are $(l_1^2-l_1)/2$ such
points, and they are all in $\tetra 2k$.  Similarly, the segments
$\overline{P'Q}$ and $\overline{Q'B}$ give
$(l_2^2-l_2)/2 + (l_3^2-l_3)/2$ more points of $\tetra 2k$ which
are covered by translates $\tetra 2k + \vv$ with $\vv \succ \vzero$.

As in the proof of Lemma~3, let $T$ be the set of points in
$\tetra 2k$ that are not in $\tetra 2k + \vv$ for any $\vv \succ \vzero$;
then $|T| = \indx Z^2:L|$.  The calculations above show that
$$|T| \le |\tetra 2k| + (l_1+l_2+l_3)/2 - (l_1^2 + l_2^2 + l_3^2)/2.$$
Here $|\tetra 2k| = (k+1)(k+2)/2$, and $l_1+l_2+l_3$ is just $k+2$.
Given $l_1+l_2+l_3$, we minimize $l_1^2 + l_2^2 + l_3^2$ by making
the numbers $l_1,l_2,l_3$ as close to equal as possible; in this
case, this means that the minimum occurs when two of them are $a$
and the third is $b$.  Therefore,
$$|T| \le {(k+1)(k+2) \over 2} + {k+2 \over 2} - {2a^2 + b^2 \over 2},$$
which simplifies to $|T| \le \lfloor (k+2)^2 /3 \rfloor$, as desired. \QED
\undefine\1

\proclaim{Corollary 14 \rm (mainly Wong and Coppersmith~\cite{\WonCop})}
The largest possible size for the directed
Cayley graph of an Abelian group on two generators having diameter~$k$
is $\lfloor (k+\nobreak 2)^2 /3 \rfloor$. \QNED

Again it is hard to be historically accurate here, because different
authors have presented results in quite different ways; see
the Bermond-Comellas-Hsu survey~\cite{\BerComHsu}
for more information and references.

Let $a$, $b$, and $T_k$ be as in the proof of Theorem~13; then
the set $T_k$ gives a suitable layout for a network realizing
this Cayley graph.  In addition to the mesh connections from $(i,j)$
to $(i+1,j)$ and $(i,j+1)$ within $T_k$, one will also need
wraparound connections from $(a+b-1,j)$ to $(0,j+b)$ for $j < a$,
from $(a-1,j+a)$ to $(0,j)$ for $j<b$, from $(i,a+b-1)$ to $(i+b,0)$
for $i < a$, and from $(i+a,a-1)$ to $(i,0)$ for $i<b$.

In the case $a=b$, one can use an alternate layout in the form
of a $3a\times a$ rectangle, with wraparound connections from
$(3a,j)$ to $(1,j)$ and from $(i,a)$ to $((i+a)\bmod a,1)$.  This
is just a variant of the twisted toroidal mesh where the long dimension
is twisted by a factor of $1/3$ rather than $1/2$; it is
convenient for construction and for embedding a rectangular
grid without boundary connections into the network (although this
is not particularly useful in the directed case).  If $a\ne b$, then
one gets a rectangle with some missing nodes or extra nodes along
part of one edge, and the cross-connections are slightly messier.

If $k \equiv 1 \pmod 3$, so that $a = b = (k+2)/3$, then one can see
from the proof of Theorem~13 that the lattice $L$ with generating
vectors $(a,a)$ and $(a+b,-b)$ is the unique lattice attaining the
bound in the theorem, and hence the Cayley graph attaining the bound in
Corollary~14 is also unique.  However, if $k \not\equiv 1 \pmod 3$,
so that $a$ and $b$ differ by $1$, then there are two more lattices
attaining the bound: the lattice $\tilde L$ with generating vectors
$(a,b)$ and $(2a,-a)$, and the mirror image with generating vectors
$(b,a)$ and $(-a,2a)$.  The latter two give Cayley graphs that are
isomorphic to each other, but not to the Cayley graph of $\Z^2/L$ (if
$k > 1$), because the Cayley graph of $\Z^2/L$ has cycles of length
$2a$ while that of $\Z^2/\tilde L$ does not.  Therefore, if $k > 0$ and
$k \not\equiv 1 \pmod 3$, then there are exactly two Cayley graphs
meeting the bound of Corollary~14.

If $k > 1$ and $k \equiv 1 \pmod 3$, then the optimal group~$\Z^2/L$ is not
cyclic; it is isomorphic to $\Z_{3a} \times \Z_a$ by an isomorphism sending
$\ve_1$ and $\ve_2$ to $(1,0)$ and $(3a-1,1)$.  On the other hand,
if $k \not\equiv 1 \pmod 3$, then $(2a+b,a-b)$ is in~$L$ and
$a - b = \pm 1$, so $\ve_2$ is a multiple of~$\ve_1$ in~$\Z^2/L$ and
hence $\Z^2/L$~is cyclic.  Similarly, $\Z^2/\tilde L$ is cyclic,
since $(3a,b-a) \in \tilde L$.  One can get the corresponding
Cayley graphs directly from the cyclic group $\Z_{\lfloor (k+2)^2/3 \rfloor}$
by using the generator pairs $\{1,(2a+b)/(b-a)\}$ and
$\{1,3a/(a-b)\}$, respectively.

\subhead Three generators, undirected \endsubhead
For $d = 3$, we must consider three-dimensional lattice tilings by the
regular octahedron $\Rocta 3k$ and its discrete approximation $\octa 3k$.
These shapes do not tile space perfectly, and the best possible efficiency
for a lattice covering of $\R^3$ by $\Rocta 3k$ appears to be still
open (although there is a good guess, as we shall see).  So we will
apply our results in reverse, using computed results about the
degree-diameter problem to obtain information about lattice
tilings by octahedra.

The best three-dimensional toroidal mesh with diameter $k$ is
$\Z_{2b_0+1} \times \Z_{2b_1+1} \times \Z_{2b_2+1}$,
where $b_i = \lfloor (k+i)/3 \rfloor$; this has
about $(8/27)k^3$ vertices for large $k$.  This corresponds to the
covering of $\R^3$ by $\Rocta 31$ using the cubic lattice
$(2/3)\Z^3$; this covering has efficiency $2/9$.

It turns out that a good lattice to use for coverings with regular
octahedra is the body-centered cubic lattice, defined most simply as
the set $L_{\text{bcc}}$ of points in $\Z^3$ whose coordinates are all odd
or all even.  If $\vx$ is an arbitrary point of $\R^3$, then $\vx$
lies in or on one of the unit cubes with vertices in $\Z^3$; two
opposite corners of this cube will be in $L_{\text{bcc}}$, say
$\vv$ and $\vw$.  Then each coordinate of $\vx$ lies between (inclusively)
the corresponding coordinates of $\vv$ and~$\vw$, so, letting
$\dist$ be the $\l^1$ metric on $\R^3$, we have
$\dist(\vv,\vx) + \dist(\vx,\vw) = \dist(\vv,\vw) = 3$.
Hence, either $\dist(\vv,\vx) \le 3/2$ or $\dist(\vw,\vx) \le 3/2$.
This shows that $|\Rocta 3{3/2}| + L_{\text{bcc}} = \R^3$.
Now, $L_{\text{bcc}}$~has generators $(2,0,0)$, $(0,2,0)$, and $(1,1,1)$,
giving a matrix with determinant $4$, while the volume of
$\Rocta 3{3/2}$ is $9/2$, so this lattice covering of $\R^3$
has efficiency $8/9$.  A fundamental region for the lattice can be
obtained by truncating each of the corners of the octahedron,
giving an Archimedean solid whose faces are eight regular
hexagons and six squares.

The same reasoning shows that, for $k \ge 1$,
one can get a lattice covering of~$\Z^3$
by~$\octa 3k$ using the slightly distorted body-centered cubic lattice
$L_{\text{bcc}}(a_1,a_2,a_3)$
with generating vectors $(2a_1,0,0)$, $(0,2a_2,0)$, and $(a_1,a_2,a_3)$,
where $a_i = \lfloor (2k+i)/3 \rfloor$.  This gives a Cayley graph
of size $4a_1a_2a_3$, or approximately $(32/27)k^3$ for large~$k$.
This is an improvement over the best
toroidal mesh of diameter~$k$; it is about $4$~times as good
for large~$k$.

One can lay out the Cayley graph for $\Z^3/L_{\text{bcc}}(a_1,a_2,a_3)$
in the form of a $2a_1 \times 2a_2 \times a_3$ mesh.  Opposite
$2a_i \times a_3$ sides are connected to each other as in the
usual toroidal mesh, but the toroidal connections between
the top and bottom $2a_1 \times 2a_2$ sides are twisted in two
directions: node $(j_1,j_2,a_3)$ is connected to node
$(j_1 \pm a_1, j_2 \pm a_2, 1)$, where the signs are chosen
to give numbers between $1$ and $2a_i$, inclusive.
Routing algorithms and embeddings of rectangular grids work here just as
they did in the two-dimensional version.

Two questions now arise.  First, can one improve the efficiency
by making small adjustments to the discrete lattice, as we did
in the two-generator cases?  Second, can one get better results by
using a completely different lattice?  The answers to these
questions are not immediately clear, so we will approach the
problem from another direction.

One can write computer programs to examine various groups, choose
all (or at least many) possible sets of a certain number of generators
for the group, and compute the resulting diameters.  M.~Dinneen has
performed many such computations, some using exhaustive search of
generator sets and others using random sampling, on a number of
different kinds of groups, resulting in new best-known graphs
for the degree-diameter problem; see, for instance, Dinneen and
Hafner~\cite{\DinHaf}.  Some of Dinneen's earlier unpublished
computations were for Abelian (usually cyclic) groups of diameter
up to $10$ on various numbers of generators.

The authors have written a program to extend these calculations.
The program does an
exhaustive search of generating sets for each Abelian group,
but avoids examining many generating sets which give
Cayley graphs isomorphic to ones already examined; for instance,
in the case of a cyclic group~$\Z_n$,
one may assume that the first generator is a divisor of $n$.
Here `exhaustive search' means that all Abelian groups of size up to
$|\octa 3k| = (4k^3+6k^2+8k+3)/3$ were examined, so the results
definitely give the largest possible Cayley graph of an Abelian group
with diameter~$k$.  The program uses bit manipulations adapted from
(but simpler than) those of Dougherty and Janwa~\cite{\DouJan}, which
gave algorithms for diameter computations for Cayley graphs of Abelian
groups of exponent~$2$.

It turns out that, for each~$k$ for which the calculation has been
done so far (up to~$14$), the best Abelian Cayley graph has been obtained
from a cyclic group.
The results of the computation are shown in Table~1.
This extends (and corrects an erroneous final entry in) a similar
table given by Chen and Jia~\cite{\CheJia}.

\midinsert
\centerline{
\vbox{\offinterlineskip
\halign{
&\vrule#&\hfil$\,\,\vphantom{|\Rocta3k|}#\,\,$\hfil\cr
\noalign{\hrule}
height2pt&\omit&&\omit&&\omit&&\omit&&\omit&&\omit&&\omit&&\omit&\cr
&k&&|\octa 3k|&&\omit\hfil\,\,Toroidal\,\,\hfil&&
    \omit\hfil\,\,Twisted\,\,\hfil&&n_c&&\omit\hfil\,\,Generators\,\,\hfil&&
      n_c/|\octa 3k|&&n_c/\volume(\Rocta 3{k+3/2})&\cr
height2pt&\omit&&\omit&&\omit&&\omit&&\omit&&\omit&&\omit&&\omit&\cr
\noalign{\hrule}
height2pt&\omit&&\omit&&\omit&&\omit&&\omit&&\omit&&\omit&&\omit&\cr
&0&&1&&1&&&&1&&&&1&&.222222&\cr
&1&&7&&3&&4&&7&&1,2,3&&1&&.336000&\cr
&2&&25&&9&&16&&21&&1,2,8&&.840000&&.367347&\cr
&3&&63&&27&&48&&55&&1,5,21&&.873016&&.452675&\cr
&4&&129&&45&&108&&117&&1,16,22&&.906977&&.527423&\cr
&5&&231&&75&&192&&203&&1,7,57&&.878788&&.554392&\cr
&6&&377&&125&&320&&333&&1,9,73&&.883289&&.592000&\cr
&7&&575&&175&&500&&515&&1,46,56&&.895652&&.628944&\cr
&8&&833&&245&&720&&737&&1,11,133&&.884754&&.644700&\cr
&9&&1159&&343&&1008&&1027&&1,13,157&&.886109&&.665371&\cr
&10&&1561&&441&&1372&&1393&&1,92,106&&.892377&&.686940&\cr
&11&&2047&&567&&1792&&1815&&1,15,241&&.886663&&.696960&\cr
&12&&2625&&729&&2304&&2329&&1,17,273&&.887238&&.709953&\cr
&13&&3303&&891&&2916&&2943&&1,154,172&&.891008&&.724015&\cr
&14&&4089&&1089&&3600&&3629&&1,19,381&&.887503&&.730892&\cr
height2pt&\omit&&\omit&&\omit&&\omit&&\omit&&\omit&&\omit&&\omit&\cr
\noalign{\hrule}
}}
}
\botcaption{Table 1}
Best undirected Cayley graphs of cyclic groups, three generators.
\endcaption
\endinsert

The first column
is the desired diameter $k$.  The
second column gives the largest size one could hope for
of an undirected Cayley graph of an Abelian group on $3$~generators.
The next two columns give the sizes attained by
the best possible ordinary toroidal mesh and the twisted
toroidal mesh described above.  Next comes~$n_c$, the
computed largest~$n$ such that $\Z_n$ has three generators
giving it an undirected diameter of~$k$.  Then comes a triple of
generators of $\Z_{n_c}$ attaining this diameter (this is not always
unique, but only one generator set is given here).  The
final two columns gives the efficiencies of the corresponding
lattice coverings of~$\Z^3$ by~$\octa 3k$ and of~$\R^3$
by~$\Rocta 3{k+3/2}$ (see Proposition~4).

Some interesting observations can be made from Table~1.  First,
note that the twisted toroidal meshes do almost as well as the
optimal cyclic groups.  Also note that the numbers in the second-to-last
column
do seem to be getting close to $8/9$ for larger~$k$; this provides evidence
that the body-centered cubic lattice gives the best lattice covering
of $\R^3$ by $\octa 31$.

One can confirm this more strongly by reconstructing the lattices~$L$
for which $\Z^3/L$~gives these optimal cyclic groups.  For instance,
look at $k=10$, for which we have the cyclic group~$\Z_{1393}$ with
generating set $\{1,92,106\}$.  There is a unique homomorphism
from~$\Z^3$ to $\Z_{1393}$ which sends $\ve_1,\ve_2,\ve_3$ to $1,92,106$,
and the desired lattice $L$ is just the kernel of this homomorphism;
this means that
$$L = \{(x_1,x_2,x_3) \in \Z^3\colon\quad x_1 + 92x_2 + 106x_3 \equiv 0
\pmod{1393}\}.$$
One can easily find three vectors in~$L$, namely $(1393,0,0)$,
$(92,{-}1,0)$, and $(106,0,{-}1)$; the matrix with these three vectors
as rows has determinant $1393 = \indx Z^3:L|$, so these vectors generate $L$.
Now one can perform elementary operations to reduce these vectors
to a smaller set of generators for~$L$, such as $(7,7,7)$,
$(8,{-}7,6)$, and $(6,8,{-}7)$.  These vectors are quite close
to the vectors $(7,7,7)$, $(7,{-}7,7)$, and $(7,7,{-}7)$, which generate
a scaled-up body-centered cubic lattice (in fact, the latter lattice
gives the twisted toroidal mesh of size~$1372$ mentioned in the table).
Similarly, one finds that the other lattices corresponding to
the generators in Table~1 are almost body-centered cubic.

There are definite patterns in Table~1; every third~$k$ gives groups
and generators of the same form.  These patterns can be generalized,
giving the following result.

\proclaim{Theorem 15} For all $k \ge 0$, there is an undirected Cayley
graph on three generators of an Abelian (in fact, cyclic) group
which has diameter $k$ and size $n$, where
$$n = \cases (32k^3+48k^2+54k+27)/27 &\text{if $k \equiv 0 \pmod 3$,} \\
(32k^3+48k^2+78k+31)/27 &\text{if $k \equiv 1 \pmod 3$,} \\
(32k^3+48k^2+54k+11)/27 &\text{if $k \equiv 2 \pmod 3$.} \endcases$$
\endproclaim

\demo{Proof} We will show the existence of lattices $L_k \subseteq
\Z^3$ such that $\Z^3/L_k$ is cyclic, $\octa 3k + L_k = \Z^3$,
and $\indx \Z^3:L|$ is the $n$ specified in the theorem.

Let $a = \lceil 2k/3\rceil$.  For each $k$, we define $L_k$ by specifying
three generating vectors $\vv_1,\vv_2,\vv_3$ for it, as follows:
$$\vv_1,\vv_2,\vv_3 = \cases
(a{+}1,a,a),\,(a,{-}a,a{+}1),\,(a{+}1,a{-}1,{-}a{-}1)
&\text{if $k \equiv 0 \pmod 3$,} \\
(a,a,a),\,(a{+}1,{-}a,a{-}1),\,(a{-}1,a{+}1,{-}a)
&\text{if $k \equiv 1 \pmod 3$,} \\
(a,a,a{-}1),\,(a{-}1,{-}a,a),\,(a,a{-}1,{-}a)
&\text{if $k \equiv 2 \pmod 3$.} \endcases$$

A simple determinant computation shows that $\indx \Z^3:L_k|$ is
$(2a^2+a+1)(2a+1)$, $4a^3+3a$, or $(2a^2-a+1)(2a-1)$ in the respective
cases $k\equiv 0$, $k \equiv 1$, or $k \equiv 2 \pmod 3$.
Since $a$ is respectively $2k/3$, $(2k+1)/3$, or $(2k+2)/3$,
the index $\indx \Z^3:L_k|$ works out to be the desired value $n$.

For $k\equiv 0 \pmod 3$, the following vectors are in~$L_k$:
$$\align \vv_2 + \vv_3 &= (2a{+}1,{-}1,0), \\
\vv_1 + (2a{-}1)\vv_2 + 2a\vv_3 &= (4a^2{+}2a{+}1,0,{-}1). \endalign$$
Hence, we have $\ve_2 = (2a{+}1)\ve_1$ and $\ve_3 = (4a^2{+}2a{+}1)\ve_1$ in
$\Z^3/L_k$, so $\ve_1$ generates $\Z^3/L_k$.  So $\Z^3/L_k$ is isomorphic
to~$\Z_n$, via an isomorphism taking $\ve_1,\ve_2,\ve_3$ to
$1,2a{+}1,4a^2{+}2a{+}1$.  Similarly, for $k\equiv 1 \pmod 3$ we have
$$\align a\vv_2+(a{-}1)\vv_3 &= (2a^2{-}a{+}1,{-}1,0), \\
(a{+}1)\vv_2+a\vv_3 &= (2a^2{+}a{+}1,0,{-}1), \endalign$$
so $\Z^3/L_k$ is isomorphic to $\Z_n$ with generators
$1,2a^2{-}a{+}1,2a^2{+}a{+}1$; and for $k\equiv 2 \pmod 3$ we have
$$\align \vv_2 + \vv_3 &= (2a{-}1,{-}1,0), \\
\vv_1 + (2a{-}1)\vv_2 + 2a\vv_3 &= (4a^2{-}2a{+}1,0,{-}1), \endalign$$
so $\Z^3/L_k$ is isomorphic to~$\Z_n$ with generators
$1,2a{-}1,4a^2{-}2a{+}1$.

It remains to show that $\octa 3k + L_k = \Z^3$.  We will do only the
case $k \equiv 1 \pmod 3$ here; the other two cases are handled
by the same method, but with a few more subcases because of less
symmetry.

For $k=1$ one just has to show that $\Z_7$ with generators $1,2,4$ has
diameter $1$, and this is trivial to do directly; so we may assume
$k>1$ and hence $a > 1$.

Let $\vv_4 = \vv_1-\vv_2-\vv_3 = ({-}a,a{-}1,a{+}1)$.  Then the vectors
$\pm \vv_i$ for $i=1,2,3,4$ give one member of $L_k$ strictly within
each of the eight octants of $\Z^3$, and all of the coordinates of
these vectors have absolute value at most $a{+}1$.

We must show that each $\vx \in \Z^3$ is in $\octa 3k + L_k$.  This
is equivalent to showing that there is a member~$\vw$ of~$L_k$
such that $\vx - \vw \in \octa 3k$, which in turn is equivalent
to $\dist(\vx,\vw) \le k$, where $\dist$~is the $\l^1$ (Manhattan)
metric on~$\Z^3$.  Note that if $\vx,\vy,\vz$ are such that
each coordinate of~$\vy$ is between (inclusively)
the corresponding coordinates
of $\vx$ and~$\vz$, then $\dist(\vx,\vy)+\dist(\vy,\vz) =
\dist(\vx,\vz)$.  From now on, we will state this situation more
briefly as ``$\vy$~lies between $\vx$ and~$\vz$.''

Suppose we are given $\vx \in \Z^3$.  The idea is to repeatedly
reduce $\vx$ by adding members of~$L_k$ to it, until one reaches
a vector which is within $\l^1$\snug-distance $k$ of $\vzero$
or some other known member of~$L_k$.

The first thing we will do is reduce $\vx$ to a vector whose coordinates
all have absolute value at most $a{+}1$.  Suppose $\vx$ does not already
have this property.  Let $\vv$ be one of the vectors $\pm \vv_i$
($i=1,2,3,4$) such that the coordinates of $\vv$ have the same
signs as the corresponding coordinates of $\vx$; if a coordinate
of $\vx$ is $0$, then either sign is allowed for the corresponding
coordinate of $\vv$.  Now look at $\vx' = \vx - \vv$.  If a coordinate
of $\vx$ has absolute value $\le a{+}1$, then the corresponding
coordinate of $\vx'$ will also have absolute value $\le a{+}1$,
because of the sign matching and the fact that the
coordinates of $\vv$ have absolute value $\le a{+}1$.  If a coordinate
of $\vx$ has absolute value $> a{+}1$, then the corresponding
coordinate of $\vx'$ will be strictly smaller in absolute value.
Therefore, repeating this procedure will lead after finitely many steps
to a vector whose coordinates all have absolute value at most $a{+}1$.

If this new $\vx$ lies between $\vzero$ and one of the vectors
$\pm \vv_i$, then we have $\dist(\vzero,\vx)+\dist(\vx,\pm \vv_i)
= \dist(\vzero,\pm \vv_i)$.  But all of the vectors $\pm \vv_i$
satisfy $\dist(\vzero,\pm \vv_i) = 2k+1$; since
$\dist(\vzero,\vx)$ and $\dist(\vx,\pm \vv_i)$ are both integers,
one of them must be at most $k$, so we are done with this $\vx$.

We now break into cases depending on which octant the new $\vx$ lies
in.  Since $L_k$ is centrosymmetric, we only need to handle
the octants containing $\vv_1$, $\vv_2$, $\vv_3$, and $\vv_4$.
Also, $L_k$ is invariant under cyclic permutations of the three
coordinates, since these leave $\vv_1$ fixed and permute
$\vv_2,\vv_3,\vv_4$; hence, we may assume that the new $\vx$ is in
the octant of $\vv_1$ or the octant of $\vv_2$.

First suppose that $\vx$ is now in the octant of $\vv_1$ (all
three coordinates nonnegative).  If $\vx$ is between $\vzero$ and $\vv_1$,
we are done.  If two or more of the coordinates of $\vx$ are
equal to $a{+}1$, say (by cyclic symmetry of $L_k$)
$\vx = (a{+}1,a{+}1,r)$, then we have $\dist(\vx,\vv_1) \le k$
unless $k = 4$ and $r = 0$, in which case
$\dist(\vx,\vv_1+\vv_3) = 4 = k$.

If $\vx$ has exactly one coordinate equal to $a{+}1$, say
$\vx = (a{+}1,r,s)$ with $0 \le r,s \le a$,  then we can
subtract $\vv_1$ from $\vx$ to get $\vx' = (1,r{-}a,s{-}a)$, which
is in the octant containing ${-}\vv_4$.  If $\vx'$ lies between
$\vzero$ and $-\vv_4$, we are done.  If not, then $r = 0$.
Now let $\vx'' = \vx' + \vv_4 = ({-}a{+}1,{-}1,s{+}1)$, which
lies between $\vzero$ and $-\vv_3$ unless $s = a$, in which
case $\vx = (a{+}1,0,a)$ and $\dist(\vx,\vv_2) = a{+}1 \le k$.

The procedure when $\vx$ is in the octant of $\vv_2$ is similar.
Either $\vx=(r,s,t)$ lies between $\vzero$ and $\vv_2$, or $s={-}a{-}1$,
or $t \ge a$.  In the latter cases, let $\vx' = \vx - \vv_2$.
When $s={-}a{-}1$, we have $\vx' = (r{-}a{-}1,{-}1,t{-}a{+}1)$;
either this lies between $\vzero$ and one of the vectors
$\pm \vv_i$, or $\vx+\vv_3$ does.  When $s \ge {-}a$ but
$t \ge a$, try $\vx' - \vv_4$; either it lies between $\vzero$
and some $\pm \vv_i$, or it is $({-}1,1,{-}a)$, $(a,1,{-}a)$,
or $(a,1,{-}a{+}1)$.  These last three lie within $\dist$\snug-distance~$k$
of $\vv_3-\vv_1$, $\vv_3$, and either
$\vv_3+\vv_2-\vv_1$ or $\vv_3+\vv_2$, respectively. \QED

The authors conjecture that the graphs given by this theorem are
actually the largest undirected Cayley graphs of Abelian groups
on three generators for each diameter~$k$.

This conjecture would imply that the lattice covering of~$\R^3$
by~$\Rocta 3{3/2}$ using the lattice~$L_{\text{bcc}}$ is optimal;
that is, $8/9$~is the best possible efficiency for a lattice covering
by regular octahedra.  The latter statement seems quite plausible,
but remains unproved at this point.  We can prove a partial result,
though, that a ``small'' adjustment to~$L_{\text{bcc}}$ cannot
improve the covering:

\proclaim{Theorem 16} Among those lattices~$L$ for which
$\Rocta 3{3/2} + L = \R^3$, the lattice~$L_{\text{bcc}}$ is
locally optimal; that is, for any other lattice~$L$
sufficiently near~$L_{\text{bcc}}$ such that $\Rocta 3{3/2} + L = \R^3$,
the efficiency of the covering using~$L$ is less than~$8/9$.
\endproclaim

\demo{Proof} Let us use the vectors $\vv_1 = (1,1,1)$, $\vv_2 = (1,-1,1)$,
and $\vv_3 = (1,1,-1)$ as generating vectors for~$L_{\text{bcc}}$; then
a nearby lattice~$L$ will be generated by nearby vectors $\vv'_1 =
(a_1,b_1,c_1)$, $\vv'_2 = (a_2,b_2,c_2)$, and $\vv'_3 = (a_3,b_3,c_3)$.
We can concatenate the three vectors $\vv'_1,\vv'_2,\vv'_3$ to get
a single vector~$\vv'$ in~$\R^9$; similarly, let
$\vv$ be the concatenation $\vv_1,\vv_2,\vv_3$.
Let $F(\vv')$ be the determinant
of the matrix with rows $\vv'_1,\vv'_2,\vv'_3$.  Note that
$F(\vv) = 4$; we must see that this point is a
strict local maximum of~$F(\vv')$ for those points~$\vv'$ satisfying
the constraint that $\Rocta 3{3/2} + L = \R^3$.
We compute that the gradient of~$F$ at the point~$\vv$ is
$\vg = (0,2,2,2,-2,0,2,0,-2)$.

Using the lattice $L_{\text{bcc}}$, the point $(1/2,1/2,1/2)$, in the
center of a face of~$\Rocta 3{3/2}$, is covered
by only two copies of~$\Rocta 3{3/2}$, namely $\Rocta 3{3/2}$~itself
and $\Rocta 3{3/2} + \vv_1$, and it is on the boundary (a face) of
each of these copies.  If the lattice is altered slightly so that
these two copies no longer touch, then the points in between will not
be covered by any copy.  In particular, if $L$~is near~$L_{\text{bcc}}$
but $a_1 + b_1 + c_1 > 3$, then the point $(1/2,1/2,1/2+\eps)$ for small
positive~$\eps$ will not be in $\Rocta 3{3/2} + L$.  So the
constraint $\Rocta 3{3/2} + L = \R^3$ gives us the linear inequality
$a_1 + b_1 + c_1 \le 3$.  We will rewrite this as
$$\align
\vu_1 \vdot \vv' \le 3,\qquad \text{where }&\vu_1 = (1,1,1,0,0,0,0,0,0).\\
\intertext{The same argument for points on the other faces of the octahedron
gives inequalities}
\vu_2 \vdot \vv' \le 3,\qquad \text{where }&\vu_2 = (0,0,0,1,-1,1,0,0,0), \\
\vu_3 \vdot \vv' \le 3,\qquad \text{where }&\vu_3 = (0,0,0,0,0,0,1,1,-1), \\
\vu_4 \vdot \vv' \le 3,\qquad \text{where }&\vu_4 = (-1,1,1,1,-1,-1,1,-1,-1).
\endalign$$

Next, consider the point $(1,0,1/2)$.  This is in $\Rocta 3{3/2} + \vy$
for four members~$\vy$ of~$L_{\text{bcc}}$, namely $\vzero$, $\vv_1$,
$\vv_2$, and $\vv_2+\vv_3$, and it is an edge point of each of these
four copies.  If the lattice is altered slightly, then a gap can open up
near this point even if there are no gaps between octahedra adjacent at
a face as above.

Specifically, if $\vv'$~is near~$\vv$, $\eps$~is a very small positive
number, and we define the point~$\vx$ by the linear equations
$$\align
\vx \vdot (1,-1,-1) &= \vv'_1 \cdot (1,-1,-1) + 3/2 + \eps, \\
\vx \vdot (-1,-1,1) &= (\vv'_2+\vv'_3) \cdot (-1,-1,1) + 3/2 + \eps,
 \quad \text{and}\\
\vx \vdot (-1,1,-1) &= \vv'_2 \cdot (-1,1,-1) + 3/2 + \eps,
\endalign$$
then $\vx$~will be a point near $(1,0,1/2)$ which is not in
$\Rocta 3{3/2} + \vy$ for $\vy \in \{\vv'_1,\vv'_2,\vv'_2+\vv'_3\}$.
Adding up the three given equations yields
$$
\vx \vdot (-1,-1,-1) = \vv'_1 \cdot (1,-1,-1) + \vv'_2 \cdot (-2,0,0) +
\vv'_3 \cdot (-1,-1,1) + 9/2 + 3\eps.
$$
If the right hand side of this equation is less than $-3/2$, then
$\vx$~will not be in~$\Rocta 3{3/2}$ either, and hence will not be in
$\Rocta 3{3/2} + L$.  Since $\eps$~can be arbitrarily small, in order
to have $\Rocta 3{3/2} + L = \R^3$, it is necessary to have
$$\vv'_1 \cdot (1,-1,-1) + \vv'_2 \cdot (-2,0,0) +
\vv'_3 \cdot (-1,-1,1) \ge -6.$$
This can be rewritten as
$$\align
\vu_5 \vdot \vv' \le 6,\qquad \text{where }&\vu_5 = (-1,1,1,2,0,0,1,1,-1). \\
\intertext{The same argument can be performed using the octahedra around
$(1,0,1/2)$ in the opposite order, and there are $23$~other points
on the edges of~$\Rocta 3{3/2}$ where the same configuration occurs.
But one only gets six distinct inequalities from this; the other
five are:}
\allowdisplaybreak
\vu_6 \vdot \vv' \le 6,\qquad \text{where }&\vu_6 = (0,2,0,1,-1,1,1,-1,-1), \\
\vu_7 \vdot \vv' \le 6,\qquad \text{where }&\vu_7 = (-1,1,1,1,-1,1,2,0,0), \\
\vu_8 \vdot \vv' \le 6,\qquad \text{where }&\vu_8 = (1,1,1,0,-2,0,1,-1,-1), \\
\vu_9 \vdot \vv' \le 6,\qquad \text{where }&\vu_9 = (0,0,2,1,-1,-1,1,1,-1), \\
\vu_{10} \vdot \vv' \le 6,\qquad \text{where }&\vu_{10} =
  (1,1,1,1,-1,-1,0,0,-2).
\endalign$$

Note that all ten of these inequalities are satisfied with equality
when $\vv'=\vv$.  Hence, they can be rewritten as
$\vu_i \vdot (\vv'-\vv) \le 0$ for $i =1,2,\dots,10$.

One can easily check that the vectors $\vu_1,\dots,\vu_7$ are linearly
independent; their common null space (i.e., the set of~$\vw$ such
that $\vu_i \vdot \vw = 0$ for all $i \le 7$) is generated by
the independent vectors $\vw_1 = (1,0,-1,1,0,-1,1,0,1)$ and
$\vw_2 = (-1,1,0,-1,-1,0,-1,1,0)$.  Also, we have
$$\vg = \vu_1 + \vu_2 + \vu_3 + \vu_4 = \vu_5 + \vu_8 = \vu_6 + \vu_9 =
\vu_7 + \vu_{10}.$$

Let $C$~be the closed cone consisting of all vectors~$\vt$ in the
subspace spanned by $\vu_1,\dots,\vu_7$
such that $\vu_i \vdot \vt \le 0$ for all $i \le 10$.  Then
the above equations imply that $\vg \cdot \vt \le 0$ for all
$\vt$ in~$C$, and equality can hold only when $\vt = \vzero$.
In particular, we have $\vg \cdot \vt_0 < 0$ for any unit vector~$\vt_0$
in~$C$.  The set of such~$\vt_0$ is closed and bounded, hence compact,
so there is a positive number~$\eps$ such that $\vg \cdot \vt_0 < -\eps$
for all such~$\vt_0$.
It follows that there is a neighborhood~$U$ of~$\vg$ such that,
for any $\vg'$ in~$U$ and any unit vector $\vt_0$ in~$C$,
$\vg' \vdot \vt_0 < 0$.  Since $C$~is a cone, we have
$\vg' \vdot \vt < 0$ for all $\vg' \in U$ and all nonzero $\vt \in C$.

We can compute that, for any real numbers $r$ and~$s$,
the deteriminant for the lattice given by $\vv + r \vw_1 + s \vw_2$ is
$$F(\vv + r \vw_1 + s \vw_2) = 4(1-r)(1+s)(1+r-s).$$
If $|r|+|s| < 1$, then the numbers $1-r$, $1+s$, and $1+r-s$ are
positive numbers with arithmetic mean~$1$, so their geometric mean is at
most~$1$; this means that $F(\vv + r \vw_1 + s \vw_2) \le 4$.
Equality holds only when the above three numbers are equal, which is
when $r=s=0$.

Let $U'$ be a convex neighborhood of~$\vv$ so small that $(\grad
F)(\vv') \in U$ for all $\vv' \in U'$.  Now, any vector $\vv'$
sufficiently close to~$\vv$ can be expressed as $\vv + \vt_1 + \vt_2$
where $\vt_1$ is a (small) linear combination of $\vw_1$ and $\vw_2$,
$\vt_2$ is a linear combination of $\vu_1,\dots,\vu_7$, and both
$\vv + \vt_1$ and $\vv + \vt_1 + \vt_2$ are in $U'$.
If $\vv'$ satisfies the condition $\Rocta 3{3/2} + L = \R^3$ and is near
to $\vv$, then we must have $\vu_i \vdot (\vv' - \vv) \le 0$ for
all $i \le 10$, so $\vu_i \vdot \vt_2 \le 0$ for all $i \le 10$
(since $\vu_i \vdot \vt_1 = 0$), so $\vt_2 \in C$.
We have $F(\vv + \vt_1) \le 4$, with equality holding only when $\vt_1
= \vzero$.  If $\vt_2$ is nonzero, then for any~$\vt$ on the segment
from $\vv+\vt_1$ to $\vv+\vt_1+\vt_2$ we have $\vt \in U'$,
so $(\grad F)(\vt) \in U$, so $(\grad F)(\vt) \vdot \vt_2 < 0$;
it follows that $F(\vv + \vt_1 + \vt_2) < F(\vv + \vt_1)$.
Therefore, $F(\vv') \le F(\vv)$, with equality holding only when
$\vv' = \vv$.  So $\vv$~gives a local maximum of~$F$, as desired.
\QED

It is still possible (though very unlikely) that a lattice quite
different from~$L_{\text{bcc}}$ gives a more efficient covering.
Theoretically, the search for an optimal lattice can be set up
as a large optimization problem and solved once and for all,
but this appears to be a formidable task.

One could begin this task by considering an arbitrary lattice~$L$
such that $\Rocta 3{3/2} + L = \R^3$ and this covering is reasonably
efficient (at least as efficient as the covering from~$L_{\text{bcc}}$).
Such a lattice is generated by vectors $\vv_1,\vv_2,\vv_3$, and we
can carefully choose these generators so as to limit their lengths.
In particular, we can choose~$\vv_1$ to be a nonzero member of~$L$
with minimal length.  We can then choose $\vv_2$ in~$L$ whose distance
from the subspace of~$\R^3$ spanned by~$\vv_1$ is as small as possible
(but nonzero), and adjust~$\vv_2$ by subtracting off an integer multiple
of~$\vv_1$ so as to ensure that the closest integer multiple of~$\vv_1$
to~$\vv_2$ is~$\vzero$.  One can similarly choose~$\vv_3$ to be as
close as possible to (but not in) the subspace spanned by $\vv_1$
and~$\vv_2$.  These three chosen vectors will be a set of generating
vectors for~$L$.  In order to have $\Rocta 3{3/2} + L = \R^3$, it is
necessary that the length of~$\vv_1$ be no more than the diameter
of~$\Rocta 3{3/2}$; there are similar but slightly larger bounds on the
lengths of $\vv_2$ and~$\vv_3$.  This limits our search for
$\vv_1,\vv_2,\vv_3$ to a compact subset of nine-dimensional space.
We must find the point in this subset which maximizes
$\det(\vv_1,\vv_2,\vv_3)$ subject to the constraint that
$\Rocta 3{3/2} + L = \R^3$.

This constraint looks infinitary, but it can actually be reduced to
finitely many sets of linear inequalities.  To see this, note that,
using the above upper bounds on the lengths of the vectors~$\vv_i$
along with the assumed lower bound on the lattice determinant
(the covering must be at least as efficient as that from~$L_{\text{bcc}}$),
we can get lower bounds on the lengths of the~$\vv_i$, the angles
between them, and associated quantities such as the distance from~$\vv_3$
to the plane spanned by $\vv_1$ and~$\vv_2$.  These will allow us
to get upper bounds on the absolute values of integers $a_1,a_2,a_3$
such that $\Rocta 3{3/2} + a_1\vv_1 + a_2 \vv_2 + a_3 \vv_3$
overlaps or almost touches~$\Rocta 3{3/2}$.  (In other words, we get
an upper bound on the number~$M$ from the proof of Proposition~7.)
So we only have to consider finitely many of the lattice translates
of~$\Rocta 3{3/2}$ when trying to cover the space near~$\Rocta 3{3/2}$
(which is all that is needed, by Proposition~6).

There are only finitely many configurations (specifications of arrangements
and overlaps) for these finitely many translates of~$\Rocta 3{3/2}$.
For each such configuration, the assertion that there are no
`gaps' in the coverage of the space near~$\Rocta 3{3/2}$ becomes
a list of linear inequalities like the inequalities $\vu_i \cdot \vv'
\le b$ from the proof of Theorem~16.  So we need to optimize a
cubic function (the lattice determinant) subject to a list of
linear inequalities in order to find the optimal version of each
configuration, and then compare the resulting values to find the
best configuration.

Unfortunately, there is a very large number of possible configurations
(for an example of the possibilities for complicated configurations,
see Figure~4 later in this paper), so this finite computation
appears to be beyond reach at present.  Of course, a different approach
to the problem might lead to a more feasible computation.

\smallpagebreak

One might hope to be able to use the arguments of Proposition~7 and
Theorem~8 in reverse, to get an upper bound on the efficiency of
lattice coverings of~$\R^3$ by the octahedron~$\Rocta 31$
by showing that any extremely
efficient real lattice covering would lead to integer lattice
coverings more efficient than what the computation actually found.
To do this, one would fix a value for the distance~$\rho$ from
Proposition~7, and then use the method described above to get an
upper bound on the number~$M$ from that Proposition.  If there is
actually a lattice covering of~$\R^3$ by~$\Rocta 31$ using the
lattice~$L$ generated by $\vv_1,\vv_2,\vv_3$ having
a specified large determinant (equivalently, a specified
large efficiency), then we can
round the coordinates of these vectors to the nearest multiples
of~$1/k$ to get vectors $\vv'_1,\vv'_2,\vv'_3$ generating
a lattice~$L'$.  By Proposition~7, if $1/(2k)$ is less than
$\eta = \rho/M$, then we will have $\Rocta 3{1 + 3\rho} + L' = \R^3$,
so $\Rocta 3{k + 3k\rho} + kL' = \R^3$.  But $kL'$~is an
integer lattice; if $n$~is its determinant, then this lattice
covering will yield an Abelian Cayley graph on three generators
with size~$n$ and diameter at most $k + 3k\rho$.
The fact that $L'$~is close to~$L$ means that we can get a lower
bound on~$n$ from the determinant of~$L$.  If the actual computational
search showed that there is no Abelian Cayley graph of such a size for this
diameter, then our original assumption that there was a lattice~$L$
giving a covering of that efficiency must have been false.

Unfortunately, the constants involved are such that even the large
computation done so far does not suffice to get a bound less than~$1$
for the efficiency of~$L$ (even if we are optimistic enough to assume
that $M$~is as small as 3 or~4).  It probably requires searches for
values of~$k$ larger than~$500$ in order to get actual
results from this method; such searches are completely out of range
at the moment.

\subhead Three generators, directed \endsubhead
For the directed case of three generators, we want to study lattice
coverings of~$\R^3$ by the trirectangular tetrahedron~$\Rtetra 3k$.
(Since lattice covering efficiency is affine invariant, it makes
no difference which particular tetrahedron we consider.)  One hopes
that one can discretize these coverings to get good lattice coverings
of~$\Z^3$ by~$\tetra 3k$, and hence good directed Cayley graphs.

The best three-dimensional directed toroidal mesh with diameter $k$ is
$\Z_{b_0+1} \times \Z_{b_1+1} \times \Z_{b_2+1}$,
where $b_i = \lfloor (k+i)/3 \rfloor$; this has
about $(1/27)k^3$ vertices for large $k$.  This corresponds to the
covering of $\R^3$ by $\Rtetra 31$ using the cubic lattice
$(1/3)\Z^3$; this covering has efficiency~$2/9$.

It is more difficult to find a candidate for a good covering lattice
(or, equivalently, a large subset which gives a lattice tiling) for the
tetrahedron than it was for the octahedron.  One possible method is to
try to find the three-dimensional analogue of the L-tromino used for
the triangle; this leads one to consider the tetracube shown on the
left of Figure~3.  In order for the shape to fit into $\Rtetra 31$,
the edge-length of the subcubes should be $1/4$.  It is easy to see that
this shape does indeed tile space, using the lattice generated by
$(1/2,0,0)$, $(0,1/2,0)$, and $(1/4,1/4,-1/4)$ (this is just
$(1/4)L_{\text{bcc}}$); since the shape has volume $1/16$ while
$\Rtetra 31$ has volume $1/6$, we get a lattice covering of
$\R^3$ by $\Rtetra 31$ with efficiency $3/8$.

\midinsert
\hbox to\hsize{
\hfill
\plotfigbegin
   \plotfigscales 28.5 16.5
   \plotvskip 7
   \plottext{}
   \plot 4 7
   \plotmove 1 0
   \plotline 0 1
   \plotline 0 3
   {\plotline 1 2 }
   \plotline 1 4
   \plotline 3 2
   \plotline 4 3
   {\plotline 3 4 }
   \plotline 4 1
   \plotline 3 0
   {\plotline 3 2 }
   \plotline 2 1
   \plotline 2 5
   \plotline 1 6
   {\plotline 1 4 }
   \plotline 2 7
   \plotline 3 6
   {\plotline 2 5 }
   \plotline 3 4
   \plotline 1 2
   \plotline 1 0
   \plotline 2 1
   \plotfigscales 19 11
   \plotvskip 0.5
\plotfigend
\hfill
\plotfigbegin
   \plotfigscales 19 11
   \plotvskip 11
   \plottext{}
   \plot 6 11
   \plotmove 1 3
   \plotline 3 5
   {\plotline 5 3 }
   \plotline 3 9
   \plotline 4 8
   {\plotline 5 9 }
   \plotline 4 6
   {\plotline 5 7 }
   \plotline 5 5
   {\plotline 6 6 }
   \plotline 5 3
   \plotline 4 2
   {\plotline 4 0 }
   \plotline 3 3
   {\plotline 3 1 }
   \plotline 2 2
   {\plotline 2 0 }
   \plotline 1 3
   \plotline 1 5
   {\plotline 2 6 }
   \plotline 0 6
   {\plotline 1 7 }
   \plotline 0 2
   \plotline 2 0
   \plotline 3 1
   \plotline 4 0
   \plotline 6 2
   \plotline 6 6
   \plotline 5 7
   \plotline 5 9
   \plotline 3 11
   \plotline 1 9
   {\plotline 2 8 }
   \plotline 1 7
   \plotline 2 6
   \plotline 2 8
   \plotline 3 9
\plotfigend
\hfill}
\botcaption{Figure 3}
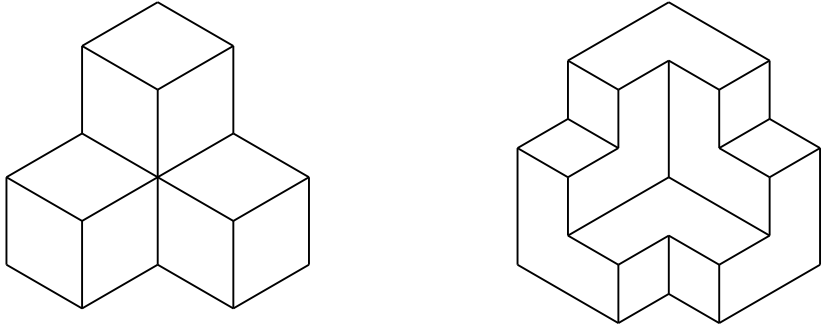
Two subsets of the tetrahedron $\Rtetra 31$ which tile space.
\endcaption
\endinsert

The discrete form of this shape, scaled by a factor $s_i$ in the
$i$\snug'th dimension, is a subset of $\Z^3$ of size $4s_1s_2s_3$
which gives a lattice tiling of $\Z^3$; this subset is included
in $\tetra 3k$, where $k = s_1+s_2+s_3+\max(s_1,s_2,s_3)-3$.
Optimizing this for a given $k \ge 1$ gives a subset of~$\tetra 3k$
which tiles and has size $4a_3a_4a_5$, where $a_i = \lfloor (k+i)/4 \rfloor$.

One can obtain another lattice covering of $\R^3$ by tetrahedra as follows.
If one cuts off the four corners of a regular tetrahedron at planes
passing through the midpoints of the edges (so one removes four half-size
regular tetrahedra), then what is left is a regular octahedron with
volume $1/2$ that of the tetrahedron.  We have a lattice giving
a covering of~$\R^3$ by this octahedron with efficiency $8/9$; the
same lattice therefore gives a covering of~$\R^3$ by the original
tetrahedron with efficiency~$4/9$.

If one uses an affine transformation to change the regular tetrahedron
to the tetrahedron $\Rtetra 31$, then the corresponding lattice will
be generated by $(1/6,1/6,1/6)$, $(1/6,{-}1/2,1/6)$, and
$(1/6,1/6,{-}1/2)$.  One fundamental region for this lattice
is an affinely distorted truncated octahedron.  Another one can be obtained
by the method of Lemma~3, using an ordering~$\prec$ which
orders vectors primarily by the sum of their coordinates; the
resulting region is shown at the right of Figure~3.  This shape
consists of $16$~cubes of edge-length~$1/6$, for a total volume
of~$2/27$, which, as expected, is $4/9$~of $\volume(\Rtetra 31) = 1/6$.

Discretizing this new shape with scale factors $s_1 \le s_2 \le s_3$
gives a subset of $\Z^3$ of size $16s_1s_2s_3$ which gives a lattice
tiling of $\Z^3$; this subset is included in $\tetra 3k$, where $k =
s_1+2s_2+3s_3-3$.  Another simple optimization shows that, for any
given $k \ge 3$, we get a subset of $\tetra 3k$ which tiles and has
size $16\hat a_3\hat a_4\hat a_6$, where $\hat a_i = \lfloor (k+i)/6
\rfloor$.  For large~$k$ (in fact, for all $k \ge 30$), this
new lattice gives a better covering of~$\Z^3$ by~$\tetra 3k$ than
the preceding one did, but for smaller $k$ the preceding one
sometimes does better.

Aguil\'o, Fiol, and Garcia~\cite{\AguFioGar} also work with this shape,
but discretize it in a rotationally symmetric way rather than
in each dimension separately; the Cayley graphs they obtain are
slightly larger than the graphs of size $16\hat a_3\hat a_4\hat a_6$
given above, but still of the form $(2/27) k^3 + O(k^2)$.

In order to see whether these lattice coverings give close-to-optimal
Cayley graphs, the authors performed a computer search for the
best (smallest-diameter) directed Abelian Cayley graphs
on three generators.  This extends similar computations performed
by Aguil\'o, Fiol, and Garcia~\cite{\AguFioGar} and by
Fiduccia, Forcade, and Zito~\cite{\FidForZit}.
The latter paper also contains a useful upper bound:
an Abelian Cayley digraph on three generators with diameter~$k$
has size at most~$3(k+3)^3/25$.  This improves the obvious
upper bound~$|\octa 3k|$ when $k > 7$.

Comparing the above figures with the output from
the authors' computations gives a slight surprise: the best
cyclic groups do substantially better than the groups from the
above coverings.  The data are shown in Table~2; here `FFZ' is
the Fiduccia-Forcade-Zito upper bound,
`Impr.'
refers to the larger of the sizes obtained from the two improved
constructions
above, `AFG' is the size attained by the
Aguil\'o-Fiol-Garcia
construction, and the remaining columns are analogous to
those of Table~1.  The computations were run on Abelian groups
of sizes up to and
including~$4871$; this means that the entries marked with an asterisk
in the~$n'_c$ column (for which the FFZ bound is
greater than~$4871$) have not been completely proven optimal, but
it is extremely likely that they are.

\topinsert
\centerline{
\vbox{\offinterlineskip
\halign{
&\vrule#&\hfil$\,\,\vphantom{|\Rtetra 3k|}{#}\,\,$\hfil\cr
\noalign{\hrule}
height2pt&\omit&&\omit&&\omit&&\omit&&\omit&&\omit&&\omit&&\omit&&\omit
&&\omit&\cr
&k&&|\tetra 3k|&&\omit\,\,FFZ\,\,&&\omit\hfil\,\,Toroidal\,\,\hfil&&
    \omit\hfil\,\,Impr.\,\,\hfil&&\omit\hfil\,\,AFG\,\,\hfil&&n'_c
    &&\omit\hfil\,\,Generators\,\,\hfil&&
      n'_c/|\tetra 3k|&&n'_c/\volume(\Rtetra 3{k+3})&\cr
height2pt&\omit&&\omit&&\omit&&\omit&&\omit&&\omit&&\omit&&\omit&&\omit
&&\omit&\cr
\noalign{\hrule}
height2pt&\omit&&\omit&&\omit&&\omit&&\omit&&\omit&&\omit&&\omit&&\omit
&&\omit&\cr
&0&&1&&&&1&&&&1&&1&&&&1&&.222222&\cr
&1&&4&&&&2&&4&&4&&4&&1,2,3&&1&&.375000&\cr
&2&&10&&&&4&&4&&7&&9&&1,3,4&&.900000&&.432000&\cr
&3&&20&&&&8&&16&&16&&16&&1,4,5&&.800000&&.444444&\cr
&4&&35&&&&12&&16&&19&&27&&1,4,17&&.771428&&.472303&\cr
&5&&56&&&&18&&32&&31&&40&&1,6,15&&.714286&&.468750&\cr
&6&&84&&&&27&&32&&50&&57&&1,13,33&&.678571&&.469136&\cr
&7&&120&&120&&36&&48&&56&&84&&2,9,35&&.700000&&.504000&\cr
&8&&165&&159&&48&&72&&86&&111&&1,31,69&&.672727&&.500376&\cr
&9&&220&&207&&64&&128&&128&&138&&1,11,78&&.627273&&.479167&\cr
&10&&286&&263&&80&&128&&134&&176&&1,17,56&&.615385&&.480655&\cr
&11&&364&&329&&100&&144&&182&&217&&1,13,119&&.596154&&.474490&\cr
&12&&455&&405&&125&&192&&243&&273&&1,14,153&&.600000&&.485333&\cr
&13&&560&&491&&150&&256&&252&&340&&1,90,191&&.607143&&.498047&\cr
&14&&680&&589&&180&&288&&333&&395&&1,35,271&&.580882&&.482394&\cr
&15&&816&&699&&216&&432&&432&&462&&1,29,97&&.566176&&.475309&\cr
&16&&969&&823&&252&&432&&441&&560&&1,215,326&&.577915&&.489867&\cr
&17&&1140&&960&&294&&500&&549&&648&&1,76,237&&.568421&&.486000&\cr
&18&&1330&&1111&&343&&576&&676&&748&&1,41,147&&.562406&&.484613&\cr
&19&&1540&&1277&&392&&600&&688&&861&&1,27,463&&.559091&&.485162&\cr
&20&&1771&&1460&&448&&768&&844&&979&&1,22,351&&.552795&&.482781&\cr
&21&&2024&&1658&&512&&1024&&1024&&1140&&1,45,196&&.563241&&.494792&\cr
&22&&2300&&1875&&576&&1024&&1036&&1305&&1,246,1030&&.567391&&.501120&\cr
&23&&2600&&2109&&648&&1024&&1228&&1440&&1,126,415&&.553846&&.491579&\cr
&24&&2925&&2361&&729&&1280&&1445&&1616&&1,56,257&&.552479&&.492608&\cr
&25&&3276&&2634&&810&&1372&&1460&&1788&&1,154,1452&&.545788&&.488703&\cr
&26&&3654&&2926&&900&&1600&&1715&&1963&&1,90,780&&.537219&&.482923&\cr
&27&&4060&&3240&&1000&&2000&&2000&&2224&&1,425,704&&.547783&&.494222&\cr
&28&&4495&&3574&&1100&&2000&&2015&&2442&&1,964,1372&&.543270&&.491826&\cr
&29&&4960&&3932&&1210&&2048&&2315&&2693&&1,39,942&&.542944&&.493103&\cr
&30&&5456&&4312&&1331&&2400&&2646&&2920&&1,540,831&&.535191&&.487520&\cr
&31&&5984&&4716&&1452&&2400&&2664&&3220&&7,30,2277&&.538102&&
.491553&\cr
&32&&6545&&5145&&1584&&2880&&3042&&\hphantom{{}^*}3591^*&&1,1519,2031&&.548663&&
.502531&\cr
&33&&7140&&5598&&1728&&3456&&3456&&\hphantom{{}^*}3850^*&&2,475,1177&&.539216&&
.495113&\cr
&34&&7770&&6078&&1872&&3456&&3474&&\hphantom{{}^*}4191^*&&1,748,2652&&.539382&&
.496437&\cr
&35&&8436&&6584&&2028&&3456&&3906&&\hphantom{{}^*}4468^*&&1,353,2789&&.529635&&
.488555&\cr
height2pt&\omit&&\omit&&\omit&&\omit&&\omit&&\omit&&\omit&&\omit&&\omit
&&\omit&\cr
\noalign{\hrule}
}}
}
\botcaption{Table 2}
Best directed Cayley graphs of cyclic groups, three generators.
\endcaption
\endinsert

Note that in three cases, $k=7,31,33$, the best cyclic Cayley graph
was not achieved using~$1$ as one of the generators.  If one is required
to use~$1$ as a generator (which may be useful when actually building
the corresponding loop network), then the best one can do is
size~$78$ for $k=7$ (with generators $1,6,49$), size~$3178$ for
$k=31$ (with generators $1,386,1295$), and size~$3794$ for
$k=33$ (with generators $1,469,2094$).

There is one other difference between this case and the undirected case:
here there are values of~$k$ for which one can do better with general
Abelian groups than with cyclic groups.  The improved values obtained
from non-cyclic groups are shown in Table~3.  A number of these optimal
graphs are actually obtained by applying Proposition~5(b) to
smaller Cayley graphs; for instance, the Abelian graph for $k=17$
is obtained this way from the cyclic graph for $k=7$, which is the
reason that these two graphs give exactly the same real-covering
efficiency (.504).

\topinsert
\centerline{
\vbox{\offinterlineskip
\halign{
&\vrule#&\hfil$\,\,\vphantom{|\Rtetra 3k|}{#}\,\,$\hfil\cr
\noalign{\hrule}
height2pt&\omit&&\omit&&\omit&&\omit&&\omit&&\omit&\cr
&k&&n'_a&&\omit\hfil\,\,Group\,\,\hfil
    &&\omit\hfil\,\,Generators\,\,\hfil&&
      n'_a/|\tetra 3k|&&n'_a/\volume(\Rtetra 3{k+3})&\cr
height2pt&\omit&&\omit&&\omit&&\omit&&\omit&&\omit&\cr
\noalign{\hrule}
height2pt&\omit&&\omit&&\omit&&\omit&&\omit&&\omit&\cr
&12&&279&&\Z_{93}\times\Z_3&&(1,0),(9,1),(10,2)&&.613187&&.496000&\cr
&17&&672&&\Z_{168}\times\Z_2\times\Z_2&&(2,1,0),(9,0,0),(35,0,1)&&
.589474&&.504000&\cr
&18&&752&&\Z_{188}\times\Z_4&&(1,0),(13,2),(14,1)&&.565414&&.487204&\cr
&19&&888&&\Z_{222}\times\Z_2\times\Z_2&&(1,0,0),(142,1,0),(180,0,1)&&
.576623&&.500376&\cr
&26&&1980&&\Z_{330}\times\Z_6&&(1,0),(123,2),(234,3)&&.541872&&.487105&\cr
&27&&2268&&\Z_{252}\times\Z_3\times\Z_3&&(2,0,0),(9,1,0),(35,0,1)&&
.558621&&.504000&\cr
&28&&2448&&\Z_{816}\times\Z_3
&&(1,0),(427,0),(564,1)&&.544605&&.493035&\cr
&29&&2720&&\Z_{680}\times\Z_2\times\Z_2
&&(1,0,0),(191,1,0),(90,0,1)&&.548387&&.498047&\cr
&30&&2997&&\Z_{333}\times\Z_3\times\Z_3
&&(1,0,0),(31,1,0),(180,0,1)&&.549304&&.500376&\cr
&35&&\hphantom{{}^*}4500^*&&\Z_{300}\times\Z_{15}
&&(1,0),(3,1),(214,7)&&.533428&&.492054&\cr
height2pt&\omit&&\omit&&\omit&&\omit&&\omit&&\omit&\cr
\noalign{\hrule}
}}
}
\botcaption{Table 3}
Best directed Cayley graphs of Abelian groups, three generators.
\endcaption
\endinsert

The values in the $n'_c$~column of Table~2 are so much larger than
those in the preceding two columns that it is clear that the real lattices
used for the preceding columns were not optimal.  This is made explicit
in the last column of the table, which gives the efficiency of the
real lattice covering obtained from the computed integer lattice covering via
Proposition~4(b).  For $k=1$ and $k=3$ these coverings are just
(scaled versions of) the two coverings we explicitly constructed above;
but later coverings obviously do substantially better.

The best real covering obtained from these computations is that for $k=7$,
with efficiency~$.504$.  As in the undirected case, we can reconstruct
generators for the lattice from the given generating set ${2,9,35}$
for~$\Z_{84}$; after simplification, the resulting generating
vectors are $({-}2,2,2)$, $(3,{-}3,3)$, and $(4,3,{-}1)$.  We now
have a computer-assisted proof that the lattice generated by
these vectors gives a lattice covering of $\R^3$ by $\Rtetra 3{10}$;
but one can obtain useful extra information (as well as, perhaps,
more satisfaction) by proving this directly.

\proclaim{Proposition 17} Let $L'_7\badadjustI$ be the lattice in $\R^3\badadjustII$
generated by the vectors $({-}2,2,2)$, $(3,{-}3,3)$, and $(4,3,{-}1)$;
then $\Rtetra 3{10} + L'_7 = \R^3$. \endproclaim

\demo{Proof} First, note that the following vectors are in $L'_7$:
$$\alignat 4
\vv_1 &= (-2,2,2) &&
      &\qquad\qquad \vv_8 &= (1,-1,5) &&= \vv_1 + \vv_2 \\
\vv_2 &= (3,-3,3) &&
      &\qquad\qquad \vv_9 &= (-1,1,7) &&= 2\vv_1 + \vv_2 \\
\vv_3 &= (4,3,-1) &&
      &\qquad\qquad \vv_{10} &= (3,4,-6) &&= \vv_3 - \vv_1 - \vv_2 \\
\vv_4 &= (6,1,-3) &&= \vv_3 - \vv_1
      &\qquad\qquad \vv_{11} &= (-7,7,1) &&= 2\vv_1 - \vv_2 \\
\allowdisplaybreak
\vv_5 &= (5,-5,1) &&= \vv_2 - \vv_1
      &\qquad\qquad \vv_{12} &= (-5,-2,8) &&= 2\vv_1 + \vv_2 - \vv_3 \\
\vv_6 &= (1,6,-4) &&= \vv_3 - \vv_2
      &\qquad\qquad \vv_{13} &= (-1,8,-2) &&= \vv_1 - \vv_2 + \vv_3 \\
\vv_7 &= (2,5,1) &&= \vv_1 + \vv_3
      &\qquad\qquad \vv_{14} &= (8,-1,-5) &&= \vv_3 - 2\vv_1
\endalignat$$
For each vector $\vv_i = (r,s,t)$, we have $1 \le r+s+t \le 8$; hence,
the translated tetrahedron $\Rtetra 3{10} + \vv_i$ intersects the
plane $x+y+z=10$.  In fact, the intersection is a triangle whose
vertices have coordinates $(10{-}s{-}t,s,t)$, $(r,10{-}r{-}t,t)$,
and $(r,s,10{-}r{-}s)$.

Figure~4 shows the upper face of $\Rtetra 3{10}$.  For each $i \le
14$, it indicates which part of this face is covered by the translate
$\Rtetra 3{10}+\vv_i$.  (The face is divided up into unit triangles,
each of which is labeled by the value(s) of $i$ for which $\Rtetra
3{10}+\vv_i$ covers that triangle.)  Clearly each unit triangle is
labeled, so the translates $\Rtetra 3{10}+\vv_i$ cover the entire
upper face of $\Rtetra 3{10}$.

\midinsert
\centerline{
\plotfigbegin
   \plotfigscales 16.5 28.5
   \plotdotspc=3\plotdotspc
   \plotlinewidth=1truept
   \plotvskip 11
   \plottext{}
   \plot 24 0
   \plotmove 2 0
   \plotline 12 10
   \plotline 22 0
   \plotline 2 0
   {\plotlinewidth=0.05truept
   \plotmove 4 0
   \plotline 13 9
   \plotline 11 9
   \plotline 20 0
   \plotline 21 1
   \plotline 3 1
   \plotline 4 0
   \plotmove 6 0
   \plotline 14 8
   \plotline 10 8
   \plotline 18 0
   \plotline 20 2
   \plotline 4 2
   \plotline 6 0
   \plotmove 8 0
   \plotline 15 7
   \plotline 9 7
   \plotline 16 0
   \plotline 19 3
   \plotline 5 3
   \plotline 8 0
   \plotmove 10 0
   \plotline 16 6
   \plotline 8 6
   \plotline 14 0
   \plotline 18 4
   \plotline 6 4
   \plotline 10 0
   \plotmove 12 0
   \plotline 17 5
   \plotline 7 5
   \plotline 12 0 }
   \plotmove 14.5 8.5
   \plotline 8 2
   \plotline 21 2
   \plotmove 4 3
   \plotline 13 3
   \plotline 8.5 7.5
   \plotmove 7.5 -0.5
   \plotline 11 3
   \plotline 14.5 -0.5
   \plotmove 3.5 -0.5
   \plotline 7 3
   \plotline 10.5 -0.5
   \plotmove 2 1
   \plotline 11 1
   \plotline 6.5 5.5
   \plotmove 13.5 -0.5
   \plotline 17 3
   \plotline 20.5 -0.5
   \plotmove 13 1
   \plotline 15 3
   \plotline 17 1
   \plotline 13 1
   \plotmove 6 5
   \plotline 15 5
   \plotline 10.5 9.5
   \plotmove 13.5 9.5
   \plotline 11 7
   \plotline 16 7
   \plotmove 9.5 -0.5
   \plotline 13 3
   \plotline 16.5 -0.5
   \plotmove 22 1
   \plotline 17 1
   \plotline 19.5 3.5
   \plotmove 9 8
   \plotline 15 8
   \plotmove 17.5 -0.5
   \plotline 20.5 2.5
   \plotmove 6.5 -0.5
   \plotline 3.5 2.5
   \plotfmla{(10,0,0)}
   \plot 0.6 -0.2
   \plotfmla{(0,10,0)}
   \plot 23.4 -0.2
   \plotfmla{(0,0,10)}
   \plot 12 10.3
   \plotfmla{1}
   \plot 14 6.67
   \plot 15 6.33
   \plot 15 5.67
   \plot 16 5.33
   \plot 12 4.67
   \plot 13 4.33
   \plot 14 4.67
   \plot 15 4.33
   \plot 16 4.67
   \plot 17 4.33
   \plot 13 3.67
   \plot 14 3.33
   \plot 15 3.67
   \plot 16 3.33
   \plot 17 3.67
   \plot 18 3.33
   \plot 10 2.67
   \plot 12 2.67
   \plot 14 2.67
   \plot 16 2.67
   \plot 18 2.67
   \plotfmla{1}
   \plot 14 7.6
   \plot 13 6.6
   \plot 12 5.6
   \plot 13 5.8
   \plot 14 5.6
   \plot 11 4.6
   \plot 10 3.6
   \plot 11 3.8
   \plot 12 3.6
   \plot 9 2.6
   \plot 11 2.6
   \plot 13 2.6
   \plot 15 2.6
   \plot 17 2.6
   \plot 19 2.6
   \plotfmla{2}
   \plot 8 4.67
   \plot 9 4.33
   \plot 10 4.67
   \plot 9 3.67
   \plotfmla{2}
   \plot 9 6.6
   \plot 8 5.6
   \plot 9 5.8
   \plot 10 5.6
   \plot 7 4.6
   \plot 11 4.2
   \plot 6 3.6
   \plot 7 3.8
   \plot 8 3.6
   \plot 10 3.2
   \plot 11 3.4
   \plot 12 3.2
   \plotfmla{3}
   \plot 11 1.67
   \plot 10 0.67
   \plotfmla{3}
   \plot 11 2.2
   \plot 10 1.6
   \plot 12 1.6
   \plot 9 0.6
   \plot 11 0.6
   \plot 12 0.4
   \plot 13 0.6
   \plotfmla{4}
   \plot 6 0.67
   \plot 7 0.33
   \plot 8 0.67
   \plotfmla{4}
   \plot 7 2.6
   \plot 6 1.6
   \plot 7 1.8
   \plot 8 1.6
   \plot 5 0.6
   \plot 9 0.2
   \plotfmla{5}
   \plot 5 2.33
   \plot 6 2.67
   \plot 8 2.67
   \plot 5 1.67
   \plot 9 1.67
   \plotfmla{5}
   \plot 7 4.2
   \plot 6 3.2
   \plot 7 3.4
   \plot 8 3.2
   \plot 7 2.2
   \plot 9 2.2
   \plot 4 1.6
   \plot 6 1.2
   \plot 7 1.4
   \plot 8 1.2
   \plot 10 1.2
   \plotfmla{6}
   \plot 17 1.67
   \plot 16 0.67
   \plot 17 0.33
   \plot 18 0.67
   \plotfmla{6}
   \plot 17 2.2
   \plot 16 1.6
   \plot 18 1.6
   \plot 15 0.6
   \plot 19 0.6
   \plotfmla{7}
   \plot 15 1.67
   \plotfmla{7}
   \plot 15 2.2
   \plot 14 1.6
   \plot 16 1.2
   \plotfmla{8}
   \plot 10 7.33
   \plot 11 7.67
   \plot 10 6.67
   \plot 11 6.33
   \plot 12 6.67
   \plot 11 5.67
   \plotfmla{8}
   \plot 11 8.6
   \plot 12 7.6
   \plot 9 6.2
   \plot 13 6.2
   \plot 8 5.2
   \plot 9 5.4
   \plot 10 5.2
   \plot 12 5.2
   \plot 13 5.4
   \plot 14 5.2
   \plotfmla{9}
   \plot 13 7.67
   \plotfmla{9}
   \plot 13 8.6
   \plot 12 7.2
   \plot 14 7.2
   \plotfmla{10}
   \plot 13 1.7
   \plot 14 0.7
   \plotfmla{10}
   \plot 13 2.2
   \plot 12 1.2
   \plot 14 1.2
   \plot 11 0.2
   \plot 12 0.8
   \plot 13 0.2
   \plot 15 0.2
   \plotfmla{11}
   \plot 19 1.7
   \plotfmla{11}
   \plot 19 2.2
   \plot 18 1.2
   \plot 20 1.5
   \plotfmla{12}
   \plot 12 9.3
   \plot 12 8.7
   \plotfmla{12}
   \plot 11 8.2
   \plot 13 8.2
   \plotfmla{13}
   \plot 20 0.7
   \plot 21 0.3
   \plotfmla{13}
   \plot 20 1.15
   \plot 19 0.2
   \plotfmla{14}
   \plot 3 0.3
   \plot 4 0.7
   \plotfmla{14}
   \plot 4 1.2
   \plot 5 0.2
   \plotvskip 1
\plotfigend}
\botcaption{Figure 4}
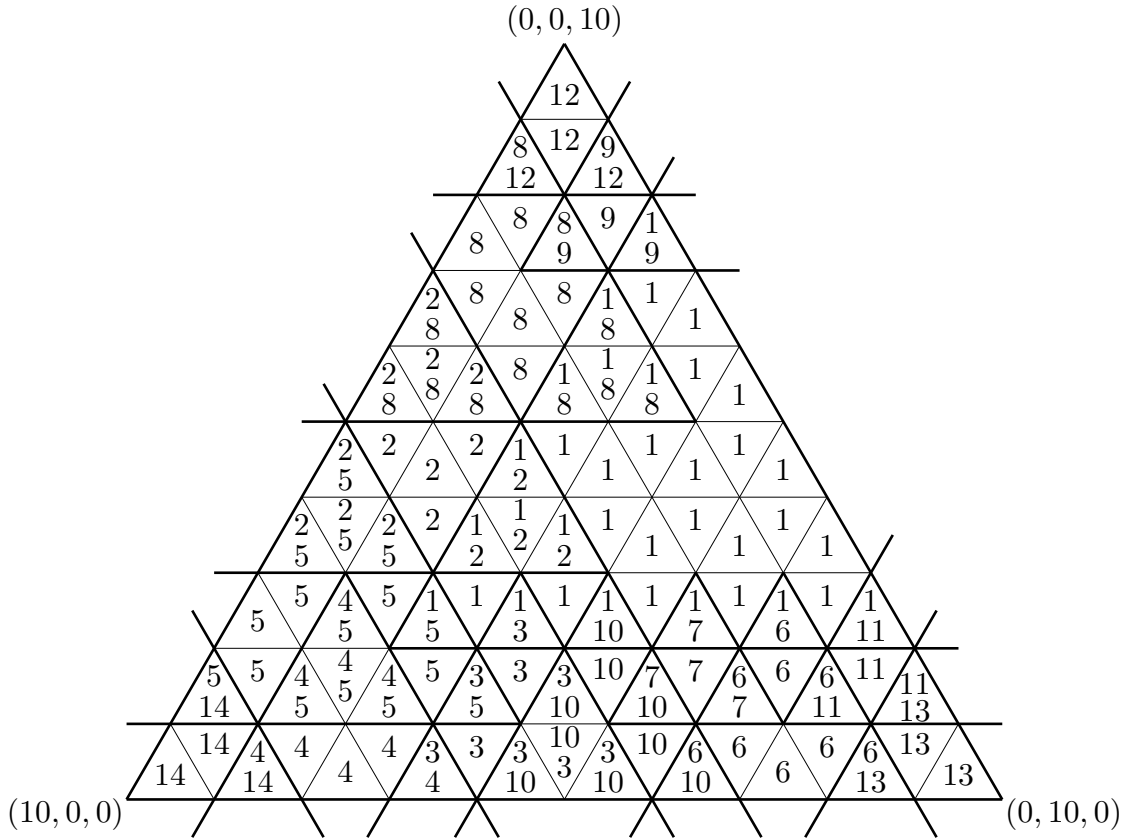
Coverage of one face of $\Rtetra 3{10}$ under $L'_7$.
\endcaption
\endinsert

In fact, since each $\vv_i$ has coordinates summing to at least $1$,
the smaller translates $\Rtetra 39 + \vv_i$, $i \le 14$, cover the
upper face of $\Rtetra 3{10}$.  For any $\vx$ in this upper face,
there is an $i$ such that $\vx \in \Rtetra 39 + \vv_i$, so
$\Rtetra31 + \vx \subseteq \Rtetra31 + \Rtetra 39 + \vv_i =
\Rtetra 3{10} + \vv_i$.  Since $\Rtetra 3{11}$ is the union of
$\Rtetra 3{10}$ and the sets $\Rtetra31 + \vx$ for $\vx$ in the upper
face of $\Rtetra 3{10}$, we have $\Rtetra 3{11} \subseteq
\Rtetra 3{10} + L'_7$, and hence $\Rtetra 3{11} + L'_7 \subseteq
\Rtetra 3{10} + L'_7$.

We now prove by induction that, for all integers $k \ge 10$, $\Rtetra
3k + L'_7 \subseteq \Rtetra 3{10} + L'_7$.  The case $k = 10$ is
trivial.  If it is true for $k$, then
$$\Rtetra 3{k+1} + L'_7 = \Rtetra 31 + \Rtetra 3k + L'_7
\subseteq \Rtetra 31 + \Rtetra 3{10} + L'_7 = \Rtetra 3{11} + L'_7
\subseteq \Rtetra 3{10} + L'_7,$$ so it is true for $k+1$.

Finally, for any $\vy \in \R^3$, there is a member $\vw$ of $L'_7$
such that the coordinates of $\vy-\vw$ are all positive (e.g.,
let $\vw$ be a large multiple of $-\vv_7$).  Then $\vy-\vw \in
\Rtetra 3k$ for some $k$, so $\vy \in \Rtetra 3k + L'_7$ and hence
$\vy \in \Rtetra 3{10} + L'_7$.  Therefore, $\Rtetra 3{10} + L'_7
= \R^3$. \QED

This covering has efficiency $\det(\vv_1,\vv_2,\vv_3)/\volume(\Rtetra 3{10})
= .504$.  Hence, Corollary~10(b) gives:

\proclaim{Corollary 18} For all $k$, there is a directed Cayley graph
of an Abelian group on three generators which has diameter $k$
and size at least $0.084k^3+O(k^2)$. \QNED

We can now use the method of (the real version of) Lemma~3 to get a
fundamental region $T'_7 \subseteq \Rtetra 3{10}$ for the lattice $L'_7$.
(Recall that any such region must have volume $\det(\vv_1,\vv_2,\vv_3) = 84$.)
To do this, just start with $\Rtetra 3{10}$, look at each of the
vectors $\vv_i$ ($i \le 14$) defined above, and delete those points of
$\Rtetra 3{10}$ which lie in $\Rtetra 3{10} + \vv_i$.  (We have $\vv_i
\succ \vzero$ for all $i$ if $\prec$ orders vectors primarily by the
sum of coefficients.)  What is left is the set shown in Figure~5; since
this is the union of $84$ unit cubes, we know that there is
no need to subtract further translates $\Rtetra 3{10} + \vw$.
This set was obtained independently by Fiduccia, Forcade, and
Zito~\cite{\FidForZit}.

\midinsert
\centerline{
\plotfigbegin
   \plotfigscales 19 11
   \plotvskip 25
   \plottext{}
   \plot 16 25
   \plotmove 5 6
   \plotline 10 11
   {\plotline 15 6 }
   \plotline 10 21
   \plotline 9 22
   {\plotline 9 24 }
   \plotline 8 21
   {\plotline 9 20 }
   \plotline 8 23
   {\plotline 7 24 }
   \plotline 9 24
   \plotline 8 25
   \plotline 7 24
   \plotline 7 18
   {\plotline 9 16 }
   \plotline 5 16
   {\plotline 7 14 }
   \plotline 5 12
   {\plotline 7 10 }
   \plotline 3 10
   {\plotline 5 8 }
   \plotline 3 6
   {\plotline 0 3 }
   \plotline 6 3
   {\plotline 6 5 }
   \plotline 5 2
   {\plotline 5 0 }
   \plotline 3 4
   {\plotline 3 2 }
   \plotline 1 2
   {\plotline 1 0 }
   \plotline 0 3
   \plotline 0 1
   \plotline 1 0
   \plotline 3 2
   \plotline 5 0
   \plotline 7 2
   {\plotline 7 6 }
   \plotline 8 1
   {\plotline 8 5 }
   \plotline 9 2
   {\plotline 9 6 }
   \plotline 11 0
   {\plotline 11 2 }
   \plotline 13 2
   {\plotline 13 6 }
   \plotline 15 0
   {\plotline 15 2 }
   \plotline 16 1
   \plotline 16 3
   {\plotline 15 4 }
   \plotline 15 2
   \plotline 14 3
   {\plotline 14 5 }
   \plotline 15 4
   \plotline 15 6
   \plotline 14 5
   \plotline 13 6
   \plotline 12 5
   {\plotline 11 6 }
   \plotline 12 3
   {\plotline 11 4 }
   \plotline 11 2
   \plotline 10 3
   {\plotline 10 5 }
   \plotline 11 4
   \plotline 11 6
   \plotline 10 5
   \plotline 9 6
   \plotline 8 5
   \plotline 7 6
   \plotline 6 5
   \plotline 5 6
   \plotline 5 8
   \plotline 7 10
   \plotline 7 14
   \plotline 9 16
   \plotline 9 20
   \plotline 10 21
\plotfigend}
\botcaption{Figure 5}
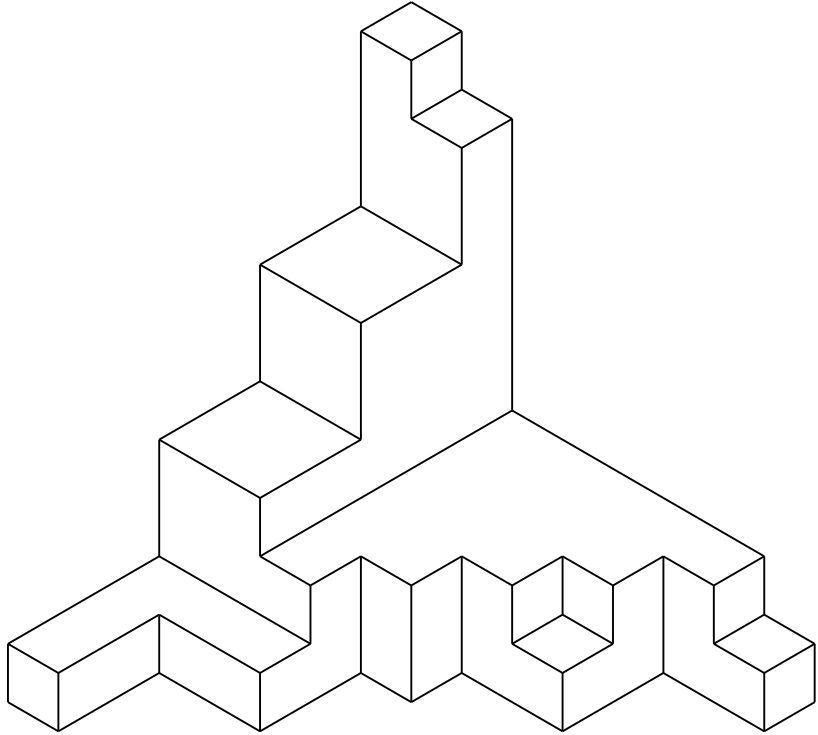
A subset of the tetrahedron $\Rtetra 3{10}$ which tiles space.
\endcaption
\endinsert

So this set $T'_7$ gives a lattice tiling of $\R^3$ using $L'_7$.
This tiling is quite unusual; the translates of $T'_7$ fit together
in a peculiar way, seeming to wind around each other.
One interesting fact is that each translate of $T'_7$ is adjacent to
(i.e., shares a boundary segment of positive area with) $28$ other
translates, a surprisingly high number.  ($T'_7$ itself is adjacent
to $T'_7 + \vv_i$ and $T'_7 - \vv_i$ for $i \le 14$.)

In many of the tilings we constructed explicitly, there was a polycube
fundamental region like $T'_7$, but there was also an alternate
fundamental region which was convex; for instance, for the optimal
covering of $\R^2$ by right triangles, one could use either
an L-tromino or a hexagon as the fundamental region.  Clearly
$L'_7$ has convex fundamental regions (e.g., its Voronoi regions),
but it turns out that they are unsuitable for the current problem:

\proclaim{Proposition 19} There is no convex fundamental region for
the lattice $L'_7$ included within the tetrahedron $\Rtetra 3{10}$.
\endproclaim

\demo{Proof} Since $T'_7$~gives a lattice tiling of~$\R^3$ by~$L'_7$,
every point of~$\R^3$, except for those lying on boundaries of the tiling,
can be translated by a vector in~$L'_7$ to a unique point of~$T'_7$.
In particular, if we look at the part of $\Rtetra 3{10}$ lying outside~$T'_7$,
then we can break it up into finitely many parts (in fact,
just cut it along the integer translates of the three coordinate planes)
which can be translated in a unique way by members of~$L'_7$ so as
to lie within~$T'_7$.

If one does this, one finds that there are parts of $T'_7$ which do
not get covered by translates of parts of $\Rtetra 3{10} \setminus
T'_7$.  Most of these uncovered parts look like inverted copies of~$\Rtetra 31$
(i.e., translates of $-\Rtetra 31$), although there are some larger ones.
In particular, the sets $(1,1,1)-\Rtetra 31$, $(8,1,1)-\Rtetra 31$,
$(1,8,1)-\Rtetra 31$, and $(1,1,8)-\Rtetra 31$ are not covered
by such translates.  This implies that each of those four sets is
disjoint (except for boundaries) from all of the translates
$\Rtetra 3{10} + \vw$ for $\vw \in L'_7 \setminus \{\vzero\}$.
It follows that any fundamental region for $L'_7$ included
within $\Rtetra 3{10}$ must include all four of these sets.

If the fundamental region is also convex, then it must contain any
convex combinations of points in those four sets; in particular,
it must include the sets $(3,1,1)-\Rtetra 31$ and $(1,3,3)-\Rtetra 31$.
But $(1,3,3) = (3,1,1) + \vv_1$, so the $\vv_1$\snug-translate
of the region overlaps the region itself in a set of positive volume,
which is impossible for a fundamental region of~$L'_7$.  Therefore, no
fundamental region of $L'_7$ within $\Rtetra 3{10}$ can be convex. \QED

There is no obvious reason why the lattice~$L'_7$ should be
exactly optimal for a lattice covering of~$\R^3$ by the
tetrahedron~$\Rtetra 3{10}$.  In the case of undirected graphs
on three generators, the lattices obtained for each~$k$ were not
optimal, but were closer and closer approximations to the
lattice~$L_{\text{bcc}}$, which does appear to be optimal.
One would expect something similar to occur in the directed case,
but it does not; the real lattice efficiencies in the last columns
of Tables 2 and~3 go up and down irregularly and do not (so far)
exceed the value~$.504$ attained by~$L'_7$.

Given this, it seems reasonable to examine~$L'_7$ and try to adjust it
slightly in order to improve its efficiency; there should be some
locally optimal lattice which $L'_7$~is an approximation to, and
we would like to find it.  Quite surprisingly, it turns out that
no adjustment is necessary.  Just before submitting the present paper,
the authors found a recent paper of Forcade and Lamoreaux~\cite{\ForLam}
proving this same result by a method slightly different from that
presented here.

\proclaim{Theorem 20 {\rm (Forcade and Lamoreaux \cite{\ForLam})}}
Among those lattices~$L$ for which
$\Rtetra 3{10} + L = \R^3$, the lattice~$L'_7$ is
locally optimal.
\endproclaim

\demo{Proof}
We use the same methods as for Theorem~16.  Recall the vectors
$\vv_1,\dots,\vv_{14}$ from Proposition~18.  The vectors
$\vv_1,\vv_2,\vv_3$ generate~$L'_7$;
a nearby lattice~$L$ will be generated by nearby vectors $\vv'_1 =
(a_1,b_1,c_1)$, $\vv'_2 = (a_2,b_2,c_2)$, and $\vv'_3 = (a_3,b_3,c_3)$.
Again concatenate $\vv'_1,\vv'_2,\vv'_3$ and $\vv_1,\vv_2,\vv_3$
to get $\vv'$ and~$\vv$ in~$\R^9$.
Let $F(\vv')$ be the determinant
of the matrix with rows $\vv'_1,\vv'_2,\vv'_3$; then we
have $F(\vv) = 84$, and we want to see that $F(\vv') < 84$
for any other~$\vv'$ near~$\vv$ for which the corresponding
lattice~$L$ satisfies $\Rtetra 3{10} + L = \R^3$.
We compute that the gradient of~$F$ at the point~$\vv$ is
$\vg = (-6,15,21,8,-6,14,12,12,0)$.

Referring back to Figure~4, we see that the point
$(1,1,8)$ is on the boundary of $\Rtetra 3{10} + \vv_i$ for
$i = 8,9,12$, as well as on the boundary of $\Rtetra 3{10}$~itself;
one can check that no other $L'_7$\snug-translate of~$\Rtetra 3{10}$
is near this point.  The nearby lattice~$L$ contains points
$\vzero$, $\vv'_8 = \vv'_1 + \vv'_2$, $\vv'_9 = 2 \vv'_1 + \vv'_2$,
and $\vv'_{12} = 2 \vv'_1 + \vv'_2 - \vv'_3$.  For any
small positive~$\eps$, the point $(a_1+a_2 - \eps, 2b_1 + b_2 - \eps,
2c_1+c_2 - c_3 - \eps)$, which is near $(1,1,8)$, will not be
in $\vv'_8 + \Rtetra 3{10}$ because its first coordinate is smaller
that that of~$\vv'_8$.  By looking at second and third coordinates
respectively, we see that this point is not in $\vv'_9 + \Rtetra 3{10}$
or $\vv'_{12} + \Rtetra 3{10}$ either.  Hence, in order to
have $\Rtetra 3{10} + L = \R^3$, this point must be in
$\Rtetra 3{10}$~itself, so we must have
$$a_1+a_2+2b_1 + b_2 +2c_1+c_2 - c_3 - 3 \eps \le 10.$$
Since $\eps$~can be arbitrarily small, we need
$$a_1+a_2+2b_1 + b_2 +2c_1+c_2 - c_3 \le 10$$ in order to have
$\Rtetra 3{10} + L = \R^3$.  So we have the constraint
$$\align
\vu_1 \vdot \vv' \le 10,\qquad \text{where }&\vu_1 = (1,2,2,1,1,1,0,0,-1). \\
\intertext{The same reasoning applied at the points
$(1,2,7)$, $(3,2,5)$, $(5,2,3)$, $(5,3,2)$, $(4,4,2)$, $(3,5,2)$,
$(2,6,2)$, $(1,7,2)$, $(8,1,1)$, $(6,3,1)$, $(3,6,1)$, and $(1,8,1)$
gives the constraints}
\allowdisplaybreak
\vu_2 \vdot \vv' \le 10,\qquad \text{where }&\vu_2 = (1,1,2,1,0,1,0,0,0), \\
\vu_3 \vdot \vv' \le 10,\qquad \text{where }&\vu_3 = (0,1,1,1,0,1,0,0,0), \\
\vu_4 \vdot \vv' \le 10,\qquad \text{where }&\vu_4 = (-1,1,0,1,0,1,0,0,0), \\
\vu_5 \vdot \vv' \le 10,\qquad \text{where }&\vu_5 = (-1,0,1,1,0,0,0,1,0), \\
\allowdisplaybreak
\vu_6 \vdot \vv' \le 10,\qquad \text{where }&\vu_6 = (0,-1,1,0,-1,0,1,1,0), \\
\vu_7 \vdot \vv' \le 10,\qquad \text{where }&\vu_7 = (-1,1,1,-1,0,0,1,1,0), \\
\vu_8 \vdot \vv' \le 10,\qquad \text{where }&\vu_8 = (1,0,1,0,-1,0,1,1,0), \\
\vu_9 \vdot \vv' \le 10,\qquad \text{where }&\vu_9 = (0,2,1,-1,-1,0,1,0,0), \\
\allowdisplaybreak
\vu_{10} \vdot \vv' \le 10,\qquad \text{where }&\vu_{10} =
(-2,-1,-1,0,0,1,1,1,0), \\
\vu_{11} \vdot \vv' \le 10,\qquad \text{where }&\vu_{11} =
(-1,0,-1,0,0,1,1,1,0), \\
\vu_{12} \vdot \vv' \le 10,\qquad \text{where }&\vu_{12} =
(-1,0,1,-1,-1,0,1,1,1), \\
\vu_{13} \vdot \vv' \le 10,\qquad \text{where }&\vu_{13} =
(0,1,2,-1,-1,-1,1,1,0).
\endalign$$

Again note that all thirteen of these inequalities are satisfied with equality
when $\vv'=\vv$.  Hence, they can be rewritten as
$\vu_i \vdot (\vv'-\vv) \le 0$ for $i =1,2,\dots,13$.

One can easily check that the vectors
$\vu_1,\vu_2,\vu_4,\vu_5,\vu_6,\vu_7,\vu_8,\vu_{10}$ are linearly
independent; their common null space is generated by
the vector $\vw = (0,0,0,1,1,-1,2,-1,1)$.
(The other five vectors~$\vu_i$ are also orthogonal
to~$\vw$, so they are linear combinations of the
eight listed above.)  Also, we have
$$\vg = \vu_1 + 4.8 \vu_2 + 6.4 \vu_4 + \vu_5 + 1.6 \vu_6 + 3.2 \vu_7 +
3.4 \vu_8 + \vu_9 + 1.8 \vu_{10} + \vu_{12}.$$

Let $C$~be the closed cone consisting of all vectors~$\vt$ in the
subspace spanned by $\vu_1,\dots,\vu_{13}$
such that $\vu_i \vdot \vt \le 0$ for all $i \le 13$.  Then
the above equations imply that $\vg \cdot \vt \le 0$ for all
$\vt$ in~$C$, and equality can hold only when $\vt = \vzero$.
Hence, as in Theorem~16, there is a neighborhood~$U$ of~$\vg$ such that,
for any $\vg'$ in~$U$ and any nonzero $\vt$ in~$C$,
$\vg' \vdot \vt < 0$.

We can compute that, for any real number $r$, the determinant for the lattice
given by $\vv + r \vw$ is
$F(\vv + r \vw) = 84-12r^2$.  Clearly this is at most~$84$,
with equality only when $r = 0$.

Now, any vector $\vv'$
close to~$\vv$ can be expressed as $\vv + \vt_1 + \vt_2$
where $\vt_1$ is a small multiple of~$\vw$ and
$\vt_2$~is a small linear combination of the vectors
$\vu_1,\dots,\vu_{13}$.
The reasoning from Theorem~16 shows that
$\vt_2$~must be in~$C$ if $\Rtetra 3{10} + L = \R^3$.
Also as in that Theorem, we find that
$F(\vv + \vt_1) \le F(\vv)$
with equality only when $\vt_1 = \vzero$, and
$F(\vv + \vt_1 + \vt_2) \le F(\vv + \vt_1)$ with equality
only when $\vt_2 = \vzero$.
Therefore, $F(\vv') \le F(\vv)$, with equality holding only when
$\vv' = \vv$.  So $\vv$~gives a local maximum of~$F$, as desired.
\QED

This and the computational evidence make it plausible that $L'_7$~actually
gives an optimal lattice covering of~$\R^3$ by~$\Rtetra 3{10}$,
and hence that the asymptotic formula in Corollary~18 is optimal.

\subhead More than three generators \endsubhead
In higher dimensions, analogues of many of the preceding constructions
exist, but they do not produce lattice coverings as efficient as
one would hope for.

For lattice coverings with the $d$-dimensional dual cube, one can
use the $d$-dimensional body-centered cubic lattice (the set of
vectors in $\Z^d$ whose coordinates are all odd or all even).  By the
same argument as for the three-dimensional case, this lattice gives
a lattice covering of $\R^d$ by $\Rocta d{d/2}$.  The efficiency
of this covering is $2^{d-1}/\volume(\Rocta d{d/2}) = 2^{d-1}d!/d^d$,
which is $2^{d-1}$ times the efficiency of the covering using the
ordinary cubic lattice $\Z^d$.

As usual, the Cayley graph corresponding to this lattice is a twisted
toroidal mesh.  For a given number $m$, one can connect
the elements of $\Z_{2m}^{d-1}\times \Z_m$ as in an ordinary
toroidal mesh, except that the wraparound connections for the
last coordinate are twisted along all of the other coordinates:
$(x_1,\dots,x_{d-1},m{-}1)$ is connected to
$(x_1{\pm}m,\dots,x_{d-1}{\pm}m,0)$.  This gives a graph of
diameter $\lfloor dm/2 \rfloor$ and size $2^{d-1}m^d$, which is
about $2^{d-1}$ times as large as the best ordinary toroidal mesh
of this diameter.

One can optimize this slightly.  Given the dimension~$d$ and the
desired diameter~$k$, let $q$ and~$r$ be the quotient and
remainder when $2k+1$ is divided by~$d$; we assume $k$~is large enough
that $q > 0$.  Then a good lattice $L$ to use is the body-centered
cubic lattice above scaled up by a factor~$q+1$ in each of the
first $r$~coordinates and a factor~$q$ in the remaining $d-r$ coordinates.
The resulting $\Z^d/L$ is isomorphic to
$(\Z_{2q+2}^r \times \Z_{2q}^{d-r})/H$ with the canonical generators,
where $H$ is the two-element subgroup $\{\vzero,(q{+}1,\dots,q{+}1,q,\dots,q)\}$
(there are $r$ $q{+}1$\snug's); it can be laid out as a twisted
toroidal mesh on $\Z_{2q+2}^r \times \Z_{2q}^{d-r-1} \times \Z_q$
or on $\Z_{2q+2}^{r-1} \times \Z_{2q}^{d-r} \times \Z_{q+1}$.
If $q$~is even and $r > 0$, this Cayley graph is isomorphic
to that of $\Z_{q+1} \times \Z_{2q+2}^{r-1} \times \Z_{2q}^{d-r}$
with the generators $\ve_2,\dots,\ve_d$ and $(1,\dots,1,q,\dots,q)$
with $r$~$1$\snug's;  if $q$~is odd or $r=0$, then it is isomorphic
to the Cayley graph of $\Z_{2q+2}^r \times \Z_{2q}^{d-r-1} \times \Z_q$
with generators $\ve_1,\dots,\ve_{d-1}$ and $(q{+}1,\dots,q{+}1,1,\dots,1)$
with $r$~$q{+}1$\snug's.  The size of this graph is slightly larger than
the size of the cyclic Cayley graph constructed by Chen and
Jia~\cite{\CheJia}, but the ratio of the two sizes tends to~$1$
for large~$k$.

For the directed case, we must consider lattice coverings by
$d$-simplices; as usual, by affine invariance, it doesn't matter which
simplex is used.  One can show that a lattice for covering with
a given $d$-simplex is given by the following generating vectors:
for each face of the simplex, take a vector which is twice the
vector from the centroid of the simplex to the centroid of that face.
(This gives $d+1$ vectors, but they sum to $\vzero$, so just take
$d$ of them.)  The efficiency of this covering works out to be
$d!2^d/(d^d(d+1))$.

Unfortunately, in both cases, the efficiency decreases exponentially
with $d$: by Stirling's formula,
$${2^{d-1}d!\over d^d} \sim \sqrt{\pi d \over 2} \left({2 \over e}\right)^d$$
and
$${d!2^d\over d^d(d+1)} \sim \sqrt{2\pi \over d} \left({2 \over e}\right)^d.$$
This seems to be the case for all known explicitly constructed
lattice coverings by these shapes (and by spheres).

On the other hand, in 1959 Rogers~\cite{\Rogers} gave a nonconstructive
proof that there exist much more efficient lattice coverings by
these shapes (or by any convex body) in high dimensions; and he gave
an even better result for the case of spheres.  More recently
Gritzmann~\cite{\Gritzmann} extended the latter result to apply to
any convex body with a sufficient number of mutually orthogonal
hyperplanes of symmetry.  (The number required is quite small: only
$\lfloor \log_2 \ln d \rfloor + 5$.)  Gritzmann's result states
that there is a constant $c$ (not depending on $d$ or on the
convex body) such that, for any convex body $K$ in $\R^d$ with
the above number of mutually orthogonal planes of symmetry,
there is a lattice covering of $\R^d$ by $K$ with density at most
$c d (\ln d)^{1+\log_2 e}$.

The regular dual $d$\snug-cube and the regular $d$\snug-simplex do have the
required symmetry planes for large enough~$d$.  This is clear
for the dual $d$\snug-cube; it has the same $d$~orthogonal planes
of symmetry as the $d$\snug-cube it is dual to.  For the regular
$d$\snug-simplex, note that the perpendicular bisector of an edge
is a hyperplane of symmetry, and that edges which do not share
a vertex have orthogonal directions (the easiest way to see this is
to look at the regular $d$\snug-simplex in $\R^{d+1}$ whose
vertices are $\ve_1,\dots,\ve_{d+1}$, and take dot products), so one
can find $\lceil d/2 \rceil$ mutually orthogonal hyperplanes of symmetry.
Therefore, we get lattice coverings of the specified density
for large $d$, and by adjusting the constant~$c$ we can make
the bound apply for all~$d$ (for these two particular shapes).
Therefore, letting $\bar c = c^{-1}$, we can use Corollary~10 to get:

\proclaim{Theorem 21} There is a constant $\bar c > 0$ (not depending on
$d$ or $k$) such that, for any fixed $d > 1$ and for all $k$, there
exist undirected Cayley graphs of Abelian groups on $d$ generators
having diameter $\le k$ and size at least
$${2^d \bar c \over d!d(\ln d)^{1+\log_2 e}} k^d + O(k^{d-1}),$$
and there exist directed Cayley graphs of Abelian groups on $d$
generators having diameter $\le k$ and size at least
$${\bar c \over d!d(\ln d)^{1+\log_2 e}} k^d + O(k^{d-1}).$$
\QNED

The coverings produced by this method are probably fairly strange.
We seem to have run into this already in three dimensions, for
the directed case; for the undirected case it apparently happens later.

\subhead Layouts with short wires \endsubhead
The obvious way to lay out a toroidal mesh is as a rectangular array
with mesh connections between adjacent nodes in the array and with long
wires connecting opposite ends of the array; these long wires may cause
communications delays.  However, there is a standard trick for
rearranging the layout so as to remove the need for long wires.  In the
one-dimensional case, instead of placing the nodes in the order
$1,2,3,\dots,n$ (where $i$ is connected to $i+1$ for
$i=1,2,\dots,n{-}1$ and $n$ is connected to $1$), one can place them in
the order $1,n,2,n{-}1,3,n{-}2,\dots$; then the maximum required wire
length is only twice the mesh spacing.  In higher dimensions, one can
apply the same trick to each dimension separately, and again the
required wire length is twice the mesh spacing.

It is not immediately obvious that this interleaving trick can be applied
to twisted toroidal meshes; a simple interleaving in each dimension
would not make the twisted cross-connections short.  But it is possible
to get short-wire layouts for the twisted meshes in a similar way.
One approach is to perform the interleaving {\sl twice} on the
long dimensions of the mesh; for instance, if the mesh has length $16$
in one of the long dimensions, then the nodes would be arranged in
the order $$1,9,16,8,2,10,15,7,3,11,14,6,4,12,13,5.$$ Then wires in this
dimension would have length at most $4$ times the mesh spacing.
Now, when one does a single interleaving on the short dimension,
the twisted cross-connections become short as well.

\midinsert
\hbox to\hsize{
\hfill
\plotfigbegin
   \plotfigscales 18 18
   \plotvskip 10.5
   \plotclosed
   \dimen0=1sp
   \plotcurrx=1\plotunitx
   \loop
      {\dimen2=1sp
      \plotcurry=2\plotunity
      \loop
         \plothere
         \advance\plotcurry by\plotunity
      \ifdim\dimen2<8sp
         \advance\dimen2 by 1sp
      \repeat}
      \advance\plotcurrx by\plotunitx
   \ifdim\dimen0<4sp
      \advance\dimen0 by 1sp
   \repeat
   \plotlinewidth=1truept
   \plotmove 1 1.5
   \plotline 1 9.5
   \plotmove 2 1.5
   \plotline 2 9.5
   \plotmove 3 1.5
   \plotline 3 9.5
   \plotmove 4 1.5
   \plotline 4 9.5
   \plotlinewidth=0.2truept
   \plotmove 0.5 2
   \plotline 4.5 2
   \plotmove 0.5 3
   \plotline 4.5 3
   \plotmove 0.5 4
   \plotline 4.5 4
   \plotmove 0.5 5
   \plotline 4.5 5
   \plotmove 0.5 6
   \plotline 4.5 6
   \plotmove 0.5 7
   \plotline 4.5 7
   \plotmove 0.5 8
   \plotline 4.5 8
   \plotmove 0.5 9
   \plotline 4.5 9
   \plottext{A}
   \plot 1 1
   \plot 1 10
   \plottext{B}
   \plot 2 1
   \plot 2 10
   \plottext{C}
   \plot 3 1
   \plot 3 10
   \plottext{D}
   \plot 4 1
   \plot 4 10
   \plottext{1}
   \plot 0 2
   \plot 5 6
   \plottext{2}
   \plot 0 3
   \plot 5 7
   \plottext{3}
   \plot 0 4
   \plot 5 8
   \plottext{4}
   \plot 0 5
   \plot 5 9
   \plottext{5}
   \plot 0 6
   \plot 5 2
   \plottext{6}
   \plot 0 7
   \plot 5 3
   \plottext{7}
   \plot 0 8
   \plot 5 4
   \plottext{8}
   \plot 0 9
   \plot 5 5
   \plottext{(a)}
   \plot 2.5 0
\plotfigend
\hfill
\plotfigbegin
   \plotfigscales 18 18
   \plotvskip 10.5
   \plotclosed
   \dimen0=1sp
   \plotcurrx=1\plotunitx
   \loop
      {\dimen2=1sp
      \plotcurry=2\plotunity
      \loop
         \plothere
         \advance\plotcurry by\plotunity
      \ifdim\dimen2<8sp
         \advance\dimen2 by 1sp
      \repeat}
      \advance\plotcurrx by\plotunitx
   \ifdim\dimen0<4sp
      \advance\dimen0 by 1sp
   \repeat
   \plotlinewidth=1truept
   \plotmove 1 1.5
   \plotline 1 9.5
   \plotmove 2 1.5
   \plotline 2 9.5
   \plotmove 3 1.5
   \plotline 3 9.5
   \plotmove 4 1.5
   \plotline 4 9.5
   \plotlinewidth=0.2truept
   \plotmove 0.5 2
   \plotline 1 2
   \plotline 1.5 1.5
   \plotmove 1.5 9.5
   \plotline 4 7
   \plotline 4.5 7
   \plotmove 0.5 3
   \plotline 1 3
   \plotline 2.5 1.5
   \plotmove 2.5 9.5
   \plotline 4 8
   \plotline 4.5 8
   \plotmove 0.5 4
   \plotline 1 4
   \plotline 3.5 1.5
   \plotmove 3.5 9.5
   \plotline 4 9
   \plotline 4.5 9
   \plotmove 0.5 5
   \plotline 1 5
   \plotline 4 2
   \plotline 4.5 2
   \plotmove 0.5 6
   \plotline 1 6
   \plotline 4 3
   \plotline 4.5 3
   \plotmove 0.5 7
   \plotline 1 7
   \plotline 4 4
   \plotline 4.5 4
   \plotmove 0.5 8
   \plotline 1 8
   \plotline 4 5
   \plotline 4.5 5
   \plotmove 0.5 9
   \plotline 1 9
   \plotline 4 6
   \plotline 4.5 6
   \plottext{A}
   \plot 1 1
   \plot 1 10
   \plottext{B}
   \plot 2 1
   \plot 2 10
   \plottext{C}
   \plot 3 1
   \plot 3 10
   \plottext{D}
   \plot 4 1
   \plot 4 10
   \plottext{a}
   \plot 1.5 1
   \plot 1.5 10
   \plottext{b}
   \plot 2.5 1
   \plot 2.5 10
   \plottext{c}
   \plot 3.5 1
   \plot 3.5 10
   \plottext{1}
   \plot 0 2
   \plot 5 3
   \plottext{2}
   \plot 0 3
   \plot 5 4
   \plottext{3}
   \plot 0 4
   \plot 5 5
   \plottext{4}
   \plot 0 5
   \plot 5 6
   \plottext{5}
   \plot 0 6
   \plot 5 7
   \plottext{6}
   \plot 0 7
   \plot 5 8
   \plottext{7}
   \plot 0 8
   \plot 5 9
   \plottext{8}
   \plot 0 9
   \plot 5 2
   \plottext{(b)}
   \plot 2.5 0
\plotfigend
\hfill
\plotfigbegin
   \def\1{\plotPSbegin currentpoint newpath moveto 0.3 0.5 0.3 1.5 0 2
         rcurveto stroke \plotPSend \plotrmove 0 2 }
   \def\2{\plotPSbegin currentpoint newpath moveto -0.3 -0.5 -0.3 -1.5 0 -2
         rcurveto stroke \plotPSend \plotrmove 0 -2 }
   \def\3{\plotPSbegin currentpoint newpath moveto -0.3 0.7 -1.3 1.7 -2 2
         rcurveto stroke \plotPSend \plotrmove -2 2 }
   \def\4{\plotPSbegin currentpoint newpath moveto 0.3 -0.7 1.3 -1.7 2 -2
         rcurveto stroke \plotPSend \plotrmove 2 -2 }
   \def\5{\plotPSbegin currentpoint newpath moveto 0.7 0.3 1.7 1.3 2 2
         rcurveto stroke \plotPSend \plotrmove 2 2 }
   \def\6{\plotPSbegin currentpoint newpath moveto -0.7 -0.3 -1.7 -1.3 -2 -2
         rcurveto stroke \plotPSend \plotrmove -2 -2 }
   \plotfigscales 18 18
   \plotvskip 10.5
   \plothphant 5
   \plotclosed
   \dimen0=1sp
   \plotcurrx=1\plotunitx
   \loop
      {\dimen2=1sp
      \plotcurry=2\plotunity
      \loop
         \plothere
         \advance\plotcurry by\plotunity
      \ifdim\dimen2<8sp
         \advance\dimen2 by 1sp
      \repeat}
      \advance\plotcurrx by\plotunitx
   \ifdim\dimen0<4sp
      \advance\dimen0 by 1sp
   \repeat
   \plotlinewidth=1truept
   \plotmove 1 2
   \1
   \1
   \1
   \plotrline 0 1
   \2
   \2
   \2
   \plotrline 0 -1
   \plotmove 2 2
   \1
   \1
   \1
   \plotrline 0 1
   \2
   \2
   \2
   \plotrline 0 -1
   \plotmove 3 2
   \1
   \1
   \1
   \plotrline 0 1
   \2
   \2
   \2
   \plotrline 0 -1
   \plotmove 4 2
   \1
   \1
   \1
   \plotrline 0 1
   \2
   \2
   \2
   \plotrline 0 -1
   \plotlinewidth=0.2truept
   \plotmove 1 2
   \plotrline 1 2
   \5
   \plotrline -1 2
   \plotrline -2 1
   \plotrline 1 -2
   \4
   \plotrline -1 -2
   \plotrline -2 -1
   \plotmove 1 4
   \plotrline 1 2
   \5
   \plotrline -1 1
   \6
   \plotrline 1 -2
   \4
   \plotrline -1 -1
   \3
   \plotmove 1 6
   \plotrline 1 2
   \plotrline 2 1
   \plotrline -1 -2
   \6
   \plotrline 1 -2
   \plotrline 2 -1
   \plotrline -1 2
   \3
   \plotmove 1 8
   \plotrline 1 1
   \4
   \plotrline -1 -2
   \6
   \plotrline 1 -1
   \5
   \plotrline -1 2
   \3
   \plottext{A}
   \plot 1 1
   \plottext{B}
   \plot 3 1
   \plottext{C}
   \plot 4 1
   \plottext{D}
   \plot 2 1
   \plottext{1}
   \plot 0 2
   \plottext{2}
   \plot 0 4
   \plottext{3}
   \plot 0 6
   \plottext{4}
   \plot 0 8
   \plottext{5}
   \plot 0 9
   \plottext{6}
   \plot 0 7
   \plottext{7}
   \plot 0 5
   \plottext{8}
   \plot 0 3
   \plottext{(c)}
   \plot 2.5 0
\plotfigend
\hfill}
\botcaption{Figure 6}
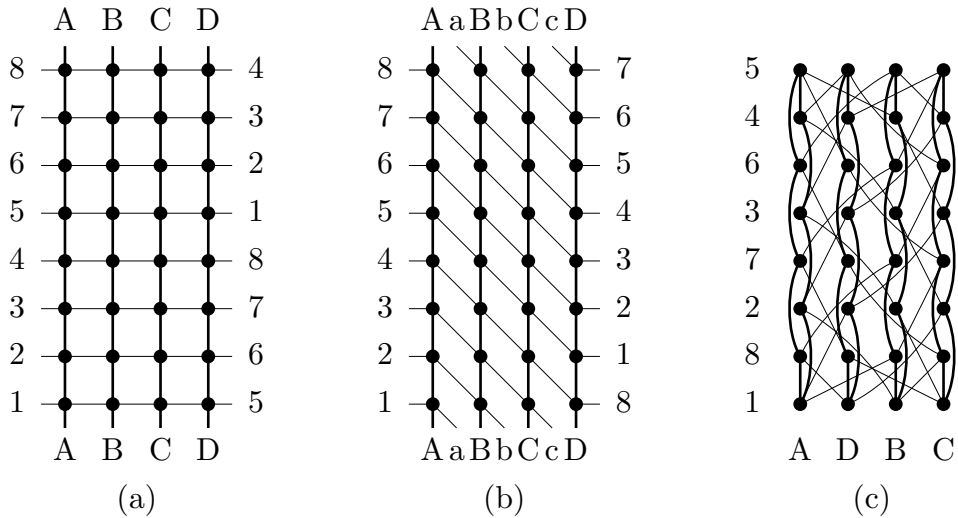
Short-wire layout for a twisted toroidal mesh.
\endcaption
\endinsert

Another method is shown in Figure~6.  Here the idea is to modify the
original arrangement (a) by shearing the mesh (rotating the $i$\snug'th
level in the short dimension by $i-1$ units in each of the long
dimensions), as shown in (b), so that the twisted cross-connections
become (almost) straight.  Then one can do an ordinary interleaving in
each dimension to get the result shown in (c).  This gives a layout in
which the maximum wire length is $2\sqrt{d}$ times the mesh spacing.

Some of the other Abelian Cayley graphs we have considered can be
treated similarly, especially the ones which differ only slightly
from twisted toroidal meshes.  For the optimal two-generator
undirected Abelian Cayley graph, if one starts with the almost-rectangular
layout shown in Figure~1 (a $(k{+}1)\times(k{+}1)$ square next to a
$k\times k$ square), and performs a shear as in Figure~6, then the
result is a $2k \times (k{+}1)$ rectangle with one node left over;
this can then be interleaved to get a short-wire layout.  A more difficult
case is the two-generator directed graph from Theorem~13 and Corollary~14;
here one can start with the natural L-shaped layout and perform
shears on separate parts to obtain a rectangle with dimensions
$(a+2b) \times a$ (made up of three subrectangles with different
shear patterns) where the necessary cross-connections are almost
straight across, and hence interleaving will give a good layout.

\subhead Generators of order 2 \endsubhead
The undirected Cayley graphs produced so far all have even degree (twice
the number of generators).  If one is interested in undirected Cayley
graphs of odd degree, one will have to use a generator of order~2.

Using $d$~unrestricted generators plus one order-2 generator, one can get
an undirected Abelian Cayley graph of a given diameter which is about
twice as large as one can get using $d$~unrestricted generators alone.
More precisely, if $n_a(d,k)$ is the size of the largest Abelian Cayley
graph of diameter~$k$ using $d$~generators and $n_a^+(d,k)$ is the
size of the largest such graph using $d$~generators plus one order-2
generator, then
$$2 n_a(d,k-1) \le n_a^+(d,k) \le 2 n_a(d,k).$$
To see this, first let $G$~be generated by $g_1,\dots,g_d$ and~$\rho$
where $\rho$~has order~2.  If $G$~has diameter at most~$k$ using these
generators, and $H$~is the subgroup of size~2 generated by~$\rho$,
then $G/H$~is generated by the images $g_i+H$ for $1 \le i \le d$
with diameter at most~$k$, so $|G/H| \le n_a(d,k)$, so
$|G| \le 2n_a(d,k)$; hence, $n_a^+(d,k) \le 2n_a(d,k)$.
On the other hand, if $G$~is generated by $g_1,\dots,g_d$ with
diameter at most~$k-1$, then $G \times \Z_2$ is generated by
$(g_i,0)$ for $i=1,\dots,d$ and $(0,1)$, and has diameter
at most~$k$ using these generators; this shows that
$2 n_a(d,k-1) \le n_a^+(d,k)$.

We can also study $n_a^+(d,k)$ using the same methods that were used
for $n_a(d,k)$.  The appropriate free (universal) group to use here
is the group $\Z^d \times \Z_2$, with the canonical generators
$(\ve_i,0)$ for $i=1,\dots,d$ and $(\vzero,1)$.  The set of elements of
this group which can be written as a word of length at most~$k$ in the
generators is precisely $W_k = (\octa dk \times {0}) \cup (\octa d{k-1} \times
{1})$.  For any Abelian group~$G$ with generators $g_1,\dots,g_d$ and~$\rho$
($\rho$~of order~2), there is a unique homomorphism from
$\Z^d \times \Z_2$ to~$G$ taking $(\ve_i,0)$ to~$g_i$ and
$(\vzero,1)$ to~$\rho$; the Cayley graph of~$G$ using these generators
has diameter at most~$k$ if and only if the homomorphism maps~$W_k$
onto~$G$.  So the obvious upper limit for the size of~$G$
is $|W_k| = |\octa dk| + |\octa d{k-1}|$.

We are now led to study quotient groups $(\Z^d \times \Z_2)/N$
where $N$~is a (normal) subgroup of $\Z^d \times \Z_2$ of finite index;
we want such an~$N$ of largest possible index such that
$W_k +N = \Z^d \times \Z_2$.

One simple possibility is that $N \subseteq \Z^d \times \{0\}$; in this case
the resulting group is just $(\Z^d / N_0) \times \Z_2$
where $N = N_0 \times \{0\}$.  It is easy to see that the diameter
of this group is precisely one more than the diameter of~$\Z^d/N_0$
using the canonical $d$~generators.

Note that $N_0$~is a $d$\snug-dimensional lattice; let $\vv_1,\dots,\vv_d$
be a list of generators for this lattice.  Now let
$N'$ be the subgroup of $\Z^d \times \Z_2$ generated by
$(\vv_i,1)$ for $i=1,\dots,d$.  Then we have
$$\indx \Z^d \times \Z_2 : N'| = 2\indx \Z^d : N_0 |
= \indx \Z^d \times \Z_2 : N|.$$
Furthermore, the diameter of $(\Z^d \times \Z_2)/N'$ is at most
one more than the diameter of~$\Z^d/N_0$, which means that it is
no larger than the diameter of $(\Z^d \times \Z_2)/N$.

This shows that, when trying to determine
$n_a^+(d,k)$, we may restrict ourselves to studying
subgroups~$N$ of $\Z^d \times \Z_2$ of finite index which are {\sl not}\/
included in~$\Z^d \times \{0\}$.

So choose $\vg \in \Z^d$ such that $(\vg,1) \in N$.
The subgroup $N \cap (\Z^d \times \{0\})$ is of index~2 in~$N$
and hence of finite index in~$\Z^d \times \Z_2$.
So we have $N \cap (\Z^d \times \{0\}) = L \times \{0\}$ for
some $d$\snug-dimensional lattice~$L$.  Note that
$(2\vg,0) = 2(\vg,1) \in N$, so
$2\vg \in L$.  (Normally $\vg$~will not be in~$L$; if $\vg \in L$, then
$(\vg,0) \in N$, so $(\vzero,1) \in N$, so the order-2 generator
collapses to the identity in the quotient group.)
Also, we have $\indx \Z^d \times \Z_2 : N| = \indx \Z^d : L|$.

It is now easy to see that
$$(W_k + N) \cap (\Z^d \times \{0\}) =
((\octa dk + L) \cup (\octa d{k-1} + \vg + L)) \times \{0\}.$$
Hence, in order to have $W_k + N = \Z^d \times \Z_2$, it is necessary
to have
$$(\octa dk + L) \cup (\octa d{k-1} + \vg + L) = \Z^d.\tag{$*$}$$
This necessary condition is also sufficient, because
$$\align
(W_k + N) \cap (\Z^d \times \{1\})
&= ((\octa d{k-1} + L) \cup (\octa dk + \vg + L)) \times \{1\} \\
&= (((\octa dk + L) \cup (\octa d{k-1} + \vg + L)) + \vg) \times \{1\}.
\endalign$$
So our goal is to find such a lattice~$L$ and extra generator~$\vg$
(with $2\vg \in L$) so that ($*$)~is satisfied and
$\indx \Z^d : L |$ is as large as possible.

We are now ready to consider specific values of~$d$.  As usual, the case
$d=1$ is easy.  The maximal possible value of $\indx \Z : L |$ is
$|W_k| = 4k$, and this value is attained by letting $L$~be generated
by the element~$4k$, with $\vg = 2k$.  This leads to the cyclic
Cayley graph on the group~$\Z_{4k}$ with unrestricted generator~$1$
and order-2 generator~$2k$.

For $d=2$ we have a situation very similar to that in Figure~1
(lattice coverings with Aztec diamonds), but not identical because
we must use two different shapes.  The upper bound on
$\indx \Z^2 : L |$ is $|W_k| = 4k^2 + 2$.  For $k=1$ this bound is
actually attainable; it leads to the Cayley graph from~$\Z_6$
with unrestricted generators $1$ and~$2$ and order-2 generator~$3$.
But for $k>1$ the pieces $\octa 2k$ and~$\octa 2{k-1}$ do not
fit together well enough to give a perfect tiling.  The best one
can do is the lattice~$L$ generated by $(2k+1,1)$ and $(-1,2k-1)$,
with the extra generator~$\vg = (k,k)$, as shown in Figure~7.

\midinsert
\centerline{
\plotfigbegin
   \plotfigscale 18
   \plotvskip 17
   \plottext{}
   \plot 23 17
   \def\1{\plotrline 0 1 }
   \def\2{\plotrline 0 -1 }
   \def\3{\plotrline 1 0 }
   \plotmove 0 4 \2\3\2\3\3\2\3\2\3\3\3{\1}\3
   \plotmove 0 8 \2\3\3\2\3\2\3{\2}\3\3{\1}\3\2\3\2\3\3\2\3\2\3{\2}\3\3{\1}%
         \3\2\3
   \plotmove 0 12 \3\2\3\2\3{\2}\3\3{\1}\3\2\3\2\3\3\2\3\2\3{\2}\3\3{\1}\3%
         \2\3\2\3\3\2\3\2\3{\2}\3\3{\1}\3\2\3\2
   \plotmove 0 17 \2\3\2\3{\2}\3\3{\1}\3\2\3\2\3\3\2\3\2\3{\2}\3\3{\1}\3\2%
         \3\2\3\3\2\3\2\3{\2}\3\3{\1}\3\2\3\2\3
   \plotmove 7 17 \3\2\3{\2}\3\3{\1}\3\2\3\2\3\3\2\3\2\3{\2}\3\3{\1}\3\2\3%
         \2\3\3\2
   \plotmove 15 17 \3{\2}\3\3\3\2\3\2\3\3\2\3\2
   \plotmove 0 11 \1\1\3\1\3\3\1\1\3\3\1\3
   \plotmove 0 5 \1\3\1\1\3\1\3\3\1\1\3\3\1\3\1\1\3\1\3\3\1\1\3\3
   \plotmove 0 0 \3\1\3\1\1\3\1\3\3\1\1\3\3\1\3\1\1\3\1\3\3\1\1\3\3\1\3\1\1%
         \3\1\3\3\1
   \plotmove 6 0 \1\3\3\1\3\1\1\3\1\3\3\1\1\3\3\1\3\1\1\3\1\3\3\1\1\3\3\1\3%
         \1\1\3\1\3
   \plotmove 11 0 \3\3\1\1\3\3\1\3\1\1\3\1\3\3\1\1\3\3\1\3\1\1\3\1
   \plotmove 17 0 \3\1\3\3\1\1\3\3\1\3\1\1
   \plotclosed
   \dimen0=1sp
   \plotcurrx=0.5\plotunitx
   \loop
      {\dimen2=1sp
      \plotcurry=0.5\plotunity
      \loop
         \plothere
         \advance\plotcurry by\plotunity
      \ifdim\dimen2<17sp
         \advance\dimen2 by 1sp
      \repeat}
      \advance\plotcurrx by\plotunitx
   \ifdim\dimen0<23sp
      \advance\dimen0 by 1sp
   \repeat
   \plotlinewidth=0.1truept
   \plotmove 11.2 8.8
   \plotline 17.8 8.8
   \plotline 17.8 4.8
   \plotline 18.8 4.8
   \plotline 18.8 4.2
   \plotline 11.2 4.2
   \plotline 11.2 8.8
\plotfigend}
\botcaption{Figure 7}
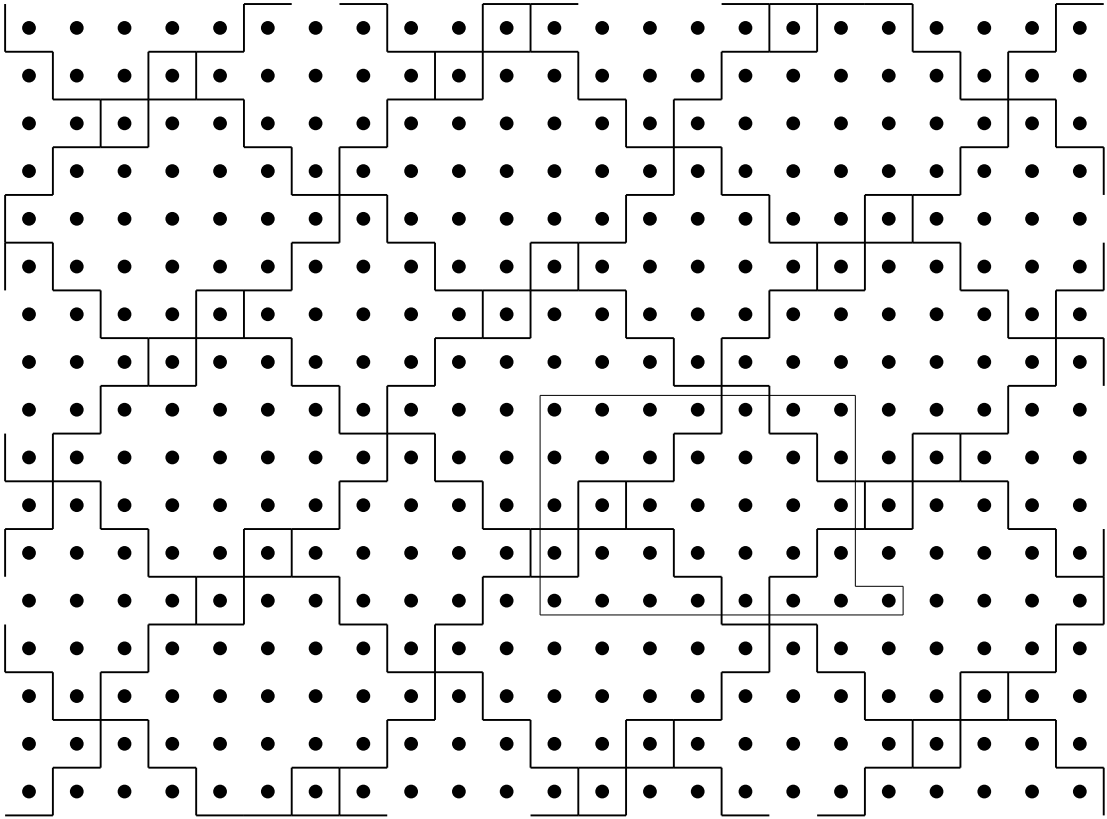
Lattice covering of $\Z^2$ using $\octa 2k$ and $\octa 2{k-1}$
(shown for $k = 3$).
\endcaption
\endinsert

The graph of diameter~$k$ resulting from this covering is the Cayley
graph of the cyclic group~$\Z_{4k^2}$ with unrestricted generators
$1$ and~$2k-1$ and order-2 generator~$2k^2$.  One can get another
Cayley graph of this size by using the lattice generated by
$(2k,0)$ and $(0,2k)$, but this graph will not be cyclic; it
comes from the group $\Z_{2k} \times \Z_{2k}$.

The outlined shape in Figure~7 (a $(2k+1) \times (2k-1)$ rectangle
with one extra point) is a fundamental region which
is convenient for an actual layout of nodes in
a network.  Note that a $2k \times 2k$ rectangle (or for that matter,
any rectangle with both sides greater than~1) will not work as
a fundamental region for this lattice.  The alternative lattice in the
previous paragraph does allow one to use a layout which is a
$2k \times 2k$ rectangle; in fact, this is just a toroidal mesh.
However, in either case one will have to make the extra connections
specified by the order-2 generator.

When one moves to $d=3$, it becomes harder to get optimal results,
so again the authors resorted to a computational search.
For $k=1$ the best graph is the Cayley graph of~$\Z_8$
with unrestricted generators $1,2,3$ and order-2 generator~$4$;
for $k=2$ the best is $\Z_{26}$ with unrestricted generators $1,2,8$
and order-2 generator~$13$.  For $3 \le k \le 10$ the optimal results,
like those for three generators alone, form a pattern of period~3,
as shown in Table~4.  (Again the best graphs are all cyclic.
This time the parameter~$a$ is defined to be the integer
nearest~$2k/3$.)

\midinsert
\centerline{
\vbox{\offinterlineskip
\halign{
&\vrule#&\hfil$\,\,\vphantom{|\Rtetra 3k|}{#}\,\,$\hfil\cr
\noalign{\hrule}
height2pt&\omit&&\omit&&\omit&&\omit&\cr
&k \bmod 3&&0&&1&&2&\cr
height2pt&\omit&&\omit&&\omit&&\omit&\cr
\noalign{\hrule}
height2pt&\omit&&\omit&&\omit&&\omit&\cr
&a&&2k/3&&(2k+1)/3&&(2k-1)/3&\cr
height2pt&\omit&&\omit&&\omit&&\omit&\cr
\noalign{\hrule}
height2pt&\omit&&\omit&&\omit&&\omit&\cr
&\gathered\text{Lattice}\\\text{generators}\endgathered
&&\gathered(2a,1,-1)\\(-1,2a,-1)\\(1,1,2a)\endgathered
&&\gathered(2a-1,-1,0)\\(1,2a,-1)\\(0,1,2a-1)\endgathered
&&\gathered(2a+1,-1,0)\\(1,2a,-1)\\(0,1,2a+1)\endgathered
&\cr
height2pt&\omit&&\omit&&\omit&&\omit&\cr
\noalign{\hrule}
height2pt&\omit&&\omit&&\omit&&\omit&\cr
&\gathered\text{Extra}\\\text{generator}\endgathered
&&(a,a+1,a-1)&&(a,a,a-1)&&(a+1,a,a)
&\cr
height2pt&\omit&&\omit&&\omit&&\omit&\cr
\noalign{\hrule}
height2pt&\omit&&\omit&&\omit&&\omit&\cr
&\gathered\text{Cyclic}\\\text{group size}\endgathered
&&\displaystyle\frac{64 k^3 + 108 k}{27}
&&\displaystyle\frac{64 k^3 + 60 k - 16}{27}
&&\displaystyle\frac{64 k^3 + 60 k + 16}{27}
&\cr
height2pt&\omit&&\omit&&\omit&&\omit&\cr
\noalign{\hrule}
height2pt&\omit&&\omit&&\omit&&\omit&\cr
&\gathered\text{Unrestricted}\\\text{generators}\endgathered
&&\gathered 1 \\ \smash{4a^3-2a^2+2a-1} \\ \smash{4a^3-2a^2+4a-1} \endgathered
&&\gathered 1 \\ 2a-1 \\ \smash{4a^2-2a+1} \endgathered
&&\gathered 1 \\ 2a+1 \\ \smash{4a^2+2a+1} \endgathered
&\cr
height2pt&\omit&&\omit&&\omit&&\omit&\cr
\noalign{\hrule}
height2pt&\omit&&\omit&&\omit&&\omit&\cr
&\gathered\text{Order-2}\\\text{generator}\endgathered
&&4a^3 + 3a
&&4a^3 - 4a^2 + 3a -1
&&4a^3 + 4a^2 + 3a +1
&\cr
height2pt&\omit&&\omit&&\omit&&\omit&\cr
\noalign{\hrule}
}}
}
\botcaption{Table 4}
Best undirected Cayley graphs of cyclic groups of diameter~$k$
using three generators plus one order-2 generator ($3 \le k \le 10$).
\endcaption
\endinsert

We can now apply the methods in the proof of Theorem~15
to show:

\proclaim{Theorem 22} For each $k \ge 3$, the cyclic undirected
Cayley graph using the group and generators specified in
Table~4 has diameter~$k$. \QNED

The authors again conjecture that these are actually the optimal
such Abelian Cayley graphs for all $k \ge 3$, not just for $3 \le k \le 10$.

For $d > 3$ one can get reasonably good results by letting $L$~be
approximately a cubic lattice, with $\vg$~near the center of one of
the cubes; this makes $L \cup (\vg + L)$ approximately a
body-centered cubic lattice.  Again, though, the efficiency
of this covering decreases exponentially with~$d$; one can do
much better using the results of Gritzmann~\cite{\Gritzmann}.

\smallpagebreak

One can also consider the possibility of using more than one generator of
order~2.  For instance, one could look at Cayley graphs of degree~$2d+2$
obtained by using $d$~unrestricted generators and two generators
of order~2.

However, this is not going to he helpful if one wants to construct
large undirected Cayley graphs of a given degree and diameter, at least
if the diameter is substantially larger than the degree.
For instance, suppose that the degree is fixed as~$2d+2$.  If one uses
$d$~unrestricted generators and two order-2 generators, then the
size of the resulting undirected Abelian Cayley graph of diameter~$k$ is at
most $4$~times the number of points in the $d$\snug-dimensional
shape~$\octa {d}k$.  (More precisely, by the methods used above
for one order-2 generator, one gets a limit of
$|\octa {d}k| + 2|\octa {d}{k-1}| + |\octa {d}{k-2}|$.)
This limit is~$O(k^{d})$, which is less than the $O(k^{d+1})$
one gets by using $d+1$~unrestricted generators.
The same argument shows that using more than two order-2 generators
cannot be useful for large~$k$; one gets larger graphs by replacing
two order-2 generators with one unrestricted generator.

If one is interested in the small-diameter case, though (especially when
the diameter is less than or equal to the degree), then order-2 generators
must be considered.  The most extreme version of this would be to make
{\sl all}\/ of the generators have order~2.  This then becomes precisely
the covering radius problem for binary linear codes; see the
surveys by Cohen et~al.~\cite{\CohKarMatSch,\CohLitLobMat} for
more on this problem.  The resulting graphs would be hypercubes
with additional diagonal connections to reduce the diameter.

\subhead Conclusions \endsubhead
We have shown that one can construct Cayley graphs of Abelian groups
which have substantially more vertices than traditional toroidal
meshes, but retain certain desirable features.  In particular, routing
on the twisted toroidal meshes is easily described in almost the same
manner as on toroidal meshes, and the twisted toroidal meshes host the
discrete non-periodic orthogonal grids used in numerical calculations
in exactly the same way that toroidal meshes do.  In addition, we have
shown how our $d$\snug-dimensional meshes can be constructed with
physical wire lengths that remain constant with increasing diameter
(and increasing number of vertices) just as the corresponding toroidal
meshes can.  We have given results which are provably optimal in 2
dimensions and probably optimal in 3 dimensions---the physically
interesting cases.

In the sequel to this paper, we will show that our methods can be
extended to cover certain types of nilpotent groups.  These groups
yield graphs with cardinalities which still increase polynomially with
diameter for a given degree, but with an exponent which is larger than
in the Abelian case.  One class of groups for which we obtain optimal
results includes the groups discussed in Draper and
Faber~\cite{\DraFab}.  In particular, we show that the particular
groups analyzed in that paper are not optimal for large diameters.


\Refs

\ref \no \AguFioGar \by F. Aguil\'o, M. A. Fiol, and C. Garcia
\paper Triple loop networks with small transmission delay
\jour Discrete Math. \vol 167/168 \yr 1997 \pages 3--16 \endref

\ref \no \AnnBau \by F. Annexstein and M. Baumslag \paper
On the diameter and bisector size of Cayley graphs \jour Math. Syst.
Theory \vol 26 \yr 1993 \pages 271--291 \endref

\ref \no \BerComHsu \by J. C. Bermond, F. Comellas, and D. F. Hsu
\paper Distributed loop computer networks: a survey
\jour J. Parallel Distrib. Comput. \vol 24 \yr 1994 \pages 2--10 \endref

\ref \no \BoeWan \by F. T. Boesch and J.-F. Wang
\paper Reliable circulant networks with minimum transmission delay
\jour IEEE Trans. Circuits Systems \vol CAS-32 \yr 1985 \pages 1286--1291
\endref

\ref \no \CheJia \by S. Chen and X.-D. Jia
\paper Undirected loop networks
\jour Networks \vol 23 \yr 1993 \pages 257--260 \endref

\ref \no \Chung \by F. R. K. Chung \paper Diameters and eigenvalues \jour
J. Amer. Math. Soc. \vol 2 \yr 1989 \pages 187--196 \endref

\ref \no \ChuFabMan \by F. R. K. Chung, V. Faber, and T. A.
Manteuffel \paper An upper bound on the diameter of a graph from
eigenvalues associated with its Laplacian \jour SIAM J. Discrete Math.
\vol 7 \yr 1994 \pages 443--457 \endref

\ref \no \CohKarMatSch \by G. D. Cohen, M. G. Karpovsky, H. F.
Mattson,~Jr., and J. R. Schatz \paper Covering radius --- survey and
recent results \jour IEEE  Trans. Inform. Theory \vol IT-31
\yr 1985 \pages 328--343 \endref

\ref \no \CohLitLobMat \by G. D. Cohen, S. N Litsyn, A. C. Lobstein, and
H. F. Mattson, Jr. \paper Covering radius 1985--1994
\jour Appl. Algebra Engrg. Comm. Comput. \vol 8 \yr 1997
\pages 173--239 \endref

\ref \no \DinHaf \by M. Dinneen and P. Hafner \paper New results for the
degree/diameter problem \jour Networks \vol 24 \yr 1994 \pages 359--367 \endref

\ref \no \DouJan \by R. Dougherty and H. Janwa \paper Covering radius
computations for binary cyclic codes \jour Math. Comp. \vol 57
\yr 1991 \pages 415--434 \endref

\ref \no \DraFab \by R. Draper and V. Faber \paper The diameter and mean
diameter of supertoroidal networks \paperinfo Supercomputing Research Center
Technical Report SRC-TR-90-004, 1990 \endref

\ref \no \Fary \by I. F\'ary \paper Sur la densit\'e des r\'eseaux de
domaines convexes \jour Bull. Soc. Math. France \vol 78 \yr 1950
\pages 152--161 \endref

\ref \no \FidForZit \by C. S. Fiduccia, R. W. Forcade, and J. S. Zito
\paper Geometry and diameter bounds of directed Cayley graphs of
Abelian groups
\jour SIAM J. Discrete Math. \vol 11 \yr 1998 \pages 157--167 \endref

\ref \no \ForLam \by R. Forcade and J. Lamoreaux \paper
Lattice-simplex coverings and the 84-shape
\jour SIAM J. Discrete Math. \vol 13 \yr 2000 \pages 194--201 \endref

\ref \no \Gritzmann \by P. Gritzmann \paper Lattice covering of
space with symmetric convex bodies \jour Mathematika \vol 32 \yr 1985
\pages 311--\allowbreak315 \endref

\ref \no \Hoylman \by D. Hoylman \paper The densest lattice packing
of tetrahedra \jour Bull. Amer. Math. Soc. \vol 76 \yr 170
\pages 135--137 \endref

\ref \no \Rogers \by C. Rogers \paper Lattice coverings of space
\jour Mathematika \vol 6 \yr 1959 \pages 33--39 \endref


\ref \no \Sabidussi \by G. Sabidussi \paper Vertex transitive graphs
\jour Monatsh. Math. \vol 68 \yr 1964 \pages 426--438 \endref

\ref \no \StaCow \by R. Stanton and D. Cowan \paper Note on a
``square'' functional equation \jour SIAM Review \vol 12 \yr 1970
\pages 277--279 \endref

\ref \no \Urakawa \by H. Urakawa \paper On the least positive eigenvalue
of the Laplacian for the compact quotient of a certain Riemannian
symmetric space \jour Nagoya Math. J. \vol 78 \yr 1980 \pages 137--152 \endref

\ref \no \WonCop \by C. K. Wong and D. Coppersmith
\paper A combinatorial problem related to multimodule memory organizations
\jour J. ACM \vol 21 \yr 1974 \pages 392--402 \endref

\ref \no \YebFioMorAle \by J. L. A. Yebra, M. A. Fiol, P. Morillo, and
I. Alegre
\paper The diameter of undirected graphs associated to plane tessellations
\jour Ars Combin. \vol 20-B \yr 1985 \pages 159--171 \endref

\endRefs
\enddocument